\input ./style/arxiv-general.cfg
\documentclass[aap,MSNbibl,seceqn,dvips]{arximspdf}
\makeatletter
   \@ifpackageloaded{graphicx}{}{\usepackage{graphicx}}
\makeatother
\usepackage[mathscr]{eucal}

\doi{10.1214/14-AAP1081}
\volume{25}
\issue{6}
\pubyear{2015}
\firstpage{3511}
\lastpage{3570}
\docsubty{FLA}

\makeatletter
\newcommand{\argmin}{\operatorname{arg}\operatorname{min}}
\newcommand{\ds}{\displaystyle}
\newcommand{\eqref}[1]{(\ref{#1})}
\newcommand{\lVert}{\Vert}
\newcommand{\rVert}{\Vert}
\newcommand{\rrvert}{\vert}
\newcommand{\llvert}{\vert}
\newtheorem{lemma}{Lemma}[section]
\newtheorem{theorem}{Theorem}[section]
\newtheorem{proposition}{Proposition}[section]
\newtheorem{corollary}{Corollary}[section]
\newproclaim{definition}{Definition}[section]
\newproclaim{assumption}{Assumption}[section]
\newproclaim{remark}{Remark}[section]
\makeatother

\begin{document}
\begin{frontmatter}

\title{Ergodic control of multi-class $M/M/N+M$ queues in~the
Halfin--Whitt regime}
\runtitle{Ergodic control in the Halfin--Whitt regime}

\begin{aug}
\author[A]{\fnms{Ari}~\snm{Arapostathis}\corref{}\thanksref{m1,T1}\ead[label=e1]{ari@ece.utexas.edu}},
\author[A]{\fnms{Anup}~\snm{Biswas}\ead[label=e2]{anupbiswas@utexas.edu}\thanksref{m1,T2}}
\and
\author[B]{\fnms{Guodong}~\snm{Pang}\thanksref{m2,T3}\ead
[label=e3]{gup3@psu.edu}}
\runauthor{A. Arapostathis, A. Biswas and G. Pang}
\affiliation{The University of Texas at Austin\thanksmark{m1}  and Pennsylvania State
University\thanksmark{m2}}
\address[A]{A. Arapostathis\\
A. Biswas\\
Department of Electrical  \\
\quad and Computer Engineering\\
The University of Texas at Austin\\
1616 Guadalupe St., UTA 7.508\\
Austin, Texas 78701\\
USA\\
\printead{e1}\\
\phantom{E-mail:\ }\printead*{e2}}
\address[B]{G. Pang\\
The Harold and Inge Marcus Department\\
\quad of Industrial and
Manufacturing  Engineering\\
College of Engineering\\
Pennsylvania State University\\
University Park, Pennsylvania 16802\\
USA\\
\printead{e3}}
\end{aug}
\thankstext{T1}{Supported in part by the Office of Naval
Research Grant N00014-14-1-0196.}
\thankstext{T2}{Supported in part by an award
from the Simons Foundation (\# 197982) to The University of Texas at Austin
and in part by
the Office of Naval Research through the Electric Ship Research and
Development Consortium.}
\thankstext{T3}{Supported in part by the Marcus Endowment Grant at the
Harold and Inge Marcus Department of Industrial and Manufacturing Engineering
at Penn State.}

%
\received{\smonth{4} \syear{2014}}
%
\revised{\smonth{11} \syear{2014}}

%
\begin{abstract}
We study a dynamic scheduling problem for a multi-class queueing network
with a large pool of statistically identical servers. The arrival processes
are Poisson, and service times and patience times are assumed to be
exponentially
distributed and class dependent. The optimization criterion is the
expected long
time average (ergodic) of a general (nonlinear) running cost function
of the
queue lengths. We consider this control problem in the Halfin--Whitt
(QED) regime, that is,
the number of servers $n$ and the total offered load~$\mathbf{r}$
scale like
$n\approx\mathbf{r}+\hat{\rho} \sqrt{\mathbf{r}}$ for some constant
$\hat\rho$. This problem was proposed in [\textit{Ann. Appl. Probab.} \textbf{14}
(2004) 1084--1134, Section~5.2].

The optimal solution of this control problem can be approximated by
that of the
corresponding ergodic diffusion control problem in the limit. We
introduce a broad
class of ergodic control problems for controlled diffusions, which
includes a large
class of queueing models in the diffusion approximation, and establish
a complete
characterization of optimality via the study of the associated HJB equation.
We also prove the asymptotic convergence of the values for the multi-class
queueing control problem to the value of the associated ergodic diffusion
control problem. The proof relies on an approximation method by
\emph{spatial truncation} for the ergodic control of diffusion
processes, where
the Markov policies follow a fixed priority policy outside a fixed
compact set.
\end{abstract}

%
\begin{keyword}[class=AMS]
\kwd[Primary ]{60K25}
\kwd[; secondary ]{68M20}
\kwd{90B22}
\kwd{90B36}
\end{keyword}
\begin{keyword}
\kwd{Multi-class Markovian queues}
\kwd{reneging/abandonment}
\kwd{Halfin--Whitt (QED) regime}
\kwd{diffusion scaling}
\kwd{long time-average control}
\kwd{ergodic control}
\kwd{stable Markov optimal control}
\kwd{spatial truncation}
\kwd{asymptotic optimality}
\end{keyword}
\end{frontmatter}

\setattribute{tocline}{skip}{\space}
\tableofcontents[alignleft,level=2]

\section{Introduction}

One of the classical problems in queueing theory is to schedule the
customers/jobs in a network
in an optimal way. These problems are known as the scheduling problems
which arise in a wide variety of applications, in particular, whenever
there are different customer classes
present in the network and competing for the same resources.
The optimal scheduling problem has a long history in the literature.
One of the appealing scheduling rules is the well-known~$c\mu$ rule.
This is a static priority policy in which
it is assumed that each class-$i$ customer has a marginal delay
cost $c_{i}$ and an average
service time ${1}/{\mu_{i}}$, and the classes are prioritized
in the
decreasing order of $c_{i}\mu_{i}$. This static priority rule has
proven asymptotically optimal in many settings
\cite{mieghem,atar-biswas,mandel-stolyar}.
In \cite{budhi-ghosh-liu}, a single-server Markov modulated
queueing network is considered and an \emph{averaged}
$c\mu$-rule is shown asymptotically optimal
for the discounted control problem.

An important aspect of queueing networks is abandonment/reneging, that is,
customers/jobs may choose to leave the system while being in the queue
before their service. Therefore, it is important to include customer
abandonment in
modeling queueing systems.
In \cite{atar-giat-shim,atar-giat-shim-2}, Atar et al.
considered a multi-class $M/M/N+M$ queueing network with
customer abandonment and proved that a modified priority policy,
referred to as $c\mu/\theta$
rule, is asymptotically optimal for the long run average cost in
the fluid scale. Dai and  Tezcan
\cite{tezcan-dai} showed the asymptotic optimality of a static
priority policy
on a finite time interval for a parallel
server model under the assumed conditions on the ordering
of the abandonment rates and running costs.
Although static priority policies are easy to implement, it may not be
optimal for
control problems of many multi-server queueing systems.
For the same multi-class $M/M/N+M$ queueing network,
discounted cost control problems are studied in
\cite{atar-mandel-rei,atar-2005,harrison-zeevi},
and asymptotically optimal controls for these problems are constructed
from the minimizer of a Hamilton--Jacobi--Bellman (HJB) equation
associated with
the controlled diffusions in the Halfin--Whitt regime.

In this article, we are interested in an ergodic control problem for a
multi-class
$M/M/N+M$ queueing network in the Halfin--Whitt regime.
The network consists of a single
pool of $n$ statistically identical servers and a buffer of infinite capacity.
There are $d$ customer classes and arrivals
of jobs/customers are $d$ independent Poisson processes
with parameters $\lambda^{n}_{i}$, $i=1,\ldots,d$.
The service rate for \mbox{class-$i$} customers
is $\mu_{i}^{n}$, $i=1,\ldots,d$.
Customers may renege from the queue if they have not started to receive
service before their patience times.
Class-$i$ customers renege from the queue at rates $\gamma_{i}^{n}>0$,
$i=1,\ldots,d$.
The scheduling policies are \emph{work-conserving},
that is, no server stays idle if any of the queues is nonempty.
We assume the system operates in the Halfin--Whitt
regime, where the arrival rates and the
number of servers are scaled appropriately in a manner that the traffic
intensity
of the system satisfies
\[
\sqrt{n} \Biggl(1- \sum_{i=1}^{d}
\frac{\lambda^{n}_{i}}{n\mu
_{i}^{n}} \Biggr) \mathop{\longrightarrow}_{n \to\infty} \hat\rho\in
\mathbb{R}.
\]
In this regime, the system operations achieve both high quality
(high server levels) and high efficiency (high servers' utilization),
and hence it is also referred to as the Quality-and-Efficiency-Driven
(QED) regime;
see, for example,
\cite{halfin-whitt,garnett-mandel-reiman,atar-mandel-rei,gamarnik-zeevi,gurvich}
on the many-server regimes.
We consider an ergodic cost function given by
\[
\limsup_{T\to\infty} \frac{1}{T} \mathbb{E} \biggl[\int
_{0}^{T}r \bigl(\hat{Q}^{n}(s) \bigr)\,
\mathrm{d} {s} \biggr],
\]
where the \emph{running cost} $r$ is a nonnegative, convex function
with polynomial
growth and $\hat{Q}^{n} = (\hat{Q}^{n}_{1},\ldots,\hat{Q}^{n}_{d})^{\mathsf{T}}$
is the diffusion-scaled queue length process.
It is worth mentioning that in addition to the running cost above
which is based on the queue-length, we can add an idle-server cost provided
that it has at most polynomial growth.
For such, a running cost structure the same analysis goes through.
The control is the allocation of servers to different classes of customers
at the service completion times.
The value function is defined to be the infimum of the above cost over
all \emph{admissible} controls (among all work-conserving scheduling policies).
In this article, we are interested in the existence and uniqueness of
asymptotically optimal
stable stationary Markov controls for the ergodic control problem, and
the asymptotic behavior of the value functions as $n$ tends to infinity.
In \cite{atar-mandel-rei}, Section~5.2, it is stated that analysis of this
type of problems is important for modeling call centers.

\subsection{Contributions and comparisons}
The usual methodology for studying these problems
is to consider the associated continuum model, which is the controlled diffusion
limit in a heavy-traffic regime, and to study the ergodic control problem
for the controlled diffusion.
Ergodic control problems governed by controlled
diffusions have been well studied in literature \cite{ari-bor-ghosh,borkar}
for models that fall in these two categories: (a) the running cost is
\emph{near-monotone}, which is defined by the requirement that
its value outside a compact set
exceeds the optimal average cost, thus penalizing unstable behavior
(see Assumption~3.4.2 in \cite{ari-bor-ghosh} for details), or
(b) the controlled diffusion is uniformly stable, that is, every
stationary Markov
control is stable and the collection of invariant probability
measures corresponding to
the stationary Markov controls is tight.
However, the ergodic control problem at hand does not fall under any of these
frameworks.
First, the running cost we consider here is not near-monotone
because the total queue length can be $0$
when the total number of customers in the system are $\mathscr{O}(n)$.
On the other hand, it is not at all clear that the controlled diffusion
is uniformly stable (unless one imposes nontrivial hypotheses on the
parameters), and this remains an open problem.
One of our main contributions in this article is that we solve the
ergodic control problem for a broad class of nondegenerate controlled
diffusions,
that in a certain way can be viewed as a mixture of the two categories
mentioned above.
As we show in Section~\ref{S-ergodic}, stability of the diffusion
under any optimal stationary Markov control
occurs due to certain interplay between the drift and the running cost.
The model studied in Section~\ref{S-ergodic} is far more general
than the queueing problem described, and thus it is
of separate interest for ergodic control.
We present a comprehensive study of this broad class of ergodic control
problems that includes existence of a solution to the ergodic HJB
equation, its stochastic representation and verification of optimality
(Theorem~\ref{T-HJB2}), uniqueness of the solution in a certain class
(Theorem~\ref{T-unique}), and convergence of the vanishing discount
method (Theorem~\ref{T-HJB3}).
These results extend the well-known results for near-monotone running
costs.
The assumptions in these theorems are verified for the multi-class queueing
model and the corresponding characterization of optimality is obtained
(Corollary~\ref{C-unique}), which includes growth estimates for the
solution of the HJB.

We also introduce a new approximation technique, \emph{spatial truncation},
for the controlled diffusion processes; see Section~\ref{S-truncation}.
It is shown that if we freeze the Markov controls to a fixed stable Markov
control outside a compact set, then we can still obtain
nearly optimal controls in this class of Markov controls for large
compact sets.
We should keep in mind that this property is not true in general.
This method can also be thought of as an approximation by a class of
controlled diffusions that are uniformly stable.

We remark that for a fixed control, the controlled diffusions
for the queueing model can be regarded as a special case of the piecewise
linear diffusions considered in \cite{dieker-gao}. It is shown in
\cite{dieker-gao} that these diffusions are stable under
\emph{constant} Markov controls.
The proof is via a suitable Lyapunov function.
We conjecture that uniform stability holds for the controlled diffusions
associated with the queueing model.
For the same multi-class Markovian model, Gamarnik and Stolyar
show that the stationary distributions of the queue lengths are
tight under any work-conserving policy \cite{gamarnik-stolyar}, Theorem~2.
We also wish to remark here that we allow~$\hat\rho$ to be negative, assuming abandonment rates are strictly positive,
while in \cite{gamarnik-stolyar}, $\hat\rho>0$ and abandonment rates
can be zero.

Another important contribution of this work is the convergence of the value
functions associated with the sequence of multi-class queueing models
to the
value of the ergodic control problem,
say $\varrho_{*}$, corresponding to the controlled diffusion model.
It is not obvious that one can have asymptotic optimality from the
existence of
optimal stable controls for the HJB equations of controlled diffusions.
This fact is relatively
straightforward when the cost under consideration is discounted.
In that situation, the tightness of paths on a finite time horizon is
sufficient to prove asymptotic optimality \cite{atar-mandel-rei}.
But we are in a situation where any finite time behavior of the stochastic
process plays no role in the cost.
In particular, we need to establish the convergence of the
controlled steady states.
Although uniform stability of stationary distributions for this multi-class
queueing model in the case where $\hat\rho>0$ and abandonment rates
can be zero
is established in \cite{gamarnik-stolyar}, it is not obvious that
the stochastic model considered here
has the property of uniform stability.
Therefore, we use a different method to establish the asymptotic
optimality.
First, we show that the value functions are asymptotically bounded
below by
$\varrho_{*}$.
To study the upper bound, we construct
a sequence of Markov scheduling policies
that are uniformly stable (see Lemma~\ref{lem-uni}).
The key idea used in establishing such
stability results is a spatial truncation technique,
under which the Markov policies follow a fixed priority policy outside
a given compact set.
We believe these techniques can also be used to study ergodic control problems
for other many-server queueing models.

The scheduling policies we consider in this paper allow preemption,
that is, a
customer in service can be interrupted for the server to serve a
customer of a
different class and her service will be resumed later.
In fact, the asymptotic optimality is shown within the class of the
work-conserving preemptive policies.
In \cite{atar-mandel-rei}, both preemptive and nonpreemptive policies
are studied,
where a nonpreemptive scheduling control policy is constructed from
the HJB
equation associated with preemptive policies and thus is shown to be
asymptotically optimal.
However, as far as we know,
the optimal nonpreemptive scheduling problem under
the ergodic cost remains open.

For a similar line of work in uncontrolled settings, we refer the
reader to
\cite{gamarnik-zeevi,gurvich}.
Admission control of the single class $M/M/N+M$ model with an ergodic cost
criterion in the Halfin--Whitt regime is studied in \cite{KW10}.
For controlled problems and for finite server models,
asymptotic optimality is obtained in \cite{budhi-ghosh-lee} in the
conventional heavy-traffic regime.
The main advantage in \cite{budhi-ghosh-lee} is the uniform
exponential stability
of the stochastic processes, which is obtained by using properties of the
Skorohod reflection map.
A recent work studying ergodic control of a multi-class single-server
queueing network is \cite{kim-ward}.

To summarize our main contributions in this paper:
\begin{itemize}
\item[--]
We introduce a new class of ergodic control problems and a framework to
solve them.
\item[--]
We establish an approximation technique by spatial truncation.
\item[--]
We provide, to the best of our knowledge,
the first treatment of ergodic control problems at the
diffusion scale for many server models.
\item[--]
We establish asymptotic optimality results.
\end{itemize}

\subsection{Organization}
In Section~\ref{S-notation}, we summarize the notation used in the paper.
In Section~\ref{S-main}, we introduce the multi-class many server
queueing model and
describe the Halfin--Whitt regime.
The ergodic control problem under the heavy-traffic setting
is introduced in Section~\ref{S-diffcontrol},
and the main results on asymptotic convergence are stated as
Theorems~\ref{T-lowerbound} and \ref{T-upperbound}.
Section~\ref{S-ergodic} introduces a class of controlled diffusions
and associated ergodic control problems, which contains the queueing models
in the diffusion scale.
The key structural assumptions are in Section~\ref{S-assumptions} and
these are verified for a generic class of queueing models in
Section~\ref{S3.3},
which are characterized by piecewise linear controlled diffusions.
Section~\ref{S3.4} concerns the existence of optimal controls under the
general hypotheses, while Section~\ref{S3.5} contains a comprehensive
study of the HJB equation.
Section~\ref{Proofs} is devoted to the proofs of the results
in Section~\ref{S3.5}.
The spatial truncation technique is introduced and studied in
Section~\ref{S-truncation}.
Finally, in Section~\ref{S-optimality} we prove the results
of asymptotic optimality.

\subsection{Notation}\label{S-notation}
The standard Euclidean norm in $\mathbb{R}^{d}$ is denoted by $\vert\cdot\vert$.
The set of nonnegative real numbers is denoted by $\mathbb{R}_{+}$,
$\mathbb{N}$ stands for the set of natural numbers, and $\mathbb{I}$ denotes
the indicator function.
By $\mathbb{Z}^{d}_{+}$ we denote the set of $d$-vectors of
nonnegative integers.
The closure, the boundary and the complement
of a set $A\subset\mathbb{R}^{d}$ are denoted
by $\overline{A}$, $\partial{A}$ and $A^{c}$, respectively.
The open ball of radius $R$ around $0$ is denoted by $B_{R}$.
Given two real numbers $a$ and $b$, the minimum (maximum) is denoted by
$a\wedge b$
($a\vee b$), respectively.
Define $a^{+}:=a\vee0$ and $a^{-}:=-(a\wedge0)$.
The integer part of a real number $a$ is denoted by $\lfloor a\rfloor$.
We use the notation
$e_{i}$, $i=1,\ldots,d$, to denote the vector with $i$th entry equal
to $1$
and all other entries equal to $0$.
We also let $e:=(1,\ldots,1)^{\mathsf{T}}$.
Given any two vectors $x,y\in\mathbb{R}^{d}$ the inner
product is denoted by $x\cdot y$.
By $\delta_{x}$ we denote the Dirac mass at $x$.
For any function $f\dvtx\mathbb{R}^{d}\to\mathbb{R}$
and domain $D\subset\mathbb{R}$ we define the oscillation of $f$ on
$D$ as follows:
\[
\mathop{\operatorname{osc}}_D(f) :=\sup \bigl\{f(x)-f(y) \dvtx x, y\in D
\bigr\}.
\]

For a nonnegative function $g\in\mathcal{C}(\mathbb{R}^{d})$,
we let $\mathscr{O}(g)$ denote the space of functions
$f\in\mathcal{C}(\mathbb{R}^{d})$ satisfying
$\sup_{x\in\mathbb{R}^{d}} \frac{\vert f(x)\vert}{1+g(x)}<\infty$.
This is a Banach space under the norm
\[
\lVert f\rVert_{g} := \sup_{x\in\mathbb{R}^{d}}
\frac{\vert f(x)\rrvert}{1+g(x)}.
\]
We also let $\mathfrak{o}(g)$ denote the subspace of $\mathscr{O}(g)$
consisting
of those functions $f$ satisfying
\[
\limsup_{\vert x\rrvert\to\infty} \frac{\vert f(x)\rrvert
}{1+g(x)} = 0.
\]
By a slight abuse of notation, we also denote by $\mathscr{O}(g)$ and
$\mathfrak{o}(g)$
a generic member of these spaces.
For two nonnegative functions $f$ and $g$, we use the notation $f\sim g$
to indicate that $f\in\mathscr{O}(g)$ and $g\in\mathscr{O}(f)$.

We denote by $L^{p}_{\mathrm{loc}}(\mathbb{R}^{d})$, $p\ge1$, the
set of
real-valued functions
that are locally $p$-integrable and by
$\mathscr{W}_{\mathrm{loc}}^{k,p}(\mathbb{R}^{d})$ the set of
functions in $L^{p}_{\mathrm{loc}}(\mathbb{R}^{d})$
whose $i$th weak derivatives, $i=1,\ldots,k$, are in
$L^{p}_{\mathrm{loc}}(\mathbb{R}^{d})$.
The set of all bounded continuous functions is denoted by $\mathcal{C}_{b}(\mathbb{R}^{d})$.
By $\mathcal{C}_{\mathrm{loc}}^{k,\alpha}(\mathbb{R}^{d})$ we
denote the set of functions that are
$k$-times continuously differentiable and whose $k$th derivatives are locally
H\"{o}lder continuous with exponent $\alpha$.
We define $\mathcal{C}^{k}_{b}(\mathbb{R}^{d})$, $k\ge0$, as the set
of functions
whose $i$th derivatives, $i=1,\ldots,k$, are continuous and bounded
in $\mathbb{R}^{d}$ and denote by
$\mathcal{C}^{k}_{c}(\mathbb{R}^{d})$ the subset of $\mathcal
{C}^{k}_{b}(\mathbb{R}^{d})$ with compact support.
For any path $X(\cdot)$,
we use the notation $\Delta X(t)$ to denote the jump at time $t$.
Given any Polish space $\mathcal{X}$, we denote by $\mathcal
{P}(\mathcal{X})$ the set of
probability measures on $\mathcal{X}$ and we endow $\mathcal
{P}(\mathcal{X})$ with the
Prokhorov metric.
For $\nu\in\mathcal{P}(\mathcal{X})$ and a Borel measurable map
$f\dvtx\mathcal{X}\to\mathbb{R}$,
we often use the abbreviated notation
\[
\nu(f) := \int_{\mathcal{X}} f \,\mathrm{d} {\nu}.
\]
The quadratic variation of a square integrable martingale is
denoted by $\langle\cdot,\cdot\rangle$ and the optional quadratic
variation by
$[ \cdot,\cdot]$. For presentation purposes we use the time variable
as the subscript for the diffusion processes.
Also $\kappa_{1},\kappa_{2},\ldots$ and $C_{1},C_{2},\ldots$
are used as generic constants whose values might vary from place to place.

\section{The controlled system in the Halfin--Whitt regime}\label{S-main}

\subsection{The multi-class Markovian many-server model}
\begin{figure}

\includegraphics{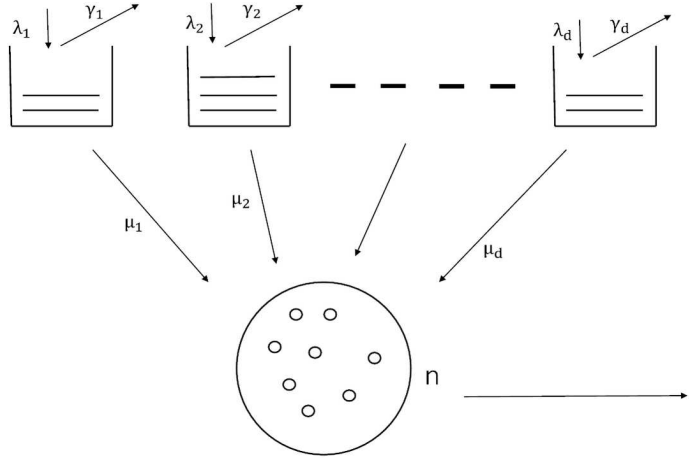}

\caption{A schematic model of the system.}\label{fig1}
\end{figure}

Let $(\Omega,\mathfrak{F},\mathbb{P})$ be a given complete probability
space and all the
stochastic variables introduced below are defined on it.
The expectation w.r.t. $\mathbb{P}$ is denoted by $\mathbb{E}$.
We consider a multi-class Markovian many-server queueing system which
consists of $d$ customer classes and $n$ parallel servers capable of
serving all customers (see Figure~\ref{fig1}).

The system buffer is assumed to have infinite capacity.
Customers of class $i\in\{1,\ldots,d\}$ arrive according to a Poisson process
with rate $\lambda^{n}_{i}>0$. Customers enter the queue of their respective
classes upon arrival if not being processed. Customers of each class are
served in the first-come-first-serve (FCFS) service discipline.
While waiting in queue, customers can abandon the system.
The service times and patience times of customers are class-dependent
and both are assumed to be exponentially distributed, that is, class $i$
customers are served at rate $\mu^{n}_{i}$ and renege at rate $\gamma^{n}_{i}$.
We assume that customer arrivals, service and abandonment of all classes
are mutually independent.

\subsubsection*{The Halfin--Whitt regime}
We study this queueing model in the Halfin--Whitt regime
[or the Quality-and-Efficiency-Driven (QED) regime].
Consider a sequence of such systems indexed by $n$, in which the
arrival rates
$\lambda_{i}^{n}$ and the number of servers $n$ both increase appropriately.
Let $\mathbf{r}^{n}_{i} :={\lambda_{i}^{n}}/{\mu^{n}_{i}}$
be the mean offered load of class $i$ customers.
The traffic intensity of the $n$th system is given by
$\rho^{n} = n^{-1} \sum_{i=1}^{d} \mathbf{r}^{n}_{i}$.
In the Halfin--Whitt regime, the parameters are assumed to satisfy the
following:
as $n\rightarrow\infty$,
%
\begin{eqnarray}
\frac{\lambda^{n}_{i}}{n} &\to & \lambda_{i} > 0,\qquad \mu_{i}^{n}
\to \mu_{i} > 0,\qquad \gamma_{i}^{n} \to
\gamma_{i} > 0,
\nonumber
\\
\label{HWpara}
\frac{\lambda^{n}_{i} - n \lambda_{i}}{\sqrt n} & \to & \hat{\lambda }_{i},\qquad
\sqrt{n} \bigl( \mu^{n}_{i} - \mu_{i} \bigr) \to
\hat{\mu}_{i},
\\
\frac{\mathbf{r}^{n}_{i}}{n} &\to & \rho_{i} := \frac{\lambda_{i}}{\mu_{i}} < 1,\qquad
\sum_{i=1}^{d} \rho_{i} = 1.
\nonumber
\end{eqnarray}
This implies that
\[
{\sqrt n} \bigl(1- \rho^{n} \bigr) \to\hat\rho:=\sum
_{i=1}^{d} \frac{\rho_{i} \hat{\mu}_{i} - \hat{\lambda}_{i}}{\mu_{i}} \in \mathbb{R}.
\]
The above scaling is common in multi-class multi-server models
\cite{harrison-zeevi,atar-mandel-rei}.
Note that we do not make any assumption on the sign of $\hat{\rho}$.

\subsubsection*{State descriptors}
Let $X_{i}^{n} = \{X_{i}^{n}(t) \dvtx t\ge0\}$ be the total number of class
$i$ customers in the system, $Q_{i}^{n} = \{Q_{i}^{n}(t) \dvtx t\ge0\}$
the number of class $i$ customers in the queue and
$Z_{i}^{n} = \{Z_{i}^{n}(t) \dvtx t\ge0\}$
the number of class $i$ customers in service.
The following basic relationships hold for these processes:
for each $t\ge0$ and $i =1,\ldots,d$,
%
\begin{eqnarray}
X_{i}^{n}(t)& =& Q_{i}^{n}(t) +
Z_{i}^{n}(t),
\nonumber
\\[-8pt]
\label{con1}
\\[-8pt]
\nonumber
Q_{i}^{n}(t) &\ge & 0,\qquad Z_{i}^{n}(t)
\ge0\quad \mbox{and}\quad e\cdot Z^{n}(t) \le n.
\end{eqnarray}
We can describe these processes using a collection
$ \{A_{i}^{n}, S_{i}^{n}, R_{i}^{n},i = 1,\ldots,d \}$
of independent rate-$1$ Poisson processes.
Define
\begin{eqnarray*}
\tilde{A}_{i}^{n}(t)& := & A_{i}^{n}
\bigl(\lambda^{n}_{i} t \bigr),
\\
\tilde{S}_{i}^{n}(t)& :=&  S_{i}^{n}
\biggl(\mu_{i}^{n} \int_{0}^{t}
Z_{i}^{n}(s) \,\mathrm{d} {s} \biggr),
\\
\tilde{R}_{i}^{n}(t)& :=&  R_{i}^{n}
\biggl(\gamma_{i}^{n} \int_{0}^{t}
Q_{i}^{n}(s) \,\mathrm {d} {s} \biggr).
\end{eqnarray*}
Then the dynamics take the form
%
\begin{equation}
\label{xp1}
X_{i}^{n}(t) = X_{i}^{n}(0)
+ \tilde{A}_{i}^{n}(t) - \tilde{S}_{i}^{n}(t)
- \tilde{R}_{i}^{n}(t),\qquad t\ge0, i = 1,\ldots,d.
\end{equation}

\subsubsection*{Scheduling control}
Following \cite{atar-mandel-rei,harrison-zeevi}, we only consider
work-conserving
policies that are nonanticipative and allow preemption.
When a server becomes free and there are no customers waiting in any queue,
the server stays idle, but if there are customers of multiple classes
waiting in the queue, the server has to make a decision on the customer
class to serve. Service preemption is allowed,
that is, service of a customer class can be interrupted at any time to
serve some
other class of
customers and the original service is resumed at a later time.
A scheduling control policy determines the processes $Z^{n}$, which
must satisfy
the constraints in~\eqref{con1} and the work-conserving constraint,
that is,
\[
e\cdot Z^{n}(t) = \bigl(e \cdot X^{n}(t) \bigr) \wedge n,\qquad t
\ge 0.
\]
Define the action set $\mathbb{A}^{n}(x)$ as
\[
\mathbb{A}^{n}(x) := \bigl\{ a \in\mathbb{Z}^{d}_{+}
\dvtx a \le x \mbox{ and } e\cdot a = (e \cdot x) \wedge n \bigr\}.
\]
Thus, we can write $Z^{n}(t) \in\mathbb{A}^{n}(X^{n}(t))$ for each
$t\ge0$.
We also assume that all controls are nonanticipative.
Define the $\sigma$-fields
\[
\mathcal{F}^{n}_{t}  :=  \sigma \bigl\{
X^{n}(0), \tilde{A}_{i}^{n}(t),
\tilde{S}_{i}^{n}(t), \tilde{R}_{i}^{n}(t)
\dvtx i =1,\ldots,d, 0 \le s \le t \bigr\} \vee\mathcal{N}
\]
and
\[
\mathcal{G}^{n}_{t}
 := \sigma \bigl\{\delta\tilde{A}_{i}^{n}(t,r), \delta
\tilde{S}_{i}^{n}(t,r), \delta\tilde{R}_{i}^{n}(t,r)
\dvtx i =1,\ldots,d, r\ge0 \bigr\},
\]
where
\begin{eqnarray*}
\delta\tilde{A}_{i}^{n}(t,r)& :=&  \tilde{A}_{i}^{n}(t+r)
-\tilde{A}_{i}^{n}(t),
\\
\delta\tilde{S}_{i}^{n}(t,r)& :=&  S_{i}^{n}
\biggl(\mu_{i}^{n}\int_{0}^{t}
Z_{i}^{n}(s) \,\mathrm{d} {s}+\mu_{i}^{n}r
\biggr) -\tilde{S}_{i}^{n}(t),
\\
\delta\tilde{R}_{i}^{n}(t,r)& :=&  R_{i}^{n}
\biggl(\gamma_{i}^{n}\int_{0}^{t}
Q_{i}^{n}(s) \,\mathrm{d} {s} +\gamma_{i}^{n}r
\biggr)-\tilde{R}_{i}^{n}(t),
\end{eqnarray*}
and $\mathcal{N}$ is the collection of all $\mathbb{P}$-null sets.
The filtration $\{\mathcal{F}^{n}_{t}, t\ge0\}$ represents the information
available up to time $t$ while $\mathcal{G}^{n}_{t}$ contains
the information about future increments of the processes.

We say that a working-conserving control policy is \emph{admissible} if:
\begin{longlist}[(iii)]
\item[(i)]
$Z^{n}(t)$ is adapted to $\mathcal{F}^{n}_{t}$,

\item[(ii)]
$\mathcal{F}^{n}_{t}$ is independent of $\mathcal{G}^{n}_{t}$ at each time
$t\ge0$,

\item[(iii)]
for each $i=1,\ldots,d$, and $t\ge0$, the process
$\delta\tilde{S}_{i}^{n}(t,\cdot)$
agrees in law with $S_{i}^{n}(\mu_{i}^{n} \cdot)$, and the process
$\delta\tilde{R}_{i}^{n}(t,\cdot)$
agrees in law with $R_{i}^{n}(\gamma_{i}^{n} \cdot)$.
\end{longlist}
We denote the set of all admissible control policies
$(Z^n, \mathcal{F}^{n}, \mathcal{G}^{n})$ by $\mathfrak{U}^{n}$.

\subsection{The ergodic control problem in the Halfin--Whitt regime}
\label{S-diffcontrol}

Define the diffusion-scaled processes
\[
\hat{X}^{n} = \bigl(\hat{X}^{n}_{1},\ldots,
\hat{X}^{n}_{d} \bigr)^{\mathsf{T}},\qquad \hat{Q}^{n} =
\bigl(\hat{Q}^{n}_{1},\ldots,\hat{Q}^{n}_{d}
\bigr)^{\mathsf{T}} \quad \mbox{and}\quad \hat{Z}^{n} = \bigl(
\hat{Z}^{n}_{1},\ldots,\hat{Z}^{n}_{d}
\bigr)^{\mathsf{T}},
\]
by
%
\begin{eqnarray}
\hat{X}^{n}_{i}(t) & := & \frac{1}{\sqrt{n}}
\bigl(X^{n}_{i}(t)- \rho_{i} n t \bigr),
\nonumber
\\[-2pt]
\label{dc1}
\hat{Q}^{n}_{i}(t) & :=& \frac{1}{\sqrt{n}}
Q^{n}_{i}(t),
\\[-2pt]
\hat{Z}^{n}_{i}(t) & :=& \frac{1}{\sqrt{n}}
\bigl(Z^{n}_{i}(t)- \rho_{i} n t \bigr)
\nonumber
\end{eqnarray}
for $t\ge0$.
By \eqref{xp1}, we can express $\hat{X}^{n}_{i}$ as
%
\begin{eqnarray}
\hat{X}^{n}_{i}(t) & =& \hat{X}^{n}_{i}(0)
+\ell_{i}^{n} t - \mu_{i}^{n} \int
_{0}^{t} \hat{Z}_{i}^{n}(s)
\,\mathrm{d} {s} - \gamma_{i}^{n} \int_{0}^{t}
\hat{Q}^{n}_{i}(s) \,\mathrm{d} {s}
\nonumber
\\[-9pt]
\label{dc2}
\\[-9pt]
&&{}+ \hat{M}_{A,i}^{n}(t) - \hat{M}_{S,i}^{n}(t)
- \hat{M}_{R,i}^{n}(t),
\nonumber
\end{eqnarray}
where $\ell^{n}= (\ell_{1}^{n},\ldots,\ell_{d}^{n})^{\mathsf{T}}$
is defined\vspace*{-1pt} as
\[
\ell_{i}^{n} :=\frac{1}{\sqrt{n}} \bigl(
\lambda_{i}^{n} - \mu_{i}^{n}
\rho_{i} n \bigr),
\]
and
%
\begin{eqnarray}
\hat{M}_{A,i}^{n}(t) & :=& \frac{1}{\sqrt{n}}
\bigl(A_{i}^{n} \bigl(\lambda ^{n}_{i} t
\bigr) - \lambda^{n}_{i} t \bigr),
\nonumber
\\[-2pt]
\label{dc3}
\hat{M}_{S,i}^{n}(t) & :=& \frac{1}{\sqrt{n}}
\biggl( S_{i}^{n} \biggl(\mu_{i}^{n}
\int_{0}^{t} Z_{i}^{n}(s)
\,\mathrm{d} {s} \biggr) - \mu_{i}^{n} \int
_{0}^{t} Z_{i}^{n}(s)
\,\mathrm{d} {s} \biggr),
\\[-2pt]
\hat{M}_{R,i}^{n}(t) & :=& \frac{1}{\sqrt{n}} \biggl(
R_{i}^{n} \biggl(\gamma_{i}^{n} \int
_{0}^{t} Q_{i}^{n}(s)
\,\mathrm{d} {s} \biggr) - \gamma_{i}^{n} \int
_{0}^{t} Q_{i}^{n}(s)
\,\mathrm{d} {s} \biggr)
\nonumber
\end{eqnarray}
are square integrable martingales w.r.t. the filtration $\{\mathcal
{F}^{n}_{t}\}$.

Note that
\[
\ell_{i}^{n} = \frac{1}{\sqrt{n}} \bigl(
\lambda_{i}^{n} - \lambda_{i} n \bigr) -
\rho_{i} \sqrt{n} \bigl(\mu_{i}^{n} -
\mu_{i} \bigr) \mathop{\longrightarrow}_{n\to\infty} \ell_{i}
:= \frac{(\hat{\lambda}_{i}- \rho_{i} \hat{\mu}_{i})}{\mu_{i}}.
\]

Define
\[
\mathcal{S} := \bigl\{u \in\mathbb{R}^{d}_{+} \dvtx e
\cdot u = 1 \bigr\}.
\]
For $Z^{n}\in\mathfrak{U}^{n}$ we define, for $t\ge0$ and for adapted
$\hat{U}^{n}(t)\in\mathcal{S}$,
%
\begin{eqnarray}
\hat{Q}^{n}(t) & :=& \bigl(e \cdot\hat{X}^{n}(t)
\bigr)^{+}\hat{U}^{n}(t),
\nonumber
\\[-8pt]
\label{con4}
\\[-8pt]
\nonumber
\hat{Z}^{n}(t) & :=& \hat{X}^{n}(t) - \bigl(e\cdot
\hat{X}^{n}(t) \bigr)^{+} \hat{U}^{n}(t).
\end{eqnarray}
If $\hat{Q}^{n}(t)=0$, we define $\hat{U}^{n}(t):=e_{d} = (0,\ldots,0,1)^{\mathsf{T}}$.
Thus, $\hat{U}^{n}_{i}$ represents the fraction of class-$i$ customers
in the
queue when the
total queue size is positive.
As we show later, it is convenient to view $\hat{U}^{n}(t)$ as the control.
Note that the controls are nonanticipative and preemption is allowed.

\subsubsection{The cost minimization problem}

We next introduce the running cost function for the control problem.
Let $r\dvtx\mathbb{R}^{d}_{+} \rightarrow\mathbb{R}_{+}$ be a
given function satisfying
%
\begin{equation}
\label{cost1}
c_{1} \vert x\rrvert^{m} \le r(x) \le
c_{2} \bigl(1+\llvert x\rrvert^{m} \bigr)\qquad \mbox{for some } m
\ge1,
\end{equation}
and some positive constants $c_{i}$, $i=1,2$.
We also assume that $r$ is locally Lipschitz.
This assumption includes linear and convex running cost functions. For example,
if we let $h_{i}$ be the holding cost rate for class $i$ customers,
then some
of the typical running cost functions are the following:
\[
r(x) = \sum_{i=1}^{d} h_{i}
x_{i}^{m}, \qquad m\ge1.
\]
These running cost functions evidently satisfy the condition in \eqref{cost1}.

Given the initial state $X^{n}(0)$ and a work-conserving scheduling
policy $Z^{n} \in\mathfrak{U}^{n}$, we
define the diffusion-scaled cost function as
%
\begin{equation}
\label{costd1}
J \bigl(\hat{X}^{n}(0),\hat{Z}^{n} \bigr) :=
\mathop{\limsup}_{T \to\infty} \frac{1}{T} \mathbb{E} \biggl[ \int
_{0}^{T} r \bigl(\hat{Q}^{n}(s) \bigr)
\,\mathrm{d} {s} \biggr],
\end{equation}
where the running cost function $r$ satisfies \eqref{cost1}.
Note that the running cost is defined using the scaled version of $Z^n$.
Then the associated cost minimization problem becomes
%
\begin{equation}
\label{E-Vn}
\hat{V}^{n} \bigl(\hat{X}^{n}(0) \bigr) :=
\inf_{Z^{n}\in\mathfrak{U}^{n}} J \bigl(\hat{X}^{n}(0),
\hat{Z}^{n} \bigr).
\end{equation}

We refer to $\hat{V}^{n}(\hat{X}^{n}(0))$ as the
\emph{diffusion-scaled value function} given
the initial state $\hat{X}^{n}(0)$ in the $n$th system.

From \eqref{con4}, it is easy to see that by
redefining $r$ as $r(x,u)=r((e\cdot x)^{+}u)$ we can rewrite the
control problem as
\[
\hat{V}^{n} \bigl(\hat{X}^{n}(0) \bigr) = \inf\tilde{J}
\bigl( \hat{X}^{n}(0), \hat{U}^{n} \bigr),
\]
where
%
\begin{equation}
\label{cost34} \tilde{J} \bigl(\hat{X}^{n}(0),\hat{U}^{n}
\bigr) := \limsup_{T\to\infty} \frac{1}{T} \mathbb{E} \biggl[
\int_{0}^{T} r \bigl(\hat{X}^{n}(s),
\hat{U}^{n}(s) \bigr) \,\mathrm{d} {s} \biggr],
\end{equation}
and the infimum is taken over all admissible pairs $(\hat{X}^{n},\hat{U}^{n})$
satisfying \eqref{con4}.

For simplicity, we assume that the initial condition $\hat{X}^{n}(0)$
is deterministic
and $\hat{X}^{n}(0)\to x$ as $n\to\infty$ for some $x\in\mathbb{R}^{d}$.

\subsubsection{The limiting controlled diffusion process}
As in \cite{atar-mandel-rei,harrison-zeevi}, one formally deduces that,
provided $\hat{X}^{n}(0) \to x$, there exists a limit $X$ for $\hat
{X}^{n}$ on
every finite time interval, and the limit process $X$ is a
$d$-dimensional diffusion process with independent components, that is,
%
\begin{equation}
\label{dc4}
\mathrm{d}X_{t} = b(X_{t},U_{t})
\,\mathrm{d} {t} + \Sigma\,\mathrm{d}W_{t},
\end{equation}
with initial condition $X_{0} = x$.
In \eqref{dc4}, the drift
$b(x,u)\dvtx\mathbb{R}^{d}\times\mathcal{S}\rightarrow\mathbb{R}^{d}$ takes the form
%
\begin{equation}
\label{dc5}
b(x,u) = \ell-R \bigl(x-(e\cdot x)^{+}u \bigr) -(e
\cdot x)^{+}\Gamma u,
\end{equation}
with
\begin{eqnarray*}
\ell & :=& (\ell_{1},\ldots,\ell_{d}) ^{\mathsf{T}},
\\
R & :=& \operatorname{diag}(\mu_{1},\ldots,\mu_{d}),
\\
\Gamma & :=& \operatorname{diag}(\gamma_{1},\ldots,
\gamma_{d}).
\end{eqnarray*}
The control $U_{t}$ lives in $\mathcal{S}$ and is nonanticipative,
$W(t)$ is a $d$-dimensional standard Wiener process independent of
the initial condition $X_{0}=x$, and the covariance matrix is given by
\[
\Sigma\Sigma^{\mathsf{T}}= \operatorname{diag}(2 \lambda _{1},
\ldots,2 \lambda_{d}).
\]
A formal derivation of the drift in \eqref{dc5} can be obtained from
\eqref{dc2} and \eqref{con4}.
A~detailed description of equation \eqref{dc4} and related results
are given in Section~\ref{S-ergodic}.
Let $\mathfrak{U}$ be the set of all admissible controls for the
diffusion model (for a definition see Section~\ref{S-ergodic}).

\subsubsection{The ergodic control problem in the diffusion scale}
Define $\tilde{r}\dvtx\mathbb{R}^{d}_{+} \times\mathbb{R}^{d}_{+}
\rightarrow\mathbb{R}_{+}$ by
\[
\tilde{r}(x,u) :=r \bigl({(e\cdot x)^{+}}u \bigr),
\]
where $r$ is the same function as in \eqref{costd1}.
In analogy with \eqref{cost34} we define the ergodic cost associated
with the controlled diffusion process $X$ and the running cost function
$\tilde{r}(x,u)$ as
\[
J(x,U) :=\limsup_{T\rightarrow\infty} \frac{1}{T}
\mathbb{E}^{U}_{x} \biggl[\int_{0}^{T}
\tilde{r}(X_{t},U_{t}) \,\mathrm{d} {t} \biggr],\qquad U \in
\mathfrak{U}.
\]
We consider the ergodic control problem
%
\begin{equation}
\label{dcp2}
\varrho_{*}(x) = \inf_{U \in\mathfrak{U}} J(x,U).
\end{equation}
We call $\varrho_{*}(x)$ the optimal
value at the initial state $x$ for the
controlled diffusion process $X$.
It is shown later that $\varrho_{*}(x)$ is independent of $x$.
A detailed treatment
and related results corresponding to the ergodic control problem are given
in Section~\ref{S-ergodic}.

We next state the main results of this section, the proof of which can be
found in Section~\ref{S-optimality}.

\begin{theorem}\label{T-lowerbound}
Let $\hat{X}^{n}(0)\to x\in\mathbb{R}^{d}$ as $n\to\infty$.
Also assume that \eqref{HWpara} and~\eqref{cost1} hold.
Then
\[
\liminf_{n\to\infty} \hat{V}^{n} \bigl(
\hat{X}^{n}(0) \bigr) \ge\varrho _{*}(x),
\]
where $\varrho_{*}(x)$ is given by \eqref{dcp2}.
\end{theorem}

\begin{theorem}\label{T-upperbound}
Suppose the assumptions of Theorem~\ref{T-lowerbound} hold.
In addition, assume that $r$ in \eqref{costd1} is convex.
Then
\[
\limsup_{n\to\infty} \hat{V}^{n} \bigl(
\hat{X}^{n}(0) \bigr) \le\varrho _{*}(x).
\]
\end{theorem}

Thus, we conclude that for any convex running cost function $r$,
Theorems~\ref{T-lowerbound} and~\ref{T-upperbound} establish the
asymptotic convergence of the ergodic control
problem for the queueing model.

\section{A broad class of ergodic control problems for diffusions}
\label{S-ergodic}

\subsection{The controlled diffusion model}

The dynamics are modeled by a
controlled diffusion process $X = \{X_{t}, t\ge0\}$
taking values in the $d$-dimensional Euclidean space $\mathbb{R}^{d}$, and
governed by the It\^{o} stochastic differential equation
%
\begin{equation}
\label{E-sde}
\mathrm{d} {X}_{t} = b(X_{t},U_{t})
\,\mathrm{d} {t} + \sigma(X_{t}) \,\mathrm{d} {W}_{t}.
\end{equation}
All random processes in \eqref{E-sde} live in a complete
probability space $(\Omega,\mathfrak{F},\mathbb{P})$.
The process $W$ is a $d$-dimensional standard Wiener process independent
of the initial condition $X_{0}$.
The control process $U$ takes values in a compact, metrizable set
$\mathbb{U}$, and
$U_{t}(\omega)$ is jointly measurable in
$(t,\omega)\in[0,\infty)\times\Omega$.
Moreover, it is \emph{nonanticipative}:
for $s < t$, $W_{t} - W_{s}$ is independent of
\[
\mathfrak{F}_{s} :=\mbox{the completion of } \sigma\{
X_{0},U_{r},W_{r}, r\le s\} \mbox{ relative to }(\mathfrak{F},\mathbb{P}).
\]
Such a process $U$ is called an \emph{admissible control},
and we let $\mathfrak{U}$ denote the set of all admissible controls.

We impose the following standard assumptions on the drift $b$
and the diffusion matrix $\sigma$
to guarantee existence and uniqueness of solutions to equation \eqref{E-sde}.
\begin{longlist}[(A3)]
\item[{(A1)}]
\emph{Local Lipschitz continuity}:
The functions
\[
b = \bigl[b^{1},\ldots,b^{d} \bigr]^{\mathsf{T}}\dvtx
\mathbb {R}^{d}\times\mathbb{U}\to\mathbb{R}^{d}\quad \mbox{and}\quad
\sigma= \bigl[\sigma^{ij} \bigr] \dvtx\mathbb{R}^{d}\to
\mathbb{R}^{d\times d}
\]
are locally Lipschitz in $x$ with a Lipschitz constant $C_{R}>0$
depending on
$R>0$.
In other words,
for all $x,y\in B_{R}$ and $u\in\mathbb{U}$,
\[
\bigl\llvert b(x,u) - b(y,u) \bigr\rrvert + \bigl\lVert\sigma(x) - \sigma (y)
\bigr\rVert \le C_{R} \llvert x-y\rrvert.
\]
We also assume that $b$ is continuous in $(x,u)$.

\item[{(A2)}]
\emph{Affine growth condition}:
$b$ and $\sigma$ satisfy a global growth condition of the form
\[
\bigl\llvert b(x,u) \bigr\rrvert^{2}+ \bigl\lVert\sigma(x) \bigr
\rVert^{2} \le C_{1} \bigl(1 + \llvert x
\rrvert^{2} \bigr)\qquad \forall(x,u)\in\mathbb {R}^{d}\times
\mathbb{U},
\]
where $\lVert\sigma\rVert^{2} :=
\operatorname{trace} (\sigma\sigma^{\mathsf{T}} )$.

\item[(A3)]
\emph{Local nondegeneracy}:
For each $R>0$, it holds that
\[
\sum_{i,j=1}^{d} a^{ij}(x)
\xi_{i}\xi_{j} \ge C^{-1}_{R}\llvert
\xi\rrvert^{2} \qquad\forall x\in B_{R},
\]
for all $\xi=(\xi_{1},\ldots,\xi_{d})^{\mathsf{T}}\in\mathbb{R}^{d}$,
where $a:=\sigma\sigma^{\mathsf{T}}$.
\end{longlist}

In integral form, \eqref{E-sde} is written as
%
\begin{equation}
\label{E3.2}
X_{t} = X_{0} + \int_{0}^{t}
b(X_{s},U_{s}) \,\mathrm{d} {s} + \int_{0}^{t}
\sigma(X_{s}) \,\mathrm{d} {W}_{s}.
\end{equation}
The third term on the right-hand side of \eqref{E3.2} is an It\^o
stochastic integral.
We say that a process $X=\{X_{t}(\omega)\}$ is a solution of \eqref{E-sde},
if it is $\mathfrak{F}_{t}$-adapted, continuous in~$t$, defined for all
$\omega\in\Omega$ and $t\in[0,\infty)$, and satisfies \eqref{E3.2} for
all $t\in[0,\infty)$ a.s.
It is well known that under (A1)--(A3), for any admissible control
there exists a unique solution of \eqref{E-sde}
\cite{ari-bor-ghosh}, Theorem~2.2.4.

We define the family of operators
$L^{u}\dvtx\mathcal{C}^{2}(\mathbb{R}^{d})\to\mathcal{C}(\mathbb
{R}^{d})$,
where $u\in\mathbb{U}$ plays the role of a parameter, by
%
\begin{equation}
\label{E3.3}
L^{u} f(x) :=\tfrac{1}{2} a^{ij}(x)
\partial_{ij} f(x) + b^{i}(x,u) \partial_{i} f(x),\qquad u\in\mathbb{U}.
\end{equation}
We refer to $L^{u}$ as the \emph{controlled extended generator} of
the diffusion.
In \eqref{E3.3} and elsewhere in this paper, we have adopted
the notation $\partial_{i}:=\frac{\partial }{\partial{x}_{i}}$ and
$\partial_{ij}:=\frac{\partial^{2} }{\partial{x}_{i}\,\partial{x}_{j}}$.
We also use the standard summation rule that
repeated subscripts and superscripts are summed from $1$ through $d$.
In other words, the right-hand side of \eqref{E3.3} stands for
\[
\frac{1}{2} \sum_{i,j=1}^{d}
a^{ij}(x) \frac{\partial^{2}f}{\partial{x}_{i}\,\partial{x}_{j}}(x) +\sum_{i=1}^{d}
b^{i}(x,u) \frac{\partial f }{\partial{x}_{i}}(x).
\]

Of fundamental importance in the study of functionals of $X$ is
It\^o's formula.
For $f\in\mathcal{C}^{2}(\mathbb{R}^{d})$ and with $L^{u}$ as
defined in \eqref{E3.3},
it holds that
%
\begin{equation}
\label{E3.4} f(X_{t}) = f(X_{0}) + \int
_{0}^{t}L^{U_{s}} f(X_{s})
\,\mathrm{d} {s} + M_{t},\qquad \mbox{a.s.},
\end{equation}
where
\[
M_{t} :=\int_{0}^{t} \bigl\langle
\nabla f(X_{s}), \sigma(X_{s}) \,\mathrm{d}
{W}_{s} \bigr\rangle
\]
is a local martingale.
Krylov's extension of It\^o's formula \cite{krylov}, page 122,
extends~\eqref{E3.4} to functions $f$ in the local
Sobolev space $\mathscr{W}_{\mathrm{loc}}^{2,p}(\mathbb{R}^{d})$,
$p\ge d$.

Recall that a control is called \emph{Markov} if
$U_{t} = v(t,X_{t})$ for a measurable map $v\dvtx\mathbb
{R}_{+}\times\mathbb{R}^{d}\to\mathbb{U}$,
and it is called \emph{stationary Markov} if $v$ does not depend on
$t$, that is, $v\dvtx\mathbb{R}^{d}\to\mathbb{U}$.
Correspondingly, \eqref{E-sde}
is said to have a \emph{strong solution}
if given a Wiener process $(W_{t},\mathfrak{F}_{t})$
on a complete probability space $(\Omega,\mathfrak{F},\mathbb{P})$, there
exists a process $X$ on $(\Omega,\mathfrak{F},\mathbb{P})$, with
$X_{0}=x_{0}\in\mathbb{R}^{d}$,
which is continuous,
$\mathfrak{F}_{t}$-adapted, and satisfies \eqref{E3.2} for all $t$ a.s.
A strong solution is called \emph{unique},
if any two such solutions $X$ and $X'$ agree
$\mathbb{P}$-a.s., when viewed as elements of $\mathcal{C}
([0,\infty),\mathbb{R}^{d} )$.
It is well known that under assumptions (A1)--(A3),
for any Markov control $v$,
\eqref{E-sde} has a unique strong solution \cite{Gyongy-96}.

Let $\mathfrak{U}_{\mathrm{SM}}$ denote the set of stationary Markov controls.
Under $v\in\mathfrak{U}_{\mathrm{SM}}$, the process $X$ is strong Markov,
and we denote its transition function by $P^{t}_{v}(x,\cdot)$.
It also follows from the work of \cite{Bogachev-01,Stannat-99} that under
$v\in\mathfrak{U}_{\mathrm{SM}}$, the transition probabilities of $X$
have densities which are locally H\"older continuous.
Thus, $L^{v}$ defined by
\[
L^{v} f(x) :=\tfrac{1}{2} a^{ij}(x)
\partial_{ij} f(x) + b^{i} \bigl(x,v(x) \bigr)
\partial_{i} f(x),\qquad v\in\mathfrak {U}_{\mathrm{SM}},
\]
for $f\in\mathcal{C}^{2}(\mathbb{R}^{d})$,
is the generator of a strongly-continuous
semi-group on $\mathcal{C}_{b}(\mathbb{R}^{d})$, which is strong Feller.
We let $\mathbb{P}_{x}^{v}$ denote the probability measure and
$\mathbb{E}^{v}_{x}$ the expectation operator on the canonical space
of the
process under the control $v\in\mathfrak{U}_{\mathrm{SM}}$,
conditioned on the
process $X$ starting from $x\in\mathbb{R}^{d}$ at $t=0$.

We need the following definition.

\begin{definition}
A function $h\dvtx\mathbb{R}^{d}\times\mathbb{U}\to\mathbb{R}$
is called \emph{inf-compact}
on a set $A\subset\mathbb{R}^{d}$ if the set
$\bar{A}\cap \{x \dvtx\min_{u\in\mathbb{U}} h(x,u)\le\beta
 \}$
is compact (or empty) in $\mathbb{R}^{d}$ for all $\beta\in\mathbb{R}$.
When this property
holds for $A\equiv\mathbb{R}^{d}$, then we simply say that $h$ is inf-compact.
\end{definition}

Recall that control $v\in\mathfrak{U}_{\mathrm{SM}}$ is called \emph{stable}
if the associated diffusion is positive recurrent.
We denote the set of such controls by $\mathfrak{U}_{\mathrm{SSM}}$,
and let $\mu_{v}$ denote the unique invariant probability
measure on $\mathbb{R}^{d}$ for the diffusion under the control $v\in
\mathfrak{U}_{\mathrm{SSM}}$.
We also let $\mathcal{M}:=\{\mu_{v} \dvtx v\in\mathfrak
{U}_{\mathrm{SSM}}\}$.
Recall that $v\in\mathfrak{U}_{\mathrm{SSM}}$ if and only if there
exists an inf-compact function
$\mathcal{V}\in\mathcal{C}^{2}(\mathbb{R}^{d})$, a bounded domain
$D\subset\mathbb{R}^{d}$, and
a constant $\varepsilon>0$ satisfying
\[
L^{v}\mathcal{V}(x) \le-\varepsilon \qquad\forall x\in D^{c}.
\]
We denote by $\tau(A)$ the \emph{first exit time} of a process
$\{X_{t}, t\in\mathbb{R}_{+}\}$ from a set $A\subset\mathbb
{R}^{d}$, defined by
\[
\tau(A) :=\inf\{t>0 \dvtx X_{t}\notin A\}.
\]
The open ball of radius $R$ in $\mathbb{R}^{d}$, centered at the origin,
is denoted by $B_{R}$, and we let $\tau_{R} :=\tau(B_{R})$,
and ${\breve\tau}_{R}:=\tau(B^{c}_{R})$.

We assume that the running cost function $r(x,u)$
is nonnegative, continuous and locally Lipschitz
in its first argument uniformly in $u\in\mathbb{U}$.
Without loss of generality, we let $\kappa_{R}$
be a Lipschitz constant of $r( \cdot,u)$ over $B_{R}$.
In summary, we assume that
\begin{longlist}[(A4)]
\item[{(A4)}]
$r\dvtx\mathbb{R}^{d}\times\mathbb{U}\to\mathbb{R}_{+}$ is
continuous and satisfies,
for some constant $C_{R}>0$
\[
\bigl\llvert r(x,u)-r(y,u) \bigr\rrvert \le C_{R} \llvert x-y\rrvert
\qquad\forall x,y\in B_{R}, \forall u\in\mathbb{U},
\]
and all $R>0$.
\end{longlist}

In general, $\mathbb{U}$ may not be a convex set.
It is therefore often useful to enlarge the control set to
$\mathcal{P}(\mathbb{U})$.
For any $v(\mathrm{d}{u})\in\mathcal{P}(\mathbb{U})$ we can
redefine the drift and the
running cost as
%
\begin{equation}
\label{relax}
\bar{b}(x,v) :=\int_\mathbb{U}b(x,u)v(
\mathrm{d} {u})\quad \mbox{and}\quad \bar{r}(x,v) :=\int_\mathbb{U}r(x,u)v(
\mathrm{d} {u}).
\end{equation}
It is easy to see that the drift and running cost
defined in \eqref{relax} satisfy all the
aforementioned conditions {(A1)}--{(A4)}.
In what follows, we assume that
all the controls take values in $\mathcal{P}(\mathbb{U})$.
These controls are generally referred to as \emph{relaxed} controls.
We endow the set of relaxed stationary Markov controls with the following
topology: $v_{n}\to v$ in $\mathfrak{U}_{\mathrm{SM}}$ if and only if
\[
\int_{\mathbb{R}^{d}}f(x)\int_\mathbb{U}g(x,u)v_{n}(\mathrm{d} {u}\vert x) \,\mathrm{d} {x} \mathop{\longrightarrow}_{n\to\infty} \int
_{\mathbb{R}^{d}}f(x)\int_\mathbb{U}g(x,u)v(\mathrm{d}
{u}\vert x) \,\mathrm{d} {x}
\]
for all $f\in L^{1}(\mathbb{R}^{d})\cap L^{2}(\mathbb{R}^{d})$
and $g\in\mathcal{C}_{b}(\mathbb{R}^{d}\times\mathbb{U})$.
Then $\mathfrak{U}_{\mathrm{SM}}$ is a compact metric space under
this topology
\cite{ari-bor-ghosh}, Section~2.4.
We refer to this topology as the \emph{topology of Markov controls}.
A control is said to be \emph{precise} if it takes value in $\mathbb{U}$.
It is easy to see that any precise control $U_{t}$ can also be understood
as a relaxed control by $U_{t}(\mathrm{d}{u})=\delta_{U_{t}}$.
Abusing the notation, we denote the drift and running cost by $b$ and~$r$,
respectively, and the action of a relaxed control
on them is understood as in~\eqref{relax}.

\subsection{Structural assumptions}
\label{S-assumptions}

Assumptions~\ref{Ass-1} and \ref{Ass-2}, described below,
are in effect throughout the analysis, unless otherwise stated.

\begin{assumption}\label{Ass-1}
For some open set $\mathcal{K}\subset\mathbb{R}^{d}$,
the following hold:
\begin{longlist}[(ii)]
\item[(i)]
The running cost $r$ is inf-compact on $\mathcal{K}$.

\item[(ii)]
There exist inf-compact functions $\mathcal{V}\in\mathcal
{C}^{2}(\mathbb{R}^{d})$
and $h \in\mathcal{C}(\mathbb{R}^{d}\times\mathbb{U})$, such that
%
\begin{eqnarray}
L^{u}\mathcal{V}(x) & \le & 1 - h(x,u) \qquad\forall(x,u)\in
\mathcal{K}^{c}\times\mathbb{U},
\nonumber
\\[-8pt]
\label{E-KA1}
\\[-8pt]
\nonumber
L^{u}\mathcal{V}(x) & \le & 1 + r(x,u) \qquad\forall(x,u)\in
\mathcal{K}\times \mathbb{U}.
\end{eqnarray}
\end{longlist}
Without loss of generality, we assume that $\mathcal{V}$ and $h$ are
nonnegative.
\end{assumption}

\begin{remark}
In the statement of Assumption~\ref{Ass-1}, we refrain from using any constants
in the interest of notational economy.
There is no loss of generality in doing so, since the functions
$\mathcal{V}$ and $h$ can always be scaled to eliminate
unnecessary constants.
\end{remark}

The next assumption is not a structural one, but rather
the necessary requirement that the value of the ergodic control
problem is finite.
Otherwise, the problem is vacuous.
For $U\in\mathfrak{U}$, define
%
\begin{equation}
\label{E-rU}
\varrho_{U}(x) :=\limsup_{T\to\infty}
\frac{1}{T} \mathbb{E}^{U}_{x} \biggl[\int
_{0}^{T} r(X_{s},U_{s})
\,\mathrm{d} {s} \biggr].
\end{equation}

\begin{assumption}\label{Ass-2}
There exists $U\in\mathfrak{U}$ such that $\varrho_{U}(x)<\infty$
for some $x\in\mathbb{R}^{d}$.
\end{assumption}

Assumption~\ref{Ass-2} alone does not imply that $\varrho_{v}<\infty$
for some $v\in\mathfrak{U}_{\mathrm{SSM}}$.
However, when combined with Assumption~\ref{Ass-1}, this is
the case as the following lemma asserts.

\begin{lemma}\label{L-basic}
Let Assumptions~\ref{Ass-1} and \ref{Ass-2} hold.
Then there exists $u_{0}\in\mathfrak{U}_{\mathrm{SSM}}$ such that
$\varrho_{u_{0}}<\infty$.
Moreover, there exists a nonnegative inf-compact function
$\mathcal{V}_{0}\in\mathcal{C}^{2}(\mathbb{R}^{d})$, and a positive
constant $\eta$ such that
%
\begin{equation}
\label{E-lyap}
L^{u_{0}}\mathcal{V}_{0}(x) \le\eta- r
\bigl(x,u_{0}(x) \bigr)\qquad \forall x\in \mathbb{R}^{d}.
\end{equation}
Conversely, if \eqref{E-lyap} holds, then Assumption~\ref{Ass-2} holds.
\end{lemma}

\begin{pf}
The first part of the result follows from Theorem~\ref{T-bst}(e)
and \eqref{E-KA3} whereas the converse part follows from Lemma~\ref{L3.2}.
These proofs are stated later in the paper.
\end{pf}

\begin{remark}
There is no loss of generality in using only the constant
$\eta$ in Assumption~\ref{Ass-2}, since
$\mathcal{V}_{0}$ can always be scaled to achieve this.

We also observe that for $\mathcal{K}=\mathbb{R}^{d}$ the problem
reduces to an
ergodic control problem with near-monotone
cost, and for $\mathcal{K}=\varnothing$ we obtain an ergodic control
problem under a
uniformly stable controlled diffusion.
\end{remark}

\subsection{Piecewise linear controlled diffusions}
\label{S3.3}

The controlled diffusion process in \eqref{dc4} belongs to a large class
of controlled diffusion processes, called piecewise linear controlled
diffusions \cite{dieker-gao}.
We describe this class of controlled diffusions
and show that it satisfies the assumptions in Section~\ref{S-assumptions}.

\begin{definition}
A square matrix $R$ is said to be an $M$-matrix if it can be written as $R=sI-N$
for some $s>0$ and nonnegative matrix $N$ with property that $\rho
(N)\le s$, where
$\rho(N)$ denotes the spectral radius of $N$.
\end{definition}

Let $\Gamma=[\gamma^{ij}]$ be a given matrix whose diagonal elements
are positive,
$\gamma^{id}=0$ for $i=1,\ldots,d-1$, and the remaining elements are in
$\mathbb{R}$.
(Note that for the queueing model, $\Gamma$ is a positive diagonal matrix.
Our results below hold for the more general $\Gamma$.)
Let $\ell\in\mathbb{R}^{d}$ and $R$ be a nonsingular $M$-matrix.
Define
%
\begin{equation}
\label{q-drift} b(x,u) :=\ell-R \bigl(x-{(e\cdot x)^{+}}u \bigr) -{(e
\cdot x)^{+}}\Gamma u,
\end{equation}
with
$u\in\mathcal{S}:= \{u \in\mathbb{R}^{d}_{+} \dvtx e\cdot
u=1 \}$.
Assume that
\[
e^{\mathsf{T}} R \ge 0^{\mathsf{T}}.
\]
We consider the following controlled diffusion in $\mathbb{R}^{d}$:
%
\begin{equation}
\label{eg-sde1}
\mathrm{d} {X}_{t} = b(X_{t},U_{t})
\,\mathrm{d} {t}+\Sigma\,\mathrm{d} {W}_{t},
\end{equation}
where $\Sigma$ is a constant matrix such that $\Sigma\Sigma^{\mathsf
{T}}$ is invertible.
It is easy to
see that~\eqref{eg-sde1} satisfies conditions (A1)--(A3).

Analysis of these types of diffusion approximations is an established tradition
in queueing
systems. It is often easy to deal with the limiting object and it also helps
to obtain information on the behavior of the actual queueing model.

We next introduce the running cost function.
Let $r\dvtx\mathbb{R}^{d}\times\mathcal{S}\to[0,\infty)$
be locally Lipschitz with polynomial growth and
%
\begin{equation}
\label{eg-cost}
c_{1} \bigl[(e\cdot x)^{+}
\bigr]^{m} \le r(x,u) \le c_{2} \bigl(1+ \bigl[(e\cdot
x)^{+} \bigr]^{m} \bigr),
\end{equation}
for some $m\ge1$ and positive constants $c_{1}$ and $c_{2}$ that
do not depend on $u$.
Some typical examples of such running costs are
\[
r(x,u) = \bigl[(e\cdot x)^{+} \bigr]^{m}\sum
_{i=1}^{d} h_{i} u_{i}^{m}\qquad \mbox{with } m\ge1,
\]
for some positive vector
$(h_{1},\ldots,h_{d})^{\mathsf{T}}$.

\begin{remark}
The controlled dynamics in \eqref{q-drift} and running cost
in \eqref{eg-cost} are clearly more general than the model described
in Section~\ref{S-diffcontrol}.
In \eqref{eg-sde1}, $X$
denotes the diffusion approximation for the number customers in the
system in the
Halfin--Whitt regime and its $i$th
component $X^{i}$ denotes the diffusion approximation of
the number of class $i$ customers.
Therefore, $(e\cdot X)^{+}$ denotes the total number of~customers in
the queue.
For $R$ and $\Gamma$ diagonal as in \eqref{dc5},
the diagonal entries of $R$ and $\Gamma$ denote the service and
abandonment rates, respectively, of the customer classes.
The $i$th coordinate of $U$ denotes the fraction of class-$i$ customers
waiting in the queue.
Therefore, the vector-valued process $X_{t}-(e\cdot X_{t})^{+}U_{t}$ denotes
the diffusion approximation of the numbers of customers in service
from different customer classes.
\end{remark}

\begin{proposition}\label{eg-prop}
Let $b$ and $r$ be given by \eqref{q-drift}
and \eqref{eg-cost}, respectively.
Then~\eqref{eg-sde1} satisfies Assumptions~\ref{Ass-1} and~\ref{Ass-2},
with $h(x)=c_{0} \llvert x\rrvert^{m}$ and
%
\begin{equation}
\label{eg-cK}
\mathcal{K} := \bigl\{x \dvtx \delta\llvert x\rrvert < (e\cdot
x)^{+} \bigr\}
\end{equation}
for appropriate positive constants
$c_{0}$ and $\delta$.
\end{proposition}

\begin{pf}
We recall\vspace*{1pt} that if $R$ is a nonsingular
$M$-matrix, then there exists a positive definite matrix $Q$
such that $QR+R^{\mathsf{T}}Q$ is strictly positive definite
\cite{dieker-gao}.
Therefore, for some positive constant $\kappa_{0}$ it holds that
\[
\kappa_{0} \llvert y\rrvert^{2} \le y^{\mathsf{T}}
\bigl[QR+R^{\mathsf{T}}Q \bigr]y \le\kappa_{0}^{-1}\llvert y
\rrvert^{2} \qquad\forall y \in\mathbb{R}^{d}.
\]
The set $\mathcal{K}$ in \eqref{eg-cK}, where $\delta>0$ is chosen later,
is an open convex cone, and the
running cost function $r$ is inf-compact on $\mathcal{K}$.
Let $\mathcal{V}$ be a nonnegative
function in $\mathcal{C}^{2}(\mathbb{R}^{d})$
such that $\mathcal{V}(x)= [x^{\mathsf{T}}Qx]^{{m}/{2}}$ for $\llvert x\rrvert\ge1$,
where the constant $m$ is as in~\eqref{eg-cost}.

Let $\llvert x\rrvert\ge1$ and $u\in\mathcal{S}$.\vspace*{-2pt} Then
\begin{eqnarray*}
\nabla\mathcal{V}(x)\cdot b(x,u)& = &\ell\cdot\nabla\mathcal{V}(x) -
\frac{m [x^{\mathsf{T}}Qx ]^{{m}/{2}-1}}{2}x^{\mathsf{T}} \bigl[QR+R^{\mathsf{T}}Q \bigr]x
\\[-2pt]
&&{}+m \bigl[x^{\mathsf{T}}Qx \bigr]^{{m}/{2}-1}Qx \cdot(R- \Gamma) (e
\cdot x)^{+}u
\\[-2pt]
& \le & \ell\cdot\nabla\mathcal{V}(x)-m \bigl[x^{\mathsf{T}}Qx
\bigr]^{{m}/{2}-1} \biggl(\frac{\kappa_{0}}{2}\llvert x\rrvert^{2}- C
\llvert x\rrvert {(e\cdot x)^{+}} \biggr)
\end{eqnarray*}
for some positive constant $C$.
If we choose $\delta=\frac{\kappa_{0}}{4C}$, then on
$\mathcal{K}^{c}\cap\{\llvert x\rrvert\ge1\}$ we have the estimate
%
\begin{equation}
\label{eg-1}
\nabla\mathcal{V}(x)\cdot b(x,u) \le\ell\cdot\nabla\mathcal{V}(x)
-\frac{m\kappa_{0}}{4} \bigl[x^{\mathsf{T}}Qx \bigr]^{{m}/{2}-1}\llvert x
\rrvert^{2}.
\end{equation}
Note that $\ell\cdot\mathcal{V}$ is globally bounded for $m=1$.
For any $m\in(1,\infty)$, it follows by~\eqref{eg-1} that
\begin{eqnarray}
\quad\nabla\mathcal{V}(x)\cdot b(x,u) & \le & m \bigl( \ell^{\mathsf{T}}Qx \bigr)
\bigl[x^{\mathsf{T}}Qx \bigr]^{{m}/{2}-1} - \frac{m\kappa_{0}}{4}
\bigl[x^{\mathsf{T}}Qx \bigr]^{{m}/{2}-1}\llvert x \rrvert^{2}
\nonumber
\\[-8pt]
\label{eg-2}
\\[-8pt]
\nonumber
& \le &\frac{m \llvert\ell^{\mathsf{T}}Q\rrvert
 (\overline\lambda(Q) )^{{m}/{2}}}{\underline
\lambda(Q)} \llvert x\rrvert^{m-1} -
\frac{m\kappa_{0}
 (\underline\lambda(Q) )^{{m}/{2}}}{4\overline
\lambda(Q)} \llvert x\rrvert^{m}
\nonumber
\end{eqnarray}
for $x\in\mathcal{K}^{c}\cap\{\llvert x\rrvert\ge1\}$,
where $\underline\lambda(Q)$ and $\overline\lambda(Q)$ are the
smallest and
largest eigenvalues of $Q$, respectively.
We use Young's inequality
\[
|ab| \le\frac{|a|^{m}}{m}+\frac{m-1}{m}|b|^{{m}/({m-1})},\qquad  a, b \ge0,
\]
in \eqref{eg-2} to obtain the bound
%
\begin{equation}
\label{E3.14nu}
\nabla\mathcal{V}(x)\cdot b(x,u) \le \kappa_{1} -
\frac{m\kappa_{0}}{8\overline\lambda(Q)} \bigl(\underline\lambda(Q) \bigr)^{{m}/{2}}\llvert x
\rrvert^{m}
\end{equation}
for some constant $\kappa_{1}>0$.
A similar calculation shows for some constant $\kappa_{2}>0$
it holds that
%
\begin{equation}
\label{E3.15nu}
\nabla\mathcal{V}(x)\cdot b(x,u) \le\kappa_{2} \bigl(1
+ \bigl[(e\cdot x)^{+} \bigr]^{m} \bigr)\qquad \forall x\in
\mathcal{K}\cap \bigl\{\llvert x\rrvert\ge1 \bigr\}.
\end{equation}
Also note that we can select $\kappa_{3}>0$ large enough such that
%
\begin{equation}
\label{eg-4}
\frac{1}{2} \bigl\llvert\operatorname{trace} \bigl(\Sigma
\Sigma^{\mathsf
{T}}\nabla^{2}\mathcal{V}(x) \bigr) \bigr\rrvert \le
\kappa_{3} +\frac{m\kappa_{0}}{16\overline\lambda(Q)} \bigl(\underline\lambda(Q)
\bigr)^{{m}/{2}} \llvert x\rrvert ^{m}.
\end{equation}
Hence, by \eqref{eg-1}--\eqref{eg-4}
there exists $\kappa_{4}>0$ such that
%
\begin{equation}
\label{eg-4b}
\quad\hspace*{6pt} L^{u}\mathcal{V}(x) \le \kappa_{4} -
\frac{m\kappa_{0}}{16\overline\lambda(Q)} \bigl(\underline\lambda(Q) \bigr)^{{m}/{2}} \llvert x
\rrvert^{m} \mathbb{I}_{\mathcal{K}^{c}}(x)+\kappa_{2}
\bigl[(e \cdot x)^{+} \bigr]^{m} \mathbb{I}_{\mathcal{K}}(x)
\end{equation}
for all $x\in\mathbb{R}^{d}$.
It is evident that we can scale $\mathcal{V}$, by multiplying it with
a constant,
so that \eqref{eg-4b} takes the form
%
\begin{equation}
\label{eg-4c}
L^{u}\mathcal{V}(x) \le 1 -c_{0} \llvert x
\rrvert^{m} \mathbb{I}_{\mathcal{K}^{c}}(x)+c_{1} \bigl[(e
\cdot x)^{+} \bigr]^{m} \mathbb{I}_{\mathcal{K}}(x)
\qquad\forall x\in \mathbb {R}^{d}.
\end{equation}
By \eqref{eg-cost}, the running cost $r$ is inf-compact on $\mathcal{K}$.
It then follows from \eqref{eg-cost} and~\eqref{eg-4c} that
\eqref{E-KA1} is satisfied with $h(x):=c_{0}\llvert x\rrvert^{m}$.

We next show that \eqref{eg-sde1} satisfies Assumption~\ref{Ass-2}.
Let
\[
u_{0}(\cdot) \equiv e_{d} = (0,\ldots,0,1)^{\mathsf{T}}.
\]
Then we can write \eqref{eg-sde1} as
%
\begin{equation}
\label{eg-sde2}
\mathrm{d} {X}_{t} = \bigl(\ell-R
\bigl(X_{t}-(e \cdot X_{t})^{+}u_{0}
\bigr) -{(e\cdot x)^{+}} \Gamma u_{0} \bigr) \,\mathrm{d}
{t}+\Sigma\,\mathrm{d} {W}_{t}.
\end{equation}
It is shown in \cite{dieker-gao} that the solution $X_{t}$
in \eqref{eg-sde2} is positive
recurrent and, therefore, $u_{0}$ is a stable Markov control.
This is done by finding a suitable
Lyapunov function. In particular, in \cite{dieker-gao}, Theorem~3,
it is shown that there exists
a positive definite matrix $\tilde Q$ such that if we define
%
\begin{equation}
\label{dieker-gao}
\psi(x) :=(e\cdot x)^{2}+ \tilde{\kappa}
\bigl[x-e_{d}\phi(e\cdot x) \bigr]^{\mathsf{T}}\tilde Q
\bigl[x-e_{d} \phi(e\cdot x) \bigr],
\end{equation}
for some suitably chosen constant $\tilde{\kappa}$ and a function
$\phi\in\mathcal{C}^{2}(\mathbb{R})$, given by
\[
\phi(y)=
\cases{\ds y, &\quad$\mbox{if }y \ge0$, \vspace*{2pt}
\cr
\ds
-\tfrac{1}{2}\tilde\delta, &\quad $\mbox{if } y \le-\tilde\delta$,
\vspace*{2pt}
\cr
\mbox{smooth}, & \quad$\mbox{if } -\tilde\delta<y<0$,}
\]
where $\tilde\delta>0$ is a suitable constant and $0\le\phi'(y)\le1$,
then it holds that
%
\begin{equation}
\label{dieker-gao-ly}
L^{u_{0}} \psi(x) \le-\tilde\kappa_{1} \llvert
x \rrvert^{2},
\end{equation}
for $\llvert x\rrvert$ large enough and positive constant $\tilde\kappa_{1}$.
We define $\mathcal{V}_{0}:=\mathrm{e}^{a\psi}$ where $a$ is to be
determined later.
Note that $\llvert\nabla\psi(x)\rrvert\le\tilde\kappa_{2}(1+\llvert
x\rrvert)$
for some constant $\tilde\kappa_{2}>0$.
Hence, a straightforward calculation shows that
if we choose $a$ small enough, then for some  constant
$\tilde\kappa_{3}>0$ it holds that
\begin{eqnarray*}
L^{u_{0}} \mathcal{V}_{0}(x) &\le & \bigl(-\tilde
\kappa_{1} a\llvert x\rrvert^{2}+a^{2}\lVert
\Sigma \rVert^{2} \tilde\kappa_{2} \bigl(1+\llvert x\rrvert
\bigr)^{2} \bigr)\mathcal{V}_{0}(x)
\\
& \le & -\tilde\kappa_{3}\llvert x\rrvert^{2}
\mathcal{V}_{0}(x),
\end{eqnarray*}
for all $\llvert x\rrvert$ large enough.
Since $\mathcal{V}_{0}(x)> [(e\cdot x)^{+}]^{m}$,
$m\ge1$, for all large enough $\llvert x\rrvert$ we see that $\mathcal
{V}_{0}$ satisfies
\eqref{E-lyap} with control $u_{0}$.
Hence, Assumption~\ref{Ass-2} holds by Lemma~\ref{L-basic}.
\end{pf}

\subsection{Existence of optimal controls}
\label{S3.4}

\begin{definition}
Recall the definition of $\varrho_{U}$ in \eqref{E-rU}.
For $\beta>0$, we define
\[
\mathfrak{U}^{\beta} := \bigl\{U\in\mathfrak{U} \dvtx
\varrho_{U}(x) \le\beta \mbox{ for some } x\in\mathbb{R}^{d}
\bigr\}.
\]
We also let $\mathfrak{U}_{\mathrm{SM}}^{\beta}:=\mathfrak{U}^{\beta}\cap\mathfrak{U}_{\mathrm{SM}}$, and
\begin{eqnarray*}
\hat\varrho_{*} & :=& \inf \bigl\{\beta>0 \dvtx \mathfrak{U}^{\beta}
\ne\varnothing \bigr\},
\\
\varrho_{*} & :=& \inf \bigl\{\beta>0 \dvtx \mathfrak{U}_{\mathrm
{SM}}^{\beta}
\ne\varnothing \bigr\},
\\
\tilde{\varrho}_{*} & :=& \inf \bigl\{\pi(r) \dvtx\pi\in\mathscr{G}
\bigr\},
\end{eqnarray*}
where
\[
\mathscr{G}:= \biggl\{\pi\in\mathcal{P} \bigl(\mathbb{R}^{d}\times
\mathbb{U} \bigr) \dvtx \int_{\mathbb{R}^{d}\times\mathbb{U}}L^{u} f(x) \pi(\mathrm{d} {x},\mathrm{d} {u})=0 \ \forall f\in\mathcal{C}^{\infty}_c
\bigl(\mathbb{R}^{d} \bigr) \biggr\},
\]
and $L^{u}f(x)$ is given by \eqref{E3.3}.
It is well known that $\mathscr{G}$ is the set of ergodic occupation measures
of the controlled process in \eqref{E-sde}, and that
$\mathscr{G}$ is a closed and convex subset of $\mathcal{P}(\mathbb
{R}^{d}\times\mathbb{U})$
\cite{ari-bor-ghosh}, Lemmas~3.2.2 and 3.2.3.
We use the notation $\pi_{v}$ when we want to indicate the ergodic occupation
measure associated with the control $v\in\mathfrak{U}_{\mathrm{SSM}}$.
In other words,
\[
\pi_{v}(\mathrm{d} {x},\mathrm{d} {u}) :=\mu_{v}(\mathrm{d} {x})v(\mathrm{d} {u}\vert x).
\]
\end{definition}

\begin{lemma}\label{L3.2}
If \eqref{E-lyap} holds for some $\mathcal{V}_{0}$ and $u_{0}$, then
we have
$\pi_{u_{0}}(r)\le\eta$.
Therefore,
$\hat\varrho_{*}<\infty$.
\end{lemma}

\begin{pf}
Let $(X_{t},u_{0}(X_{t}))$ be the solution of \eqref{E-sde}.
Recall that $\tau_{R}$ is the first exit time from $B_{R}$ for $R>0$.
Then by It\^{o}'s formula
\[
\mathbb{E}^{u_{0}}_{x} \bigl[\mathcal{V}_{0}(X_{T\wedge\tau
_{R}})
\bigr]-\mathcal{V}_{0}(x) \le \eta T -\mathbb{E}^{u_{0}}_{x}
\biggl[\int_{0}^{T\wedge\tau_{R}} r \bigl(X_{s},u_{0}(X_{s})
\bigr) \,\mathrm{d} {s} \biggr].
\]
Therefore, letting $R\to\infty$ and using Fatou's lemma, we obtain
the bound
\[
\mathbb{E}^{u_{0}}_{x} \biggl[\int_{0}^{T}r
\bigl(X_{s},u_{0}(X_{s}) \bigr) \,\mathrm{d} {s}
\biggr] \le\eta T+\mathcal{V}_{0}(x)-\min_{\mathbb{R}^{d}}
\mathcal{V}_{0},
\]
and thus
\[
\limsup_{T\to\infty} \frac{1}{T} \mathbb{E}^{u_{0}}_{x}
\biggl[\int_{0}^{T} r \bigl(X_{s},u_{0}(X_{s})
\bigr) \,\mathrm{d} {s} \biggr] \le\eta.
\]
\upqed\end{pf}

In the analysis, we use a function $\tilde{h}\in\mathcal{C}(\mathbb
{R}^{d}\times\mathbb{U})$
which, roughly speaking, is of the same order as $r$ in $\mathcal
{K}\times\mathbb{U}$
and lies between
$r$ and a multiple of $r+h$ on $\mathcal{K}^{c}\times\mathbb{U}$,
with $\mathcal{K}$ as
in Assumption~\ref{Ass-1}.
The existence of such a function is guaranteed by Assumption~\ref{Ass-1}
as the following lemma shows.

\begin{lemma}\label{L-key}
Define
\[
\mathcal{H}:= (\mathcal{K}\times\mathbb{U}) \cup \bigl\{(x,u)\in
\mathbb{R}^{d}\times\mathbb{U}\dvtx r(x,u)> h(x,u) \bigr\},
\]
where $\mathcal{K}$ is the open set in Assumption~\ref{Ass-1}.
Then there exists an inf-compact function
$\tilde{h}\in\mathcal{C}(\mathbb{R}^{d}\times\mathbb{U})$
which is locally Lipschitz in its first argument uniformly w.r.t. its
second argument, and satisfies
%
\begin{equation}
\label{E-KA3}
r(x,u) \le\tilde{h}(x,u) \le\frac{k_{0}}{2} \bigl(1+h(x,u)
\mathbb{I}_{\mathcal{H}^{c}}(x,u)+r(x,u) \mathbb {I}_{\mathcal{H}}(x,u) \bigr)
\end{equation}
for all $(x,u)\in\mathbb{R}^{d}\times\mathbb{U}$,
and for some positive constant $k_{0}\ge2$.
Moreover,
%
\begin{equation}
\label{E-key}
L^{u}\mathcal{V}(x) \le1 - h(x,u)
\mathbb{I}_{\mathcal{H}^{c}}(x,u)+r(x,u) \mathbb{I}_{\mathcal{H}}(x,u)
\end{equation}
for all $(x,u)\in\mathbb{R}^{d}\times\mathbb{U}$,
where $\mathcal{V}$ is the function in Assumption~\ref{Ass-1}.
\end{lemma}

\begin{pf}
If $f(x,u)$ denotes the right-hand side of \eqref{E-KA3},
with $k_{0}= 4$, then
\begin{eqnarray*}
f(x,u)-r(x,u)& >&  h(x,u) \mathbb{I}_{\mathcal{H}^{c}}(x,u)+r(x,u)
\mathbb{I}_{\mathcal
{H}}(x,u)
\\
& \ge &  h(x,u) \mathbb{I}_{\mathcal{K}^{c}}(x) +r(x,u) \mathbb{I}_{\mathcal{K}}(x),
\end{eqnarray*}
since $r(x,u)> h(x,u)$ on $\mathcal{H}\setminus(\mathcal{K}\times
\mathbb{U})$.
Therefore, by Assumption~\ref{Ass-1}, the set
$\{(x,u) \dvtx f(x,u)-r(x,u)\le n\}$ is bounded in $\mathbb{R}^{d}
\times\mathbb{U}$ for
every $n\in\mathbb{N}$.
Hence, there exists an increasing
sequence of open balls $D_{n}$, $n=1,2,\ldots,$ centered at $0$ in
$\mathbb{R}^{d}$
such that
$f(x,u)-r(x,u) \ge n$ for all $(x,u)\in D_{n}^{c}\times\mathbb{U}$.
Let $g\dvtx\mathbb{R}^{d} \to\mathbb{R}$ be any nonnegative
smooth function such that
$n-1 \le g(x)\le n$ for $x\in D_{n+1}\setminus D_{n}$, $n=1,2,\ldots,$
and $g(x)=0$ on $D_{1}$.
Clearly, $\tilde{h}:=r + g$ is
continuous, inf-compact,
locally Lipschitz in its first argument, and satisfies \eqref{E-KA3}.
That \eqref{E-key} holds is clear from \eqref{E-KA1} and the
fact that $\mathcal{H}\supseteq\mathcal{K}\times\mathbb{U}$.
\end{pf}

\begin{remark}
It is clear from the proof of Lemma~\ref{L-key}
that we could fix the value of the constant $k_{0}$ in \eqref{E-KA3},
say $k_{0}=4$.
However, we keep the variable $k_{0}$ because this provides some flexibility
in the choice of $\tilde{h}$, and also in order
to be able to trace it along the
different calculations.
\end{remark}

\begin{remark}\label{R3.4}
Note that if $h\ge r$ and $r$ is inf-compact, then
$\mathcal{H}=\mathcal{K}\times\mathbb{U}$ and
$\tilde{h}:=r$ satisfies \eqref{E-KA3}.
Note also, that in view of \eqref{eg-cost} and Proposition~\ref{eg-prop},
for the multi-class queueing model we have
\begin{eqnarray*}
r(x,u) & \le&  c_{2} \bigl(1+ \bigl[(e\cdot x)^{+}
\bigr]^{m} \bigr)
\\
& \le & \frac{c_{2} d^{m-1}}{1\wedge c_{0}} \bigl(1+ (1\wedge c_{0}) \llvert x
\rrvert^{m} \bigr)
\\
& \le & \frac{c_{2} d^{m-1}}{1\wedge c_{0}} \biggl(1+c_{0}\llvert x\rrvert^{m}
\mathbb{I}_{\mathcal{K}^{c}}(x) +\frac{1}{\delta^{m}} \bigl[{(e\cdot
x)^{+}} \bigr]^{m} \mathbb{I}_{\mathcal
{K}}(x) \biggr)
\\
& \le & \frac{c_{2} d^{m-1}}{1\wedge c_{0}} \biggl(1+h(x) \mathbb{I}_{\mathcal{K}^{c}}(x) +
\frac{1}{c_{1}\delta^{m}}r(x,u) \mathbb{I}_{\mathcal{K}}(x) \biggr)
\\
& \le & \frac{c_{2} d^{m-1}}{1\wedge c_{0}\wedge c_{1}\delta^{m}} \bigl(1+h(x) \mathbb{I}_{\mathcal{K}^{c}}(x)+r(x,u)
\mathbb {I}_{\mathcal{K}}(x) \bigr)
\\
& \le & \frac{c_{2} d^{m-1}}{1\wedge c_{0}\wedge c_{1}\delta^{m}} \bigl(1+h(x) \mathbb{I}_{\mathcal{H}^{c}}(x,u)+ r(x,u)
\mathbb {I}_{\mathcal{H}}(x,u) \bigr).
\end{eqnarray*}
Therefore,
$\tilde{h}(x,u):=c_{2}+c_{2} d^{m-1}\llvert x\rrvert^{m}$ satisfies
\eqref{E-KA3}.
\end{remark}

\begin{remark}\label{R-lsc}
We often use the fact that if $g\in\mathcal{C}(\mathbb{R}^{d}\times
\mathbb{U})$ is bounded below,
then the map $\mathcal{P}(\mathbb{R}^{d}\times\mathbb{U})\ni\nu
\mapsto\nu(g)$ is lower semi-continuous.
This easily follows from two facts: (a) $g$ can be expressed as an
increasing limit
of bounded continuous functions, and (b) if $g$ is bounded and continuous,
then $\pi\mapsto\pi(g)$ is continuous.
\end{remark}

\begin{theorem}\label{T-bst}
Let
$\beta\in(\hat\varrho_{*},\infty)$.
Then:
\begin{longlist}[(d)]
\item[(a)]
For all $U\in\mathfrak{U}^{\beta}$ and $x\in\mathbb{R}^{d}$ such that
$\varrho_{U}(x)\le\beta$, then
%
\begin{equation}
\label{E-Tbst01}
\limsup_{t\to\infty} \frac{1}{T}
\mathbb{E}^{U}_{x} \biggl[\int_{0}^{T}
\tilde {h}(X_{s},U_{s}) \,\mathrm{d} {s} \biggr] \le
k_{0}(1+\beta).
\end{equation}
\item[(b)]
$\hat\varrho_{*}=\varrho_{*}=\tilde\varrho_{*}$.
\item[(c)]
For any $\beta\in(\varrho_{*},\infty)$, we have $\mathfrak
{U}_{\mathrm{SM}}^{\beta}\subset\mathfrak{U}_{\mathrm{SSM}}$.

\item[(d)]
The set of invariant probability measures $\mathcal{M}^{\beta}$
corresponding to controls in $\mathfrak{U}_{\mathrm{SM}}^{\beta}$ satisfies
\[
\int_{\mathbb{R}^{d}} \tilde{h} \bigl(x,v(x) \bigr) \mu_{v}(
\mathrm{d} {x}) \le k_{0}(1+\beta)\qquad \forall\mu_{v}\in
\mathcal{M}^{\beta}.
\]
In particular, $\mathfrak{U}_{\mathrm{SM}}^{\beta}$ is tight in
$\mathcal{P}(\mathbb{R}^{d})$.

\item[(e)]
There exists $(\tilde{V},\tilde{\varrho})\in\mathcal
{C}^{2}(\mathbb{R}^{d})\times\mathbb{R}_+$,
with $\tilde{V}$ inf-compact, such that
%
\begin{equation}
\label{erg-h}
\min_{u\in\mathbb{U}} \bigl[L^{u}\tilde{V}(x)+
\tilde{h}(x,u) \bigr] = \tilde{\varrho}.
\end{equation}
\end{longlist}
\end{theorem}

\begin{pf}
Using It\^o's formula, it follows by \eqref{E-key} that
%
\begin{eqnarray}
&& \frac{1}{T} \bigl(\mathbb{E}^{U}_{x} \bigl[
\mathcal{V}(X_{T\wedge\tau
_{R}}) \bigr]-\mathcal{V}(x) \bigr)
\nonumber
\\
\label{E-Tbst1}
&& \qquad \le 1-\frac{1}{T} \mathbb{E}^{U}_{x}
\biggl[\int_{0}^{T\wedge\tau_{R}} h(X_{s},U_{s})
\mathbb{I}_{\mathcal{H}^{c}}(X_{s},U_{s}) \,\mathrm {d} {s}
\biggr]
\\
&& \qquad \quad {}+ \frac{1}{T} \mathbb{E}^{U}_{x}
\biggl[\int_{0}^{T\wedge\tau_{R}} r(X_{s},U_{s})
\mathbb{I}_{\mathcal{H}}(X_{s},U_{s}) \,\mathrm {d} {s}
\biggr].
\nonumber
\end{eqnarray}
Since $\mathcal{V}$ is inf-compact, \eqref{E-Tbst1} together
with \eqref{E-KA3} implies that
%
\begin{eqnarray}
&&\frac{2}{k_{0}} \limsup_{T\to\infty} \frac{1}{T}
\mathbb{E}^{U}_{x} \biggl[\int_{0}^{T}
\tilde{h}(X_{s},U_{s}) \,\mathrm {d} {s} \biggr]
\nonumber
\\[-8pt]
\label{E-Tbst2}
\\[-8pt]
\nonumber
&&\qquad{}\le 2+ 2 \limsup_{T\to\infty} \frac{1}{T}
\mathbb{E}^{U}_{x} \biggl[\int_{0}^{T}
r(X_{s},U_{s}) \,\mathrm{d} {s} \biggr].
\end{eqnarray}
Part (a) then follows from \eqref{E-Tbst2}.

Now fix
$U\in\mathfrak{U}^{\beta}$ and $x\in\mathbb{R}^{d}$ such that
$\varrho_{U}(x)\le\beta$.
The inequality in \eqref{E-Tbst01} implies that the set
of \emph{mean empirical measures} $\{\zeta^{U}_{x,t} \dvtx t\ge1\}$,
defined by
\[
\zeta^{U}_{x,t}(A\times B) := \frac{1}{t}
\mathbb{E}^{U}_{x} \biggl[\int_{0}^{t}
\mathbb{I}_{A\times
B}(X_{s},U_{s}) \,\mathrm{d} {s}
\biggr]
\]
for any Borel sets $A\subset\mathbb{R}^{d}$ and $B\subset\mathbb
{U}$, is tight.
It is the case that any limit point of the mean empirical measures
in $\mathcal{P}(\mathbb{R}^{d}\times\mathbb{U})$ is  an ergodic
occupation measure
\cite{ari-bor-ghosh}, Lemma~3.4.6.
Then in view of Remark~\ref{R-lsc} we obtain
%
\begin{equation}
\label{ETbst-A}
\pi(r) \le\limsup_{t\to\infty} \zeta^{U}_{x,t}(r)
\le\beta
\end{equation}
for some ergodic occupation measure $\pi$.
Therefore, $\tilde\varrho_{*}\le\hat\varrho_{*}$.
Disintegrating the measure $\pi$ as
$\pi(\mathrm{d}{x},\mathrm{d}{u})=v(\mathrm{d}{u}\vert x) \mu
_{v}(\mathrm{d}{x})$,
we obtain the associated control $v\in\mathfrak{U}_{\mathrm{SSM}}$.
From ergodic theory \cite{yosida}, we also know that
\[
\limsup_{T\to\infty} \frac{1}{T} \mathbb{E}^{v}_{x}
\biggl[\int_{0}^{T}r \bigl(X_{s},v(X_{s})
\bigr) \,\mathrm{d} {s} \biggr] = \pi_{v}(r)\qquad \mbox{for almost every } x.
\]
It follows that $\varrho_{*}\le\tilde\varrho_{*}$,
and since it is clear that
$\hat\varrho_{*}\le\varrho_{*}$, equality must hold among the three
quantities.

If $v\in\mathfrak{U}_{\mathrm{SM}}^{\beta}$, then \eqref{E-Tbst2} implies
that \eqref{ETbst-A} holds with $U\equiv v$ and $\pi\equiv\pi_{v}$.
Therefore, parts~(c) and (d) follow.

Existence of $(\tilde{V},\tilde{\varrho})$, satisfying \eqref{erg-h},
follows from Assumption~\ref{Ass-2} and \cite{ari-bor-ghosh}, Theorem~3.6.6.
The inf-compactness of $\tilde{V}$ follows from the
stochastic representation of $\tilde{V}$ in \cite{ari-bor-ghosh}, Lemma~3.6.9.
This proves~(e).
\end{pf}

Existence of a stationary Markov control that is optimal is
asserted by the following theorem.

\begin{theorem}\label{T-exist}
Let $\mathscr{G}$ denote the set of ergodic occupation measures corresponding
to controls in $\mathfrak{U}_{\mathrm{SSM}}$, and $\mathscr
{G}^{\beta}$ those corresponding to
controls in $\mathfrak{U}_{\mathrm{SM}}^{\beta}$, for $\beta
>\varrho_{*}$.
Then:
\begin{longlist}[(a)]
\item[(a)]
The set $\mathscr{G}^{\beta}$ is compact in $\mathcal{P}(\mathbb
{R}^{d})$ for any
$\beta>\varrho_{*}$.

\item[(b)]
There exists $v\in\mathfrak{U}_{\mathrm{SM}}$ such that $\varrho
_{v} = \varrho_{*}$.
\end{longlist}
\end{theorem}

\begin{pf}
By Theorem~\ref{T-bst}(d), the set $\mathscr{G}^{\beta}$ is tight
for any $\beta>\varrho_{*}$.
Let $\{\pi_{n}\}\subset\mathscr{G}^{\beta}$, for some $\beta
>\varrho_{*}$, be any
convergent sequence in $\mathcal{P}(\mathbb{R}^{d})$ such that
$\pi_{n}(r) \rightarrow\varrho_{*}$ as $n\to\infty$
and denote its limit by $\pi_{*}$.
Since $\mathscr{G}$ is closed, $\pi_{*}\in\mathscr{G}$,
and since the map
$\pi\to\pi(r)$
is lower semi-continuous, it follows that $\pi_{*}(r)\le\varrho_{*}$.
Therefore, $\mathscr{G}^{\beta}$ is closed, and hence compact.
Since $\pi(r)\ge\varrho_{*}$ for all $\pi\in\mathscr{G}$,
the equality $\pi_{*}(r)=\varrho_{*}$ follows.
Also $v$ is obtained by disintegrating $\pi_{*}$.
\end{pf}

\begin{remark}
The reader might have noticed at this point
that Assumption~\ref{Ass-1} may be weakened significantly.
What is really required is the existence of an open set
$\hat{\mathcal{H}}\subset\mathbb{R}^{d}\times\mathbb{U}$ and
inf-compact functions $\mathcal{V}\in\mathcal{C}^{2}(\mathbb{R}^{d})$
and $h \in\mathcal{C}(\mathbb{R}^{d}\times\mathbb{U})$,
satisfying
\begin{longlist}[(H2)]
\item[(H1)]
$\inf_{\{u\dvtx(x,u)\in\hat{\mathcal{H}}\}} r(x,u)
\ds\mathop{\longrightarrow}_{\llvert x\rrvert\to\infty}\infty$.

\item[(H2)]
$L^{u}\mathcal{V}(x) \le1 - h(x,u) \mathbb{I}_{\hat{\mathcal{H}}^{c}}(x,u)
+r(x,u) \mathbb{I}_{\hat{\mathcal{H}}}(x,u)
\ \forall(x,u)\in\mathbb{R}^{d}\times\mathbb{U}$.
\end{longlist}
In (H1), we use the convention that the `$\inf$' of the empty set
is $+\infty$.
Also note that (H1) is equivalent to the statement that
$\{(x,u) \dvtx r(x,u)\le c\} \cap\hat{\mathcal{H}}$ is bounded in
$\mathbb{R}^{d}\times\mathbb{U}$
for all $c\in\mathbb{R}_{+}$.
If (H1)--(H2) are met, we define
$\mathcal{H}:= \hat{\mathcal{H}} \cup \{(x,u)\in\mathbb
{R}^{d}\times\mathbb{U}\dvtx
r(x,u)> h(x,u) \}$,
and following the proof of Lemma~\ref{L-key},
we assert the existence of an inf-compact
$\tilde{h}\in\mathcal{C}(\mathbb{R}^{d}\times\mathbb{U})$
satisfying \eqref{E-KA3}.
In fact, throughout the rest of the paper,
Assumption~\ref{Ass-1} is not really invoked.
We only use \eqref{E-key}, the inf-compact function $\tilde{h}$ satisfying
\eqref{E-KA3}, and, naturally, Assumption~\ref{Ass-2}.
\end{remark}

\subsection{The HJB equation}\label{S3.5}

For $\varepsilon>0$, let
\[
r_{\varepsilon}(x,u) :=r(x,u) + \varepsilon\tilde{h}(x,u).
\]
By Theorem~\ref{T-bst}(d), for any $\pi\in\mathscr{G}^{\beta}$,
$\beta>\varrho_{*}$,
we have the bound
%
\begin{equation}
\label{ET-exist1}
\pi(r_{\varepsilon}) \le \beta+ \varepsilon k_{0}(1+
\beta).
\end{equation}
Therefore, since $r_{\varepsilon}$ is near-monotone, that is,
\[
\liminf_{|x|\to\infty} \min_{u \in\mathbb{U}} r_{\varepsilon}
(x,u) > \inf_{\pi\in\mathscr{G}} \pi(r_{\varepsilon}),
\]
there exists
$\pi_{\varepsilon}\in\mathop{\argmin}_{\pi\in\mathscr{G}} \pi
(r_{\varepsilon})$.
Let $\pi_{*}\in\mathscr{G}$ be as in the proof of Theorem~\ref{T-exist}.
The sub-optimality of $\pi_{*}$ relative
to the running cost $r_{\varepsilon}$ and \eqref{ET-exist1} imply that
%
\begin{eqnarray}
\pi_{\varepsilon}(r) & \le & \pi_{\varepsilon}(r_{\varepsilon})
\nonumber
\\
\label{ET-exist2} & \le & \pi_{*}(r_{\varepsilon})
\\
& \le &  \varrho_{*}+ \varepsilon k_{0}(1+
\varrho_{*})\qquad \forall\varepsilon>0.
\nonumber
\end{eqnarray}
It follows from \eqref{ET-exist2} and Theorem~\ref{T-bst}(d)
that $ \{\pi_{\varepsilon} \dvtx \varepsilon\in(0,1)
\}$
is tight.
Since $\pi_{\varepsilon}\mapsto\pi_{\varepsilon}(r)$ is
lower semi-continuous, if $\bar{\pi}$ is any limit
point of $\pi_{\varepsilon}$ as $\varepsilon\searrow0$,
then taking limits in \eqref{ET-exist2}, we obtain
%
\begin{equation}
\label{ET-exist3}
\bar{\pi}(r) \le\limsup_{\varepsilon\searrow0} \pi
_{\varepsilon}(r) \le\varrho_{*}.
\end{equation}
Since $\mathscr{G}$ is closed, $\bar{\pi}\in\mathscr{G}$, which
implies that
$\bar{\pi}(r)\ge\varrho_{*}$.
Therefore, equality must hold in \eqref{ET-exist3}, or in other words,
$\bar{\pi}$ is an optimal ergodic occupation measure.

\begin{theorem}\label{T-HJB1}
There exists a unique function $V^{\varepsilon}\in\mathcal
{C}^{2}(\mathbb{R}^{d})$
with $V^{\varepsilon}(0)=0$,
which is bounded below in $\mathbb{R}^{d}$,
and solves the HJB
%
\begin{equation}
\label{E-HJBe}
\min_{u\in\mathbb{U}} \bigl[L^{u}V^{\varepsilon}(x)+r_{\varepsilon
}(x,u)
\bigr] = \varrho_{\varepsilon},
\end{equation}
where $\varrho_{\varepsilon}:=\inf_{\pi\in\mathscr{G}} \pi
(r_{\varepsilon})$,
or in other words, $\varrho_{\varepsilon}$ is the optimal value
of the ergodic control problem with running cost $r_{\varepsilon}$.
Also a stationary Markov control $v_{\varepsilon}$ is optimal
for the ergodic control problem relative to $r_{\varepsilon}$
if and only if it satisfies
%
\begin{equation}
\label{E-ve}
\qquad H_{\varepsilon} \bigl(x,\nabla V^{\varepsilon}(x) \bigr) = b
\bigl(x,v_{\varepsilon}(x) \bigr)\cdot\nabla V^{\varepsilon}(x) +
r_{\varepsilon} \bigl(x,v_{\varepsilon}(x) \bigr) \qquad\mbox{a.e. in }
\mathbb{R}^{d},
\end{equation}
where
%
\begin{equation}
\label{E-He}
H_{\varepsilon}(x,p) :=\min_{u\in\mathbb{U}}
\bigl[b(x,u) \cdot p + r_{\varepsilon}(x,u) \bigr].
\end{equation}
Moreover:
\begin{longlist}[(a)]
\item[(a)]
for every $R>0$, there exists $k_{R}$ such that
%
\begin{equation}
\label{s-014}
\mathop{\operatorname{osc}}_{B_{R}} V^{\varepsilon} \le
k_{R};
\end{equation}

\item[(b)]
if $v_{\varepsilon}$ is a measurable a.e. selector
from the minimizer of the Hamiltonian in~\eqref{E-He}, that is,
if it satisfies \eqref{E-HJBe}, then for any $\delta>0$,
\[
V^{\varepsilon}(x) \ge\mathbb{E}^{v_{\varepsilon}}_{x} \biggl[\int
_{0}^{{\breve\tau}_{\delta}} \bigl(r_{\varepsilon
}
\bigl(X_{s},v_{\varepsilon}(X_{s}) \bigr) -
\varrho_{\varepsilon} \bigr) \,\mathrm{d} {s} \biggr] +\inf_{B_{\delta}}
V^{\varepsilon};
\]

\item[(c)]
for any stationary control $v\in\mathfrak{U}_{\mathrm{SSM}}$ and for
any $\delta>0$,
\[
V^{\varepsilon}(x) \le \mathbb{E}^{v}_{x} \biggl[\int
_{0}^{{\breve\tau}_{\delta}} \bigl(r_{\varepsilon}
\bigl(X_{s},v(X_{s}) \bigr) -\varrho_{\varepsilon} \bigr)
\,\mathrm{d} {s} + V^{\varepsilon
}(X_{{\breve\tau}_{\delta}}) \biggr],
\]
where ${\breve\tau}_{\delta}$ is hitting time to the ball
$B_{\delta}$.
\end{longlist}
\end{theorem}

\begin{theorem}\label{T-HJB2}
Let $V^{\varepsilon}$, $\varrho_{\varepsilon}$, and $v_{\varepsilon}$,
for $\varepsilon>0$, be as in Theorem~\ref{T-HJB1}.
The following hold:
\begin{longlist}[(a)]
\item[(a)]
The function $V^{\varepsilon}$ converges
to some $V_{*}\in\mathcal{C}^{2}(\mathbb{R}^{d})$, uniformly on
compact sets,
and $\varrho_{\varepsilon}\to\varrho_{*}$, as $\varepsilon\searrow0$,
and $V_{*}$ satisfies
%
\begin{equation}
\label{E-HJB}
\min_{u\in\mathbb{U}} \bigl[L^{u}V_{*}(x)+r(x,u)
\bigr] = \varrho_{*}.
\end{equation}
Also, any limit point $v_{*}$ (in the topology of Markov controls)
as $\varepsilon\searrow0$ of the set
$\{v_{\varepsilon}\}$ satisfies
\[
L^{v_{*}}V_{*}(x)+r \bigl(x,v_{*}(x) \bigr) =
\varrho_{*} \qquad\mbox{a.e. in }\mathbb{R}^{d}.
\]
\item[(b)]
A stationary Markov control $v$ is optimal
for the ergodic control problem relative to $r$
if and only if it satisfies
%
\begin{equation}
\label{E-v}
\hspace*{5pt}H \bigl(x,\nabla V_{*}(x) \bigr) = b \bigl(x,v(x)
\bigr)\cdot\nabla V_{*}(x) + r \bigl(x,v(x) \bigr)\qquad \mbox{a.e. in }
\mathbb{R}^{d},
\end{equation}
where
\[
H(x,p) :=\min_{u\in\mathbb{U}} \bigl[b(x,u)\cdot p + r(x,u) \bigr].
\]
Moreover, for an optimal $v\in\mathfrak{U}_{\mathrm{SM}}$, we have
\[
\lim_{T\to\infty} \frac{1}{T} \mathbb{E}^{v}_{x}
\biggl[\int_{0}^{T}r \bigl(X_{s},v(X_{s})
\bigr) \,\mathrm{d} {s} \biggr] = \varrho_{*}\qquad \forall x\in
\mathbb{R}^{d}.
\]
\item[(c)]
The function $V_{*}$ has the stochastic representation
%
\begin{eqnarray}
V_{*}(x) & =&  \lim_{\delta\searrow0} \inf_{v \in\bigcup_{\beta>0} \mathfrak
{U}_{\mathrm{SM}}^{\beta}}
\mathbb{E}^{v}_{x} \biggl[\int_{0}^{{\breve\tau}_{\delta}}
\bigl(r \bigl(X_{s},v(X_{s}) \bigr)-\varrho_{*}
\bigr) \,\mathrm{d} {s} \biggr]
\nonumber
\\[-8pt]
\label{E-strep}
\\[-8pt]
\nonumber
& =& \lim_{\delta\searrow0} \mathbb{E}^{\bar{v}}_{x}
\biggl[\int_{0}^{{\breve\tau}_{\delta}} \bigl(r \bigl(X_{s},v_{*}(X_{s})
\bigr)-\varrho_{*} \bigr) \,\mathrm{d} {s} \biggr]
\end{eqnarray}
for any $\bar{v}\in\mathfrak{U}_{\mathrm{SM}}$ that satisfies
\eqref{E-v}.

\item[(d)]
If $\mathbb{U}$ is a convex set, $u\mapsto\{b(x,u)\cdot p + r(x,u)\}$
is strictly convex whenever it is not constant, and
$u\mapsto\tilde{h}(x,u)$
is strictly convex for all $x$, then any measurable minimizer of \eqref{E-HJBe}
converges pointwise, and thus in $\mathfrak{U}_{\mathrm{SM}}$, to the
minimizer of~\eqref{E-HJB}.
\end{longlist}
\end{theorem}

Theorem~\ref{T-HJB2} guarantees the existence of an optimal stable control,
which is made precise by \eqref{E-v}, for the ergodic diffusion
control problem
with the running cost function $r$.
Moreover, under the convexity property in part~(d),
the optimal stable control can
be obtained as a pointwise limit from the minimizing
selector of~\eqref{E-HJBe}.
For instance, if we let
\[
r(x,u) = (e\cdot x)^{+}\sum_{i=1}^{d}
h_{i} u_{i}^{m},\qquad m>1,
\]
then by choosing
$h$ and $\tilde{h}+|u|^{2}$ as in Proposition~\ref{eg-prop}, we see that
the approximate value function $V^{\varepsilon}$ and approximate control
$v_{\varepsilon}$
converge to the desired optimal value function $V_{*}$ and optimal
control $v_{*}$,
respectively.

Concerning the uniqueness of the solution to the HJB equation in
\eqref{E-HJB}, recall that in the near-monotone case
the existing uniqueness results are as follows:
there exists a unique solution pair $(V,\varrho)$ of \eqref{E-HJB}
with $V$ in the class of functions $\mathcal{C}^{2}(\mathbb{R}^{d})$
which are bounded
below in $\mathbb{R}^{d}$. Moreover, it satisfies $V(0)=0$ and
$\varrho\le\varrho_{*}$.
If the restriction $\varrho\le\varrho_{*}$ is removed, then in general,
there are multiple solutions.
Since in our model $r$ is not near-monotone in $\mathbb{R}^{d}$, the function
$V_{*}$ is not, in general, bounded below.
However, as we show later in Lemma~\ref{Luniq1} the negative
part of $V_{*}$ grows slower than $\mathcal{V}$, that is, it holds that
$V_{*}^{-}\in\mathfrak{o}(\mathcal{V})$, with $\mathfrak{o}(\cdot)$
as defined in Section~\ref{S-notation}.
Therefore, the second part of the theorem that follows
may be viewed as an extension of the well-known
uniqueness results that apply to ergodic control problems
with near-monotone running cost.
The third part of the theorem resembles the hypotheses of uniqueness
that apply to problems under a blanket stability hypothesis.

\begin{theorem}\label{T-unique}
Let $(\hat{V},\hat\varrho)$ be a solution of
%
\begin{equation}
\label{E-HJB-hat2} \min_{u\in\mathbb{U}} \bigl[L^{u}
\hat{V}(x)+r(x,u) \bigr] = \hat \varrho,
\end{equation}
such that $\hat{V}^{-}\in\mathfrak{o}(\mathcal{V})$ and $\hat{V}(0)=0$.
Then the following hold:
\begin{longlist}[(a)]
\item[(a)]
Any measurable selector $\hat{v}$ from the minimizer of the
associated Hamiltonian in \eqref{E-v}
is in $\mathfrak{U}_{\mathrm{SSM}}$ and $\varrho_{\hat{v}} < \infty$.

\item[(b)]
If $\hat\varrho\le\varrho_{*}$ then\vspace*{1.5pt} necessarily
$\hat\varrho=\varrho_{*}$ and $\hat{V}=V_{*}$.

\item[(c)]
If
$\hat{V}\in\mathscr{O} (\min_{u\in\mathbb{U}} \tilde
{h}(\cdot,u) )$,
then $\hat\varrho=\varrho_{*}$ and $\hat{V}=V_{*}$.
\end{longlist}
\end{theorem}

Applying these results to the multi-class queueing
diffusion model, we have the following corollary.

\begin{corollary}\label{C-unique}
For the queueing diffusion model with controlled dynamics given by
\eqref{eg-sde1}, drift given by \eqref{q-drift},
and running cost as in \eqref{eg-cost},
there exists a unique solution $V$, satisfying $V(0)=0$, to
the associated HJB in the class of functions
$\mathcal{C}^{2}(\mathbb{R}^{d})\cap\mathscr{O} (\llvert
x\rrvert^{m} )$,
whose negative part is in $\mathfrak{o} (\llvert x\rrvert^{m} )$.
This solution agrees with $V_{*}$ in Theorem~\ref{T-HJB2}.
\end{corollary}

\begin{pf}
Existence of a solution $V$ follows by Theorem~\ref{T-HJB2}.
Select $\mathcal{V}\sim\llvert x\rrvert^{m}$ as in the proof of
Proposition~\ref{eg-prop}.
That the solution $V$ is in the stated class then
follows by Lemma~\ref{Luniq1} and Corollary~\ref{C-exist}
that appear later in Sections~\ref{Proofs} and~\ref{S-truncation}, respectively.
With $h\sim\llvert x\rrvert^{m}$ as in the proof of Proposition~\ref{eg-prop},
it follows that
$\min_{u\in\mathbb{U}} \tilde{h}(x,u)\in\mathscr{O} (\llvert
x\rrvert^{m} )$.
Therefore, uniqueness follows by Theorem~\ref{T-unique}.
\end{pf}

We can also obtain the HJB equation in \eqref{E-HJB}
via the traditional vanishing discount approach
as the following theorem asserts.
Similar results are shown for a one-dimensional degenerate ergodic diffusion
control problem in \cite{Ocone-Weerasinghe} and certain multi-dimensional
ergodic diffusion control problems (allowing degeneracy and spatial periodicity)
in \cite{Arisawa-Lions}.

\begin{theorem}\label{T-HJB3}
Let $V_{*}$ and $\varrho_{*}$ be as in Theorem~\ref{T-HJB2}.
For $\alpha>0$, we define
\[
V_{\alpha}(x) :=\inf_{U\in\mathfrak{U}} \mathbb{E}^{U}_{x}
\biggl[\int_{0}^\infty\mathrm{e}^{-\alpha t}
r(X_{t},U_{t}) \,\mathrm{d} {t} \biggr].
\]
The function $V_{\alpha} - V_{\alpha}(0)$ converges, as $\alpha
\searrow0$,
to $V_{*}$, uniformly on compact subsets of $\mathbb{R}^{d}$.
Moreover, $\alpha V_{\alpha}(0) \to\varrho_{*}$, as $\alpha\searrow0$.
\end{theorem}

The proofs of the Theorems~\ref{T-HJB1}--\ref{T-HJB3}
are given in Section~\ref{Proofs}.
The following result, which follows directly from
\eqref{ET-exist2}, provides a way to find $\varepsilon$-optimal controls.

\begin{proposition}
Let $\{v_{\varepsilon}\}$ be the minimizing selector from Theorem~\ref
{T-HJB1} and
$\{\mu_{v_{\varepsilon}}\}$ be the corresponding invariant
probability measures.
Then almost surely for all $x\in\mathbb{R}^{d}$,
\begin{eqnarray*}
\lim_{T\to\infty} \frac{1}{T} \mathbb{E}^{v_{\varepsilon}}_{x}
\biggl[\int_{0}^{T} r \bigl(X_{s},v_{\varepsilon}(X_{s})
\bigr) \,\mathrm{d} {s} \biggr]& =& \int_{\mathbb{R}^{d}} r
\bigl(x,v_{\varepsilon}(x) \bigr) \mu_{v_{\varepsilon
}}(\mathrm{d} {x})
\\
& \le & \varrho_{*} + \varepsilon k_{0}(1+
\varrho_{*}).
\end{eqnarray*}
\end{proposition}

\subsection{Technical proofs}
\label{Proofs}

Recall that $r_{\varepsilon}(x,u)=r(x,u)+\varepsilon\tilde{h}(x,u)$, with
$\tilde{h}$ as in Lemma~\ref{L-key}.
We need the following lemma.

For $\alpha>0$ and $\varepsilon\ge0$, we define
%
\begin{equation}
\label{E-dcost}
V^{\varepsilon}_{\alpha}(x) :=\inf_{U\in\mathfrak{U}}
\mathbb{E}^{U}_{x} \biggl[\int_{0}^\infty
\mathrm{e}^{-\alpha t} r_{\varepsilon}(X_{t},U_{t})
\,\mathrm{d} {t} \biggr],
\end{equation}
where we set $r_{0}\equiv r$.
Clearly, when $\varepsilon=0$, we have $V^{0}_{\alpha}\equiv
V_{\alpha}$.

We quote the
following result from \cite{ari-bor-ghosh}, Theorem~3.5.6, Remark~3.5.8.

\begin{lemma}\label{L3.3}
Provided $\varepsilon>0$, then
$V^{\varepsilon}_{\alpha}$ defined above is in $\mathcal
{C}^{2}(\mathbb{R}^{d})$ and is the minimal
nonnegative solution of
\[
\min_{u\in\mathbb{U}} \bigl[L^{u} V^{\varepsilon}_{\alpha
}(x)+r_{\varepsilon}(x,u)
\bigr] = \alpha V^{\varepsilon}_{\alpha}(x).
\]
\end{lemma}

The HJB in Lemma~\ref{L3.3} is similar to the equation
in \cite{atar-mandel-rei}, Theorem~3,  which concerns
the characterization of the discounted control problem.

\begin{lemma}\label{L3.4}
Let $u$ be any precise Markov control and $L^{u}$ be the
corresponding generator.
Let $\varphi\in\mathcal{C}^{2}(\mathbb{R}^{d})$ be a nonnegative
solution of
\[
L^{u}\varphi-\alpha\varphi= g,
\]
where $g\in L^\infty_{\mathrm{loc}}(\mathbb{R}^{d})$.
Let $\kappa\dvtx\mathbb{R}_{+}\to\mathbb{R}_{+}$ be any
nondecreasing function
such that
$\lVert g\rVert_{L^\infty(B_{R})}\le\kappa(R)$ for all $R>0$.
Then for any $R>0$ there exists a
constant $D(R)$ which depends on $\kappa(4R)$,
but not on $u$, or $\varphi$, such that
\[
\mathop{\operatorname{osc}}_{B_{R}} \varphi\le D(R) \Bigl(1+\alpha\inf
_{B_{4R}} \varphi \Bigr).
\]
\end{lemma}

\begin{pf}
Define $\tilde g:=\alpha(g-2\kappa(4R))$ and
$\tilde\varphi:=2\kappa(4R)+\alpha\varphi$.
Then $\tilde g\le0$ in $B_{4R}$ and $\tilde\varphi$ solves
\[
L^{u}\tilde\varphi-\alpha\tilde\varphi= \tilde g
\qquad\mbox{in }B_{4R}.
\]
Also
\begin{eqnarray*}
\lVert\tilde{g}\rVert_{L^\infty(B_{4R})} & \le&  \alpha \bigl(2\kappa(4R)+\lVert g
\rVert_{L^\infty(B_{4R})} \bigr)
\\
& \le & 3\alpha \bigl(2\kappa(4R)- \lVert g\rVert_{L^\infty
(B_{4R})} \bigr)
\\
& =&  3 \inf_{B_{4R}} \llvert\tilde{g}\rrvert
\\
& \le & 3 |B_{4R}|^{-1}\lVert\tilde{g}\rVert_{L^{1}(B_{4R})}.
\end{eqnarray*}
Hence by \cite{ari-bor-ghosh}, Theorem~A.2.13,  there exists a positive constant
$\tilde C_{H}$ such that
\[
\sup_{x\in B_{3R}} \tilde\varphi(x) \le\tilde C_{H} \inf
_{x\in B_{3R}} \tilde\varphi(x),
\]
implying that
%
\begin{equation}
\label{p-009}
\alpha\sup_{x\in B_{3R}} \varphi(x) \le \tilde
C_{H} \Bigl(2\kappa(4R)+\inf_{x\in B_{3R}} \alpha\varphi
(x) \Bigr).
\end{equation}
We next consider the solution of
\[
L^{u}\psi= 0 \qquad\mbox{in }B_{3R},\qquad \psi= \varphi \qquad\mbox{on }\partial B_{3R}.
\]
Then
\[
L^{u}(\varphi-\psi) = \alpha\varphi+g \qquad\mbox{in } B_{3R}.
\]
If $\varphi(\hat x)=\inf_{x\in B_{3R}} \varphi(x)$, then applying the
maximum principle (\cite{ari-bor-ghosh}, Theorem~A.2.1, \cite{gilbarg-trudinger})
it follows from \eqref{p-009} that
%
\begin{equation}
\label{p-010}
\sup_{x\in B_{3R}} \llvert\varphi-\psi\rrvert \le
\hat{C} \bigl(1+\alpha\varphi(\hat x) \bigr).
\end{equation}
Again $\psi$ attains its minimum at the boundary
(\cite{ari-bor-ghosh}, Theorem~A.2.3, \cite{gilbarg-trudinger}).
Therefore, $\psi-\varphi(\hat x)$ is a nonnegative
function, and hence by the Harnack inequality, there exists a constant $C_{H}>0$
such that
\[
\psi(x)-\varphi(\hat x) \le C_{H} \bigl(\psi(\hat x)-\varphi(\hat x)
\bigr) \le C_{H}\hat{C} \bigl(1+\alpha\varphi(\hat x) \bigr)\qquad \forall
x\in B_{2R}.
\]
Thus, combining the above display with \eqref{p-010} we obtain
\[
\mathop{\operatorname{osc}}_{B_{2R}} \varphi\le\sup_{B_{2R}} (
\varphi-\psi) +\sup_{B_{2R}} \psi-\varphi(\hat x) \le
\hat{C}(1+C_{H}) \bigl(1+\alpha\varphi(\hat x) \bigr).
\]
This completes the proof.
\end{pf}

\begin{lemma}\label{L3.5}
Let $V^{\varepsilon}_{\alpha}$ be as in Lemma~\ref{L3.3}.
Then for any $R>0$, there exists a constant $k_{R}>0$ such that
\[
\mathop{\operatorname{osc}}_{B_{R}} V^{\varepsilon}_{\alpha} \le
k_{R} \qquad\mbox{for all } \alpha\in(0,1]\mbox{ and } \varepsilon\in[0,1].
\]
\end{lemma}

\begin{pf}
Recall that $\mu_{u_{0}}$ is the stationary probability distribution
for the process under the control $u_{0}\in\mathfrak{U}_{\mathrm
{SSM}}$ in Lemma~\ref{L-basic}.
Since $u_{0}$ is sub-optimal for
the $\alpha$-discounted criterion in \eqref{E-dcost}, and
$V^{\varepsilon}_{\alpha}$ is nonnegative, then
for any ball $B_{R}$, using Fubini's theorem, we obtain
\begin{eqnarray*}
\mu_{u_{0}}(B_{R}) \inf_{B_{R}}
V^{\varepsilon}_{\alpha} & \le &\int_{\mathbb{R}^{d}}V^{\varepsilon}_{\alpha}(x)
\mu _{u_{0}}(\mathrm{d} {x})
\\
& \le &\int_{\mathbb{R}^{d}} \mathbb{E}^{u_{0}}_{x}
\biggl[\int_{0}^\infty\mathrm{e}^{-\alpha t}
r_{\varepsilon} \bigl(X_{t},u_{0}(X_{t}) \bigr)
\,\mathrm{d} {t} \biggr] \mu_{u_{0}}(\mathrm{d} {x})
\\
& =& \frac{1}{\alpha} \mu_{u_{0}}(r_{\varepsilon})
\\
& \le& \frac{1}{\alpha} \bigl(\eta+ \varepsilon k_{0}(1+\eta) \bigr),
\end{eqnarray*}
where for the last inequality we used Lemma~\ref{L3.2} and
Theorem~\ref{T-bst}(a).

Therefore, we have the estimate
\[
\alpha\inf_{B_{R}} V^{\varepsilon}_{\alpha} \le
\frac{\eta+ \varepsilon k_{0}(1+\eta)}{\mu_{u_{0}}(B_{R})}.
\]
The result then follows by Lemma~\ref{L3.4}.
\end{pf}

We continue with the proof of Theorem~\ref{T-HJB1}.

\begin{pf*}{Proof of Theorem~\ref{T-HJB1}}
Consider the function
$\bar{V}^{\varepsilon}_{\alpha}:=V^{\varepsilon}_{\alpha}
-V^{\varepsilon}_{\alpha}(0)$.
In view of Lemma~\ref{L3.4} and Lemma~\ref{L3.5}, we see that
$\bar{V}^{\varepsilon}_{\alpha}$ is locally bounded uniformly in
$\alpha\in(0,1]$ and $\varepsilon\in(0,1]$.
Therefore, by standard elliptic theory, $\bar{V}^{\varepsilon
}_{\alpha}$
and its first- and second-order
partial derivatives are uniformly bounded in $L^{p}(B)$, for any $p>1$,
in any
bounded ball $B\subset\mathbb{R}^{d}$, that is, for some constant
$C_{B}$ depending on $B$ and $p$,
$ \lVert\bar{V}^{\varepsilon}_{\alpha} \rVert_{\mathscr
{W}^{2,p}(B)}\le C_{B}$
(\cite{gilbarg-trudinger}, Theorem~9.11, page~117).
Therefore, we can extract a subsequence along which
$\bar{V}^{\varepsilon}_{\alpha}$ converges.
Then the result follows from Theorems~3.6.6, Lemma~3.6.9 and Theorem~3.6.10
in \cite{ari-bor-ghosh}.
The proof of \eqref{s-014} follows from Lemma~\ref{L3.4} and
Lemma~\ref{L3.5}.
\end{pf*}

\begin{remark}\label{R-tight}
In the proof of the following lemma, and elsewhere in
the paper, we use the fact that if $\mathcal{U}\subset\mathfrak
{U}_{\mathrm{SSM}}$ is
a any set of controls such that the corresponding set
$\{\mu_{v} \dvtx v\in\mathcal{U}\}\subset\mathcal{M}$
of invariant probability measures is tight then
the map $v\mapsto\pi_{v}$ from the closure of
$\mathcal{U}$ to $\mathcal{P}(\mathbb{R}^{d}\times\mathbb{U})$
is continuous, and so is the map $v\mapsto\mu_{v}$.
In fact, the latter is continuous under the total variation norm topology
\cite{ari-bor-ghosh}, Lemma~3.2.6.
We also recall that $\mathscr{G}$ and $\mathcal{M}$ are closed and
convex subsets
of $\mathcal{P}(\mathbb{R}^{d}\times\mathbb{U})$ and $\mathcal
{P}(\mathbb{R}^{d})$.
Therefore,$\{\pi_{v} \dvtx v\in\bar{\mathcal{U}}\}$ is compact in
$\mathscr{G}$.
Note also that since $\mathbb{U}$ is compact, tightness of a set of invariant
probability measures is equivalent to tightness of the corresponding
set of ergodic occupation measures.
\end{remark}

\begin{lemma}\label{L3.6}
If $\{v_{\varepsilon} \dvtx\varepsilon\in(0,1]\}$ is a collection of
measurable selectors from the minimizer of \eqref{E-HJBe},
then the corresponding invariant probability measures
$\{\mu_{{\varepsilon}} \dvtx\varepsilon\in(0,1]\}$ are tight.
Moreover, if $v_{\varepsilon_{n}}\to v_{*}$
along some subsequence $\varepsilon_{n}\searrow0$, then the following hold:
\begin{longlist}[(a)]
\item[(a)]
$\mu_{\varepsilon_{n}}\to\mu_{v_{*}}$ as
$\varepsilon_{n}\searrow0$,

\item[(b)]
$v_{*}$ is a stable Markov control,

\item[(c)]
$\int_{\mathbb{R}^{d}}r(x,v_{*}(x)) \mu_{v_{*}}(\mathrm{d}{x}) =
\lim_{\varepsilon\searrow0} \varrho_{\varepsilon} = \varrho_{*}$.
\end{longlist}
\end{lemma}

\begin{pf}
By \eqref{E-Tbst01} and \eqref{ET-exist2},
the set of ergodic occupation measures corresponding to
$\{v_{\varepsilon} \dvtx\varepsilon\in(0,1]\}$ is tight.
By Remark~\ref{R-tight}, the same applies to the set
$\{\mu_\varepsilon\dvtx\varepsilon\in(0,1]\}$,
and also part~(a) holds.
Part~(b) follows from the equivalence of the existence of
an invariant probability measure for a controlled diffusion
and the stability of the associated stationary Markov control
(see \cite{ari-bor-ghosh}, Theorem~2.6.10).
Part~(c) then follows since equality holds in \eqref{ET-exist3}.
\end{pf}

We continue with the following lemma that asserts the continuity of the
mean hitting time of a ball with respect to the stable Markov controls.

\begin{lemma}\label{L3.7}
Let $\{v_{n} \dvtx n\in\mathbb{N}\}\subset\mathfrak{U}_{\mathrm{SM}}^{\beta}$, for some
$\beta>0$, be a collection of Markov controls such that
$v_{n}\to\hat{v}$ in the topology of Markov controls
as \mbox{$n\to\infty$}.
Let $\mu_{n}$, $\hat{\mu}$ be the invariant probability measures
corresponding to the controls~$v_{n}$,~$\hat{v}$, respectively.
Then for any $\delta>0$, it holds that
\[
\mathbb{E}^{v_{n}}_{x}[{\breve\tau}_{\delta}] \mathop{\longrightarrow}_{n\to\infty} \mathbb{E}^{\hat{v}}_{x}[{\breve
\tau}_{\delta}] \qquad\forall x\in B^{c}_{\delta}.
\]
\end{lemma}

\begin{pf}
Define $H(x):=\min_{u\in\mathbb{U}} \tilde{h}(x,u)$.
It is easy to see that
$H$ is inf-compact and locally Lipschitz.
Therefore, by Theorem~\ref{T-bst}(d) we have
\[
\sup_{n\in\mathbb{N}} \mu_{n}(H) \le k_{0}(1+\beta),
\]
and since $\mu_{n}\to\hat{\mu}$, we also have $\hat{\mu}(H) \le
k_{0}(1+\beta)$.
Then by \cite{ari-bor-ghosh}, Lemma~3.3.4, we obtain
%
\begin{equation}
\label{p-020}
\sup_{n\in\mathbb{N}} \mathbb{E}^{v_{n}}_{x}
\biggl[\int_{0}^{{\breve\tau}_{\delta}} H(X_{s})
\,\mathrm{d} {s} \biggr] +\mathbb{E}^{\hat{v}}_{x} \biggl[\int
_{0}^{{\breve\tau}_{\delta
}}H(X_{s}) \,\mathrm{d} {s}
\biggr] < \infty.
\end{equation}
Let $R$ be a positive number greater than $\llvert x\rrvert$.
Then by \eqref{p-020}, there exists a positive~$k$ such that
\[
\mathbb{E}^{v}_{x} \biggl[\int_{0}^{{\breve\tau}_{\delta
}}
\mathbb{I}_{\{H>R\}}(X_{s}) \,\mathrm{d} {s} \biggr] \le
\frac{1}{R} \mathbb{E}^{v}_{x} \biggl[\int
_{0}^{{\breve\tau
}_{\delta}}H(X_{s}) \mathbb{I}_{\{H>R\}}(X_{s})
\,\mathrm{d} {s} \biggr] \le\frac{k}{R}
\]
for $v\in\{\{v_{n}\}, \hat{v} \}$.
From this assertion and \eqref{p-020}, we see that
\[
\sup_{v\in\{\{v_{n}\}, \hat{v} \}} \mathbb{E}^{v}_{x} \biggl[
\int_{0}^{{\breve\tau}_{\delta
}}\mathbb{I}_{\{H>R\}}(X_{s})
\,\mathrm{d} {s} \biggr] \mathop{\longrightarrow}_{R\to\infty} 0.
\]
Therefore, in order to prove the lemma it is enough to show that, for
any $R>0$,
we have
\[
\mathbb{E}^{v_{n}}_{x} \biggl[\int_{0}^{{\breve\tau}_{\delta}}
\mathbb{I}_{\{H\le R\}}(X_{s}) \,\mathrm{d} {s} \biggr] \mathop{
\longrightarrow}_{n\to\infty} \mathbb{E}^{\hat{v}}_{x} \biggl[\int
_{0}^{{\breve\tau}_{\delta}} \mathbb{I}_{\{H\le R\}}(X_{s})
\,\mathrm{d} {s} \biggr].
\]
But this follows from \cite{ari-bor-ghosh}, Lemma~2.6.13(iii).
\end{pf}

\begin{lemma}\label{L3.8}
Let $(V^{\varepsilon},\varrho_{\varepsilon})$ be as in Theorem~\ref{T-HJB1},
and $v_{\varepsilon}$ satisfy \eqref{E-He}.
There exists a subsequence $\varepsilon_{n}\searrow0$,
such that $V^{\varepsilon_{n}}$ converges
to some $V_{*}\in\mathcal{C}^{2}(\mathbb{R}^{d})$, uniformly on
compact sets,
and $V_{*}$ satisfies
%
\begin{equation}
\label{EL3.7-HJB}
\min_{u\in\mathbb{U}} \bigl[L^{u}V_{*}(x)+r(x,u)
\bigr] = \varrho_{*}.
\end{equation}
Also, any limit point $v_{*}$ $($in the topology of Markov controls$)$
of the set $\{v_{\varepsilon}\}$, as $\varepsilon\searrow0$, satisfies
%
\begin{equation}
\label{EL3.7-p-40}
L^{v_{*}}V_{*}(x)+r \bigl(x,v_{*}(x)
\bigr) = \varrho_{*}\qquad \mbox{a.e. in } \mathbb{R}^{d}.
\end{equation}
Moreover, $V_{*}$ admits the stochastic representation
%
\begin{eqnarray}
V_{*}(x)& =& \inf_{v \in\bigcup_{\beta>0} \mathfrak{U}_{\mathrm{SM}}^{\beta}} \mathbb{E}^{v}_{x}
\biggl[\int_{0}^{{\breve\tau}_{\delta}} \bigl(r \bigl(X_{s},v(X_{s})
\bigr)-\varrho_{*} \bigr) \,\mathrm{d} {s} +V_{*}(X_{{\breve\tau}_{\delta}})
\biggr]
\nonumber
\\[-8pt]
\label{EL3.7-strep}
\\[-8pt]
\nonumber
& =& \mathbb{E}^{v_{*}}_{x} \biggl[\int
_{0}^{{\breve\tau}_{\delta}} \bigl(r \bigl(X_{s},v_{*}(X_{s})
\bigr)- \varrho_{*} \bigr) \,\mathrm{d} {s} +V_{*}(X_{{\breve\tau}_{\delta}})
\biggr].
\end{eqnarray}
It follows that $V_{*}$ is the unique
limit point of $V^{\varepsilon}$ as $\varepsilon\searrow0$.
\end{lemma}

\begin{pf}
From \eqref{s-014}, we see that the family
$\{V^{\varepsilon} \dvtx\varepsilon\in(0,1]\}$ is uniformly
locally bounded.
Hence, applying the theory of elliptic PDE, it follows that
$\{V^{\varepsilon} \dvtx\varepsilon\in(0,1]\}$ is uniformly
bounded in
$\mathscr{W}_{\mathrm{loc}}^{2,p}(\mathbb{R}^{d})$ for $p>d$.
Consequently, $\{V^{\varepsilon} \dvtx\varepsilon\in(0,1]\}$
is uniformly bounded
in $\mathcal{C}_{\mathrm{loc}}^{1,\gamma}$ for some $\gamma>0$.
Therefore, along some subsequence $\varepsilon_{n}\searrow0$,
$V^{\varepsilon_{n}}\to V_{*}\in\mathscr{W}^{2,p}\cap\mathcal
{C}^{1,\gamma}$,
as $n\to\infty$, uniformly on compact sets.
Also, $\lim_{\varepsilon\searrow0} \varrho_{\varepsilon}=\varrho_{*}$
by Lemma~\ref{L3.5}(c).
Therefore, passing to the limit we obtain the HJB equation in
\eqref{EL3.7-HJB}.
It is straightforward to verify that~\eqref{EL3.7-p-40} holds
\cite{ari-bor-ghosh}, Lemma~2.4.3.

By Theorem~\ref{T-HJB1}(c), taking limits as $\varepsilon\searrow0$,
we obtain
%
\begin{equation}
\label{EL3:upper}
V_{*}(x) \le \inf_{v \in\bigcup_{\beta>0} \mathfrak{U}_{\mathrm{SM}}^{\beta}}
\mathbb{E}^{v}_{x} \biggl[\int_{0}^{{\breve\tau}_{\delta}}
\bigl(r \bigl(X_{s},v(X_{s}) \bigr)-\varrho_{*}
\bigr) \,\mathrm{d} {s} +V_{*}(X_{{\breve\tau}_{\delta}}) \biggr].
\end{equation}
Also by Theorem~\ref{T-HJB1}(b) we have the bound
\[
V^{\varepsilon}(x) \ge -\varrho_{\varepsilon} \mathbb{E}^{v_{\varepsilon}}_{x}[{
\breve\tau}_{\delta}] +\inf_{B_{\delta}} V^{\varepsilon}.
\]
Using Lemma~\ref{L3.7} and taking limits as $\varepsilon_{n}\searrow0$,
we obtain the lower bound
%
\begin{equation}
\label{EL3.7-p-42} V_{*}(x) \ge-\varrho_{*}
\mathbb{E}^{v_{*}}_{x}[{\breve\tau }_{\delta}] +\inf
_{B_{\delta}} V_{*}.
\end{equation}
By Lemma~\ref{L3.6}(c) and Theorem~\ref{T-bst}(d), $v_{*}\in
\mathfrak{U}_{\mathrm{SSM}}$,
and $\pi_{v_{*}}(\tilde{h})\le k_{0}(1+\varrho_{*})$.
Define
%
\begin{equation}
\label{EL3.7-3001}
\varphi(x) :=\mathbb{E}^{v_{*}}_{x} \biggl[\int
_{0}^{{\breve\tau}_{\delta}} \tilde{h} \bigl(X_{s},v_{*}(X_{s})
\bigr) \,\mathrm{d} {s} \biggr].
\end{equation}
For $\llvert x\rrvert>\delta$, we have
\[
\mathbb{E}^{v_{*}}_{x} \bigl[\mathbb{I}_{\{\tau_{R}<{\breve\tau
}_{\delta}\}}
\varphi(X_{\tau_{R}}) \bigr] = \mathbb{E}^{v_{*}}_{x}
\biggl[\mathbb{I}_{\{\tau_{R}<{\breve
\tau}_{\delta}\}} \int_{\tau_{R}\wedge{\breve\tau}_{\delta}}^{{\breve\tau
}_{\delta}}
\tilde{h} \bigl(X_{s},v_{*}(X_{s}) \bigr)
\,\mathrm{d} {s} \biggr].
\]
Therefore, by the dominated convergence theorem and the fact that
$\varphi(x)<\infty$ we obtain
\[
\mathbb{E}^{v_{*}}_{x} \bigl[ \varphi(X_{\tau_{R}})
\mathbb{I}_{\{\tau_{R}<{\breve\tau
}_{\delta}\}} \bigr] \mathop{\longrightarrow}_{R\nearrow\infty} 0.
\]
By \eqref{EL3:upper} and \eqref{EL3.7-p-42}, we have $|V_{*}|\in
\mathscr{O}(\varphi)$.
Thus \eqref{EL3.7-p-42} and \eqref{EL3.7-3001} imply that
\[
\liminf_{R\nearrow\infty} \mathbb{E}^{v_{*}}_{x}
\bigl[V_{*}(X_{\tau_{R}}) \mathbb{I}_{\{\tau_{R}<{\breve\tau
}_{\delta}\}} \bigr] = 0,
\]
and thus
%
\begin{equation}
\label{EEE3.50}
\liminf_{R\nearrow\infty} \mathbb{E}^{v_{*}}_{x}
\bigl[V_{*}(X_{\tau_{R}\wedge{\breve\tau}_{\delta}}) \bigr] = \mathbb{E}^{v_{*}}_{x}
\bigl[V_{*}(X_{{\breve\tau}_{\delta
}}) \bigr].
\end{equation}
Applying It\^{o}'s formula to \eqref{EL3.7-p-40}, we obtain
%
\begin{equation}
\label{EEE3.51}
V_{*}(x) = \mathbb{E}^{v_{*}}_{x}
\biggl[\int_{0}^{{\breve\tau}_{\delta
}\wedge\tau_{R}} \bigl(r \bigl(X_{s},v_{*}(X_{s})
\bigr)-\varrho_{*} \bigr) \,\mathrm{d} {s} + V_{*}(X_{{\breve\tau}_{\delta}\wedge\tau_{R}})
\biggr].
\end{equation}
Taking limits as $R\to\infty$, and using the dominated convergence theorem,
we obtain \eqref{EL3.7-strep} from \eqref{EL3:upper}.
\end{pf}

Recall the definition of $\mathfrak{o}(\cdot)$ from Section~\ref{S-notation}.
We need the following lemma.

\begin{lemma}\label{Luniq1}
Let $V_{*}$ be as in Lemma~\ref{L3.8}.
It holds that $V_{*}^{-}\in\mathfrak{o}(\mathcal{V})$.
\end{lemma}

\begin{pf}
Let $v_{*}$ be as in Lemma~\ref{L3.8}.
Applying It\^{o}'s formula to \eqref{E-key} with $u\equiv v_{*}$
we obtain
%
\begin{eqnarray}
&&\hspace*{7pt}\quad\mathbb{E}^{v_{*}}_{x} \biggl[\int_{0}^{{\breve\tau}_{\delta}}
h \bigl(X_{s},v_{*}(X_{s}) \bigr)
\mathbb{I}_{\mathcal
{H}^{c}} \bigl(X_{s},v_{*}(X_{s})
\bigr) \,\mathrm{d} {s} \biggr]
\nonumber
\\[-8pt]
\label{Euniq1a}
\\[-8pt]
\nonumber
&&\hspace*{7pt} \quad\qquad\le \mathbb{E}^{v_{*}}_{x} \biggl[\int
_{0}^{{\breve\tau}_{\delta}} r \bigl(X_{s},v_{*}(X_{s})
\bigr) \mathbb{I}_{\mathcal
{H}} \bigl(X_{s},v_{*}(X_{s})
\bigr) \,\mathrm{d} {s} \biggr] +\mathbb{E}^{v_{*}}_{x}[{\breve
\tau}_{\delta}]+\mathcal{V}(x).
\end{eqnarray}
Therefore, adding the term
\[
\mathbb{E}^{v_{*}}_{x} \biggl[\int_{0}^{{\breve\tau}_{\delta}}
r \bigl(X_{s},v_{*}(X_{s}) \bigr)
\mathbb{I}_{\mathcal
{H}} \bigl(X_{s},v_{*}(X_{s})
\bigr) \,\mathrm{d} {s} \biggr] -(1+2\varrho_{*})\mathbb{E}^{v_{*}}_{x}[{
\breve\tau}_{\delta}]
\]
to both sides of \eqref{Euniq1a} and using the stochastic
representation of
$V_{*}$ we obtain
%
\begin{eqnarray}
F(x) & :=& 2k_{0}^{-1}\mathbb{E}^{v_{*}}_{x}
\biggl[\int_{0}^{{\breve
\tau}_{\delta}} \tilde{h}
\bigl(X_{s},v_{*}(X_{s}) \bigr) \,\mathrm{d} {s}
\biggr] -2(1+\varrho_{*})\mathbb{E}^{v_{*}}_{x}[{
\breve\tau}_{\delta}]
\nonumber
\\[-8pt]
\label{Euniq1b}
\\[-8pt]
\nonumber
& \le& 2V_{*}(x)+\mathcal{V}(x) - 2 \inf
_{B_{\delta}} V_{*}.
\end{eqnarray}
From the stochastic representation of $V_{*}$ we have
$V_{*}^{-}(x)\le\varrho_{*}\mathbb{E}^{v_{*}}_{x}[{\breve\tau
}_{\delta}]
- \inf_{B_{\delta}} V_{*}$.
For any $R>\delta$, we have
%
\begin{equation}
\label{Euniq1c}
\mathbb{E}^{v_{*}}_{x} \biggl[\int
_{0}^{{\breve\tau}_{\delta}} \tilde{h} \bigl(X_{s},v_{*}(X_{s})
\bigr) \,\mathrm{d} {s} \biggr] \ge \Bigl(\inf_{B_{R}^{c}\times\mathbb{U}} \tilde{h}
\Bigr) \mathbb {E}_{x}[{\breve\tau}_{R}] \qquad\forall x\in
B_{R}^{c}.
\end{equation}
It is also straightforward to show that
$\lim_{\llvert x\rrvert\to\infty}\frac{\mathbb{E}_{x}[{\breve
\tau}_{R}]}{\mathbb{E}_{x}[{\breve\tau}_{\delta}]}=1$.
Therefore, since $\tilde{h}$ is inf-compact,
it follows by \eqref{Euniq1b} and \eqref{Euniq1c} that the map
$x\mapsto\mathbb{E}^{v_{*}}_{x}[{\breve\tau}_{\delta}]$
is in $\mathfrak{o}(F)$,
which implies that $V_{*}^{-}\in\mathfrak{o}(F)$.
On the other hand, by \eqref{Euniq1b} we obtain
$F(x) \le\mathcal{V}(x) - 2\sup_{B_{\delta}} V_{*}$ for all $x$
such that
$V_{*}(x)\le0$, which implies that the restriction of $F$ to the
support of
$V_{*}^{-}$ is in $\mathscr{O}(\mathcal{V})$.
It follows that $V_{*}^{-}\in\mathfrak{o}(\mathcal{V})$.
\end{pf}

We next prove Theorem~\ref{T-HJB2}.

\begin{pf*}{Proof of Theorem~\ref{T-HJB2}}
Part~(a) is contained in Lemma~\ref{L3.8}.

To prove part~(b), let $\bar{v}$ be any control satisfying \eqref{E-v}.
By Lemma~\ref{Luniq1} the map
$\mathcal{V}+ 2V_{*}$ is inf-compact and
by Theorem~\ref{T-HJB2} and \eqref{E-key} it satisfies
\begin{eqnarray*}
L^{\bar{v}} (\mathcal{V}+ 2V_{*}) (x) & \le & 1+ 2
\varrho_{*} - r \bigl(x,\bar{v}(x) \bigr) - h \bigl(x,\bar{v}(x) \bigr)
\mathbb{I}_{\mathcal
{H}^{c}} \bigl(x,\bar{v}(x) \bigr)
\\
& \le &  2+2\varrho_{*} - 2k_{0}^{-1} \tilde{h}
\bigl(x,\bar{v}(x) \bigr) \qquad\forall x\in\mathbb{R}^{d}.
\nonumber
\end{eqnarray*}
This implies that $\bar{v}\in\mathfrak{U}_{\mathrm{SSM}}$.
Applying It\^o's formula, we obtain
%
\begin{equation}
\label{416a} \limsup_{T\to\infty} \frac{1}{T}
\mathbb{E}^{\bar{v}}_{x} \biggl[\int_{0}^{T}
\tilde{h} \bigl(X_{s},\bar{v}(X_{s}) \bigr) \,\mathrm{d}
{s} \biggr] \le k_{0}(1+\varrho_{*}).
\end{equation}
Therefore, $\pi_{\bar{v}} (\tilde{h}) < \infty$.
By \eqref{E-key}, we have
\[
\mathbb{E}^{\bar{v}}_{x} \bigl[\mathcal{V}(X_{t})
\bigr] \le \mathcal{V}(x) + t + \mathbb{E}^{\bar{v}}_{x} \biggl[
\int_{0}^{t} r \bigl(X_{s},
\bar{v}(X_{s}) \bigr) \,\mathrm{d} {s} \biggr],
\]
and since $r \le\tilde{h}$, this implies by \eqref{416a} that
%
\begin{equation}
\label{416b}
\limsup_{T\to\infty} \frac{1}{T}
\mathbb{E}^{\bar{v}}_{x} \bigl[\mathcal{V}(X_{T}) \bigr]
\le 1+k_{0}(1+\varrho_{*}).
\end{equation}
Since $V_{*}^{-}\in\mathfrak{o}(\mathcal{V})$, it follows by
\eqref{416b} that
\[
\limsup_{T\to\infty} \frac{1}{T} \mathbb{E}^{\bar{v}}_{x}
\bigl[V_{*}^{-}(X_{T}) \bigr] = 0.
\]
Therefore, by It\^o's formula, we deduce from \eqref{E-HJB} that
%
\begin{equation}
\label{E-T3.4new} \limsup_{T\to\infty} \frac{1}{T}
\mathbb{E}^{\bar{v}}_{x} \biggl[\int_{0}^{T}
r \bigl(X_{s},\bar{v}(X_{s}) \bigr) \,\mathrm{d} {s}
\biggr] \le\varrho_{*}.
\end{equation}
On the other hand, since the only limit point of the mean empirical
measures $\zeta^{\bar{v}}_{x,t}$, as $t\to\infty$, is $\pi_{\bar{v}}$,
and $\pi_{\bar{v}}(r)=\varrho_{*}$, then in view of Remark~\ref{R-lsc},
we obtain $\liminf_{t\to\infty} \zeta^{\bar{v}}_{x,t}(r)\ge
\varrho_{*}$.
This proves that equality holds in \eqref{E-T3.4new} and that the
``$\limsup$'' may be replaced with ``$\lim$.''

Conversely, suppose $v\in\mathfrak{U}_{\mathrm{SM}}$ is optimal but
does not satisfy \eqref{E-v}.
Then there exists $R>0$ and
a nontrivial nonnegative $f\in L^{\infty}(B_{R})$ such that
\[
f_{\varepsilon}(x) :=\mathbb{I}_{B_{R}}(x) \bigl(L^{v}V^{\varepsilon}(x)
+ r_{\varepsilon} \bigl(x,v(x) \bigr) - \varrho _{\varepsilon} \bigr)
\]
converges to $f$, weakly in $L^{1}(B_{R})$,
along some subsequence $\varepsilon\searrow0$.
By applying It\^o's formula to \eqref{E-HJBe}, we obtain
%
\begin{eqnarray}
&&\frac{1}{T} \bigl(\mathbb{E}^{v}_{x}
\bigl[V^{\varepsilon}(X_{T\wedge\tau
_{R}}) \bigr] -V^{\varepsilon}(x) \bigr) +
\frac{1}{T} \mathbb{E}^{v}_{x} \biggl[\int
_{0}^{T\wedge\tau_{R}} r_{\varepsilon} \bigl(X_{s},v(X_{s})
\bigr) \,\mathrm{d} {s} \biggr]
\nonumber
\\[-8pt]
\label{EEEE3.54}
\\[-8pt]
\nonumber
&& \qquad \ge\varrho_{\varepsilon}+\frac{1}{T} \mathbb{E}^{v}_{x}
\biggl[\int_{0}^{T\wedge\tau_{R}} f_{\varepsilon}
\bigl(X_{s},v(X_{s}) \bigr) \,\mathrm{d} {s} \biggr].
\end{eqnarray}
Define, for some $\delta>0$,
\[
G(x) :=\mathbb{E}^{v}_{x} \biggl[\int_{0}^{{\breve\tau}_{\delta}}
r_{\varepsilon} \bigl(X_{s},v(X_{s}) \bigr) \,\mathrm{d}
{s} \biggr].
\]
Since $V^{\varepsilon}$ is bounded from below, by Theorem~\ref
{T-HJB1}(c)
we have $V^{\varepsilon}\in\mathscr{O}(G)$.
Invoking~\cite{ari-bor-ghosh}, Corollary~3.7.3, we obtain
\[
\lim_{T\to\infty} \frac{1}{T} \mathbb{E}^{v}_{x}
\bigl[V^{\varepsilon}(X_{T}) \bigr] = 0,
\]
and
\[
\lim_{R\to\infty} \mathbb{E}^{v}_{x}
\bigl[V^{\varepsilon}(X_{T\wedge\tau_{R}}) \bigr]  = \mathbb{E}^{v}_{x}
\bigl[V^{\varepsilon}(X_{T}) \bigr].
\]
Therefore, taking limits in \eqref{EEEE3.54}, first as $R\nearrow
\infty$,
and then as $T\to\infty$, we obtain
%
\begin{equation}
\label{E-discrepancy}
\pi_{v}(r_{\varepsilon}) \ge\varrho_{\varepsilon}+
\pi_{v}(f_{\varepsilon}).
\end{equation}
Taking limits as $\varepsilon\searrow0$ in \eqref{E-discrepancy},
since $\mu_{v}$ has a strictly positive density in $B_{R}$,
we obtain
\[
\pi_{v}(r) \ge\varrho_{*} + \pi_{v}(f) >
\varrho_{*},
\]
which is a contradiction.
This completes the proof of part~(b).

The first equality \eqref{E-strep} follows by Lemma~\ref{L3.8},
taking limits as $\delta\searrow0$.
To show that the second equality holds for any optimal control,
suppose $\bar{v}$ satisfies \eqref{E-v}.
By \eqref{E-key} we have, for $\delta>0$ and $\llvert x\rrvert>\delta$,
\[
\mathbb{E}^{\bar{v}}_{x} \bigl[\mathcal{V}(X_{\tau_{R}})
\mathbb{I}_{\{\tau_{R}<{\breve\tau}_{\delta}\}} \bigr] \le \mathcal{V}(x) + \sup_{B_{\delta}}
\mathcal{V}^{-} + \mathbb {E}^{\bar{v}}_{x} \biggl[
\int_{0}^{\tau_{R}\wedge{\breve\tau}_{\delta}} \bigl(1+r \bigl(X_{s},
\bar{v}(X_{s}) \bigr) \bigr) \,\mathrm{d} {s} \biggr].
\]
It follows that (see \cite{ari-bor-ghosh}, Lemma~3.3.4)
\[
\limsup_{R\to\infty} \mathbb{E}^{\bar{v}}_{x} \bigl[
\mathcal {V}(X_{\tau_{R}}) \mathbb{I}_{\{\tau_{R}<{\breve\tau}_{\delta}\}} \bigr] < \infty,
\]
and since $V_{*}^{-}\in\mathfrak{o}(\mathcal{V})$ we must have
\[
\limsup_{R\to\infty} \mathbb{E}^{\bar{v}}_{x}
\bigl[V_{*}^{-}(X_{\tau_{R}}) \mathbb{I}_{\{\tau_{R}<{\breve\tau}_{\delta}\}}
\bigr] = 0.
\]
By the first equality in \eqref{EL3.7-strep}, we obtain
$V^{+}_{*}\in\mathscr{O}(\varphi)$,
with $\varphi$ as defined in \eqref{EL3.7-3001} with $v_{*}$ replaced
by $\bar{v}$.
Thus, in analogy to \eqref{EEE3.50}, we obtain
\[
\liminf_{R\nearrow\infty} \mathbb{E}^{\bar{v}}_{x}
\bigl[V_{*}(X_{\tau_{R}\wedge{\breve\tau}_{\delta}}) \bigr] = \mathbb{E}^{\bar{v}}_{x}
\bigl[V_{*}(X_{{\breve\tau}_{\delta
}}) \bigr].
\]
The rest follows as in the proof of Lemma~\ref{L3.8} via \eqref{EEE3.51}.

We next prove part~(d). We assume that $\mathbb{U}$ is a convex set
and that
\[
c(x,u,p) := \bigl\{b(x,u)\cdot p+ r(x,u) \bigr\}
\]
is strictly convex in $u$ if
it is not identically a constant for fixed $x$ and $p$.
We fix some point $\bar{u}\in\mathbb{U}$.
Define
\[
\mathcal{B}:= \bigl\{x\in\mathbb{R}^{d} \dvtx c(x,\cdot,p)= c(x,
\bar{u}, p) \mbox{ for all } p \bigr\}.
\]
It is easy to see that on $\mathcal{B}$ both $b$ and $r$ do not depend
on $u$.
It is also easy to check that $\mathcal{B}$ is a closed set.
Let $(V_{*},v_{*})$ be the limit of $(V^{\varepsilon},v_{\varepsilon})$,
where $V_{*}$
is the solution to \eqref{E-HJB} and $v_{*}$ is the corresponding
limit of $v_{\varepsilon}$.
We have already shown that $v_{*}$ is a stable Markov control.
We next show that it is, in fact, a precise Markov control.
By our assumption, $v_{\varepsilon}$ is the
unique minimizing selector in \eqref{E-ve} and, moreover,
$v_{\varepsilon}$ is continuous in $x$.
By the definition of $r_{\varepsilon}$ it is clear that the restriction
of $v_{\varepsilon}$ to $\mathcal{B}$ does not depend on $\varepsilon$.
Let $v_{\varepsilon}(x)=v'(x)$ on $\mathcal{B}$.
Using the strict convexity property of
$c(x,\cdot,\nabla V_{*})$
it is easy to verify that $v_{\varepsilon}$ converges to the unique minimizer
of~\eqref{E-HJB} on $\mathcal{B}^{c}$.
In fact, since $\mathcal{B}^{c}$ is open, then for any sequence
$x^\varepsilon\to x\in\mathcal{B}^{c}$
it holds that $v_{\varepsilon}(x^\varepsilon)\to v_{*}(x)$.
This follows from the definition of the minimizer and the uniform
convergence of
$\nabla V^{\varepsilon}$ to $\nabla V_{*}$.
Therefore, we see that $v_{*}$ is a
precise Markov control, $v_{*}=v'$ on $\mathcal{B}$, and
$v_{\varepsilon}\to v_{*}$ pointwise
as $\varepsilon\to0$.
It is also easy to check that pointwise convergence implies convergence
in the topology of Markov controls.
\end{pf*}

We now embark on the proof of Theorem~\ref{T-unique}.

\begin{pf*}{Proof of Theorem~\ref{T-unique}}
The hypothesis that $\hat{V}^{-}\in\mathfrak{o}(\mathcal{V})$ implies
that the map
$\mathcal{V}+ 2\hat{V}$ is inf-compact.
Also by \eqref{E-key} and \eqref{E-HJB-hat2}, it satisfies
\begin{eqnarray*}
L^{\hat{v}} (\mathcal{V}+ 2\hat{V}) (x) & \le & 1+ 2\hat\varrho - r \bigl(x,
\hat{v}(x) \bigr) - h \bigl(x,\hat{v}(x) \bigr) \mathbb{I}_{\mathcal
{H}^{c}}
\bigl(x, \hat{v}(x) \bigr)
\\
& \le&  2+2\hat\varrho- 2k_{0}^{-1}\tilde{h} \bigl(x,
\hat{v}(x) \bigr)\qquad \forall x\in\mathbb{R}^{d}.
\end{eqnarray*}
Therefore, $\int\tilde{h}(x,\hat{v}(x)) \,\mathrm{d}\pi_{\hat{v}}
< \infty$ from which it
follows that $\varrho_{\hat{v}} < \infty$.
This proves part~(a).

By \eqref{E-key}, we have
\[
\mathbb{E}^{\hat{v}}_{x} \bigl[\mathcal{V}(X_{t})
\bigr] \le \mathcal{V}(x) + t + \mathbb{E}^{\hat{v}}_{x} \biggl[
\int_{0}^{t} r \bigl(X_{s},
\hat{v}(X_{s}) \bigr) \,\mathrm{d} {s} \biggr],
\]
and since $\varrho_{\hat{v}} < \infty$, this implies that
%
\begin{equation}
\label{Luniq2a} \limsup_{T\to\infty} \frac{1}{T}
\mathbb{E}^{\hat{v}}_{x} \bigl[\mathcal{V}(X_{T}) \bigr]
\le 1+\varrho_{\hat{v}}.
\end{equation}
Since $\hat{V}^{-}\in\mathfrak{o}(\mathcal{V})$, it follows by
\eqref{Luniq2a} that
\[
\limsup_{T\to\infty} \frac{1}{T} \mathbb{E}^{\hat{v}}_{x}
\bigl[\hat{V}^{-}(X_{T}) \bigr] = 0.
\]
Therefore, by It\^o's formula, we deduce from \eqref{E-HJB-hat2} that
\[
\limsup_{T\to\infty} \biggl(\frac{1}{T} \mathbb{E}^{\hat{v}}_{x}
\bigl[\hat{V}^{+}(X_{T}) \bigr] + \frac{1}{T}
\mathbb{E}^{\hat{v}}_{x} \biggl[\int_{0}^{T}
r \bigl(X_{s},\hat{v}(X_{s}) \bigr) \,\mathrm{d} {s}
\biggr] \biggr) = \hat{\varrho}.
\]
This implies that $\varrho_{\hat{v}}\le\hat{\varrho}$ and since by
hypothesis
$\hat\varrho\le\varrho_{*}$ we must have $\hat\varrho=\varrho_{*}$.

Again by \eqref{E-key}, we have
\[
\mathbb{E}^{\hat{v}}_{x} \bigl[\mathcal{V}(X_{\tau_{R}})
\mathbb{I}_{\{\tau_{R}<{\breve\tau}_{\delta}\}} \bigr] \le \mathcal{V}(x) + \sup_{B_{\delta}}
\mathcal{V}^{-} + \mathbb {E}^{\hat{v}}_{x} \biggl[
\int_{0}^{\tau_{R}\wedge{\breve\tau}_{\delta}} \bigl(1+r \bigl(X_{s},
\hat{v}(X_{s}) \bigr) \bigr) \,\mathrm{d} {s} \biggr].
\]
It follows by \cite{ari-bor-ghosh}, Lemma~3.3.4, that
\[
\limsup_{R\to\infty} \mathbb{E}^{\hat{v}}_{x} \bigl[
\mathcal {V}(X_{\tau_{R}}) \mathbb{I}_{\{\tau_{R}<{\breve\tau}_{\delta}\}} \bigr] < \infty,
\]
and since $\hat{V}^{-}\in\mathfrak{o}(\mathcal{V})$ we must have
%
\begin{equation}
\label{Luniq2b}
\limsup_{R\to\infty} \mathbb{E}^{\hat{v}}_{x}
\bigl[\hat {V}^{-}(X_{\tau_{R}}) \mathbb{I}_{\{\tau_{R}<{\breve\tau}_{\delta}\}}
\bigr] = 0.
\end{equation}
Using \eqref{Luniq2b} and following the steps in the proof of the second
equality in
\eqref{EL3.7-strep}, we obtain
\begin{eqnarray*}
\hat{V}(x) & \ge &  \mathbb{E}^{\hat{v}}_{x} \biggl[\int
_{0}^{{\breve\tau}_{\delta}} \bigl(r \bigl(X_{s},
\hat{v}(X_{s}) \bigr)-\varrho_{*} \bigr) \,\mathrm{d} {s}
\biggr] + \inf_{B_{\delta}} \hat{V}
\\
& \ge &  V_{*}(x) - \sup_{B_{\delta}} V_{*} +
\inf_{B_{\delta}} \hat {V}.
\end{eqnarray*}
Taking limits as $\delta\searrow0$, we have $V_{*}\le\hat{V}$.
Since $L^{\hat{v}} (V_{*}-\hat{V})\ge0$ and $V_{*}(0)=\hat{V}(0)$,
we must have $\hat{V}=V_{*}$ on $\mathbb{R}^{d}$,
and the proof of part~(b) is complete.

To prove part~(c) note that by part~(a)
we have $\varrho_{\hat{v}}<\infty$.
Therefore, $\int\tilde{h} \,\mathrm{d}\pi_{\hat{v}}\le\infty$
by Theorem~\ref{T-bst}(a), which implies
that $\int\llvert\hat{V}\rrvert \,\mathrm{d}{\mu}_{\hat{v}}\le
\infty$
by the hypothesis.
Therefore,
$\mathbb{E}^{\hat{v}}_{x}(\llvert\hat{V}(X_{t})\rrvert)$ converges
as $t\to\infty$
by \cite{ichihara-sheu}, Proposition~2.6,
which of course implies that
$\frac{1}{t} \mathbb{E}^{\hat{v}}_{x}(\llvert\hat{V}(X_{t})\rrvert)$
tends to $0$ as $t\to\infty$.
Similarly, we deduce that
$\frac{1}{t} \mathbb{E}^{v_{*}}_{x}(\llvert\hat{V}(X_{t})\rrvert)$
as $t\to\infty$.
Applying It\^o's formula to \eqref{E-HJB-hat2}, with $u\equiv v_{*}$,
we obtain $\hat{\varrho}\le\varrho_{*}$.
Another application with $u\equiv\hat{v}$ results in
$\hat{\varrho}=\varrho_{\hat{v}}$.
Therefore, $\hat{\varrho}=\varrho_{*}$.
The result then follows by part~(b).
\end{pf*}

We finish this section with the proof of Theorem~\ref{T-HJB3}.

\begin{pf*}{Proof of Theorem~\ref{T-HJB3}}
We first show that
$\lim_{\alpha\searrow0} \alpha V_{\alpha}(0)= \varrho_{*}$.
Let $\tilde\mathcal{V}(t,x) :=\mathrm{e}^{-\alpha t} \mathcal{V}(x)$,
and $\tau_{n}(t):=\tau_{n}\wedge t$.
Applying It\^o's formula to \eqref{E-key}, we obtain
\begin{eqnarray*}
\mathbb{E}^{U}_{x} \bigl[\tilde\mathcal{V} \bigl(\tau
_{n}(t),X_{\tau_{n}(t)} \bigr) \bigr] & \le& \mathcal{V}(x) -
\mathbb{E}^{U}_{x} \biggl[\int_{0}^{\tau_{n}(t)}
\alpha\tilde\mathcal{V}(s,X_{s}) \,\mathrm{d} {s} \biggr]
\\
&&{}+\mathbb{E}^{U}_{x} \biggl[\int_{0}^{\tau_{n}(t)}
\mathrm {e}^{-\alpha s} \bigl(1-h(X_{s},U_{s}) \bigr)
\mathbb{I}_{\mathcal
{H}^{c}}(X_{s},U_{s}) \,\mathrm{d} {s}
\biggr]
\\
&&{}+\mathbb{E}^{U}_{x} \biggl[\int_{0}^{\tau_{n}(t)}
\mathrm {e}^{-\alpha s} \bigl(1+r(X_{s},U_{s}) \bigr)
\mathbb{I}_{\mathcal{H}}(X_{s},U_{s}) \,\mathrm{d} {s}
\biggr].
\end{eqnarray*}
It follows that
%
\begin{eqnarray}
&& \mathbb{E}^{U}_{x} \biggl[\int_{0}^{\tau_{n}(t)}
\mathrm {e}^{-\alpha s} h(X_{s},U_{s})
\mathbb{I}_{\mathcal{H}^{c}}(X_{s},U_{s}) \,\mathrm{d} {s}
\biggr]
\nonumber
\\[-8pt]
\label{E-disc02}
\\[-8pt]
\nonumber
&&\qquad\le  \frac{1}{\alpha}+\mathcal{V}(x)+ \mathbb{E}^{U}_{x} \biggl[\int
_{0}^{\tau_{n}(t)} \mathrm {e}^{-\alpha s}
r(X_{s},U_{s}) \mathbb{I}_{\mathcal{H}}(X_{s},U_{s})
\,\mathrm{d} {s} \biggr].
\end{eqnarray}
Taking limits first as $n\nearrow\infty$ and then as $t\nearrow
\infty$
in \eqref{E-disc02},
and evaluating $U$ at an optimal $\alpha$-discounted control
$v^{*}_{\alpha}$,
relative to $r$ we obtain the estimate, using also~\eqref{E-KA3},
%
\begin{equation}
\label{E-disc03}
2k_{0}^{-1} \mathbb{E}^{v^{*}_{\alpha}}_{x}
\biggl[\int_{0}^{\infty
}\mathrm{e}^{-\alpha s}
\tilde{h} \bigl(X_{s},v^{*}_{\alpha}(X_{s})
\bigr) \,\mathrm{d} {s} \biggr] \le \frac{2}{\alpha} +\mathcal{V}(x)+2V_{\alpha}(x).
\end{equation}
By \eqref{E-KA3} and \eqref{E-disc03}, it follows that
\begin{eqnarray*}
V_{\alpha}(x) &\le &  V^{\varepsilon}_{\alpha}(x)  \le
\mathbb{E}^{v^{*}_{\alpha}}_{x} \biggl[\int_{0}^{\infty}
\mathrm{e}^{-\alpha s} r_{\varepsilon} \bigl(X_{s},v^{*}_{\alpha}(X_{s})
\bigr) \,\mathrm{d} {s} \biggr]
\\
& \le &  V_{\alpha}(x) +\varepsilon k_{0} \bigl(
\alpha^{-1} +\mathcal{V}(x) + V_{\alpha
}(x) \bigr).
\end{eqnarray*}
Multiplying by $\alpha$ and taking limits as $\alpha\searrow0$ we obtain
\[
\limsup_{\alpha\searrow0} \alpha V_{\alpha}(0) \le
\varrho_{\varepsilon} \le(1+\varepsilon k_{0}) \limsup
_{\alpha\searrow0} \alpha V_{\alpha}(0) + \varepsilon
k_{0}.
\]
The same inequalities hold for the ``$\liminf$.''
Therefore, $\lim_{\alpha\searrow0} \alpha V_{\alpha}(0)=\varrho_{*}$.

Let
\[
\tilde{V} :=\lim_{\alpha\searrow0} \bigl(V_{\alpha} -
V_{\alpha}(0) \bigr).
\]
(Note that a similar result as Lemma~\ref{L3.4} holds.)
Then $\tilde{V}$ satisfies
\[
\tilde{V}(x) \le \lim_{\delta\searrow0} \mathbb{E}^{v}_{x}
\biggl[\int_{0}^{{\breve\tau}_{\delta}} \bigl(r \bigl(X_{s},v(X_{s})
\bigr)-\varrho_{*} \bigr) \,\mathrm{d} {s} \biggr]\qquad \forall v\in\bigcup
_{\beta>0}\mathfrak{U}_{\mathrm{SM}}^{\beta}.
\]
This can be obtained without the near-monotone assumption on the
running cost;
see, for example, \cite{ari-bor-ghosh}, Lemma~3.6.9 or Lemma~3.7.8.
It follows from \eqref{E-strep} that $\tilde{V}\le V_{*}$.
On the other hand, since $L^{v_{*}}(\tilde{V} - V_{*}) \ge0$,
and $\tilde{V}(0) = V_{*}(0)$, we must have $\tilde{V}=V_{*}$
by the strong maximum principle.
\end{pf*}

\section{Approximation via spatial truncations}\label{S-truncation}

We introduce an approximation technique which is in turn used
to prove the asymptotic convergence results in Section~\ref{S-optimality}.

Let $v_{0}\in\mathfrak{U}_{\mathrm{SSM}}$ be any control such that
$\pi_{v_{0}}(r)<\infty$.
We fix the control $v_{0}$ on the complement of the ball $\bar{B}_{l}$
and leave the parameter $u$ free inside.
In other words, for each $l\in\mathbb{N}$ we define
\begin{eqnarray*}
b_{l}(x,u) & := & \cases{ \ds b(x,u), &\quad $\mbox{if }(x,u)\in
\bar{B}_{l}\times\mathbb{U}$,\vspace*{2pt}
\cr
\ds b
\bigl(x,v_{0}(x) \bigr), & \quad\mbox{otherwise},}
\\
r_{l}(x,u) & := & \cases{ \ds r(x,u), &$\quad\mbox{if } (x,u)\in
\bar{B}_{l}\times\mathbb{U}$,\vspace*{2pt}
\cr
\ds r
\bigl(x,v_{0}(x) \bigr), & $\quad\mbox{otherwise}$.}
\end{eqnarray*}
We consider the family of controlled diffusions, parameterized by $l\in
\mathbb{N}$,
given by
%
\begin{equation}
\label{E-sdeR}
\mathrm{d} {X}_{t} = b_{l}(X_{t},U_{t})
\,\mathrm{d} {t} + \sigma (X_{t}) \,\mathrm{d} {W}_{t},
\end{equation}
with associated running costs $r_{l}(x,u)$.
We denote by $\mathfrak{U}_{\mathrm{SM}}(l,v_{0})$ the subset of
$\mathfrak{U}_{\mathrm{SM}}$ consisting of those controls
$v$ which agree with $v_{0}$ on $\bar{B}_{l}^{c}$.
Let $\eta_{0}:=\pi_{v_{0}}(r)$.
It is well known that there exists a nonnegative solution
$\varphi_{0}\in\mathscr{W}_{\mathrm{loc}}^{2,p}(\mathbb{R}^{d})$,
for any $p>d$, to the Poisson equation
(see \cite{ari-bor-ghosh}, Lemma~3.7.8(ii))
\[
L^{v_{0}}\varphi_{0}(x) = \eta_{0} - \tilde{h}
\bigl(x,v_{0}(x) \bigr) \qquad x\in\mathbb{R}^{d},
\]
which is inf-compact, and satisfies, for all $\delta>0$,
\[
\varphi_{0}(x) = \mathbb{E}^{v_{0}}_{x} \biggl[
\int_{0}^{{\breve
\tau}_{\delta}} \bigl(\tilde{h} \bigl(X_{s},v_{0}
(X_{s} )\bigr)-\eta_{0} \bigr) \,\mathrm{d} {s} +
\varphi_{0}(X_{{\breve\tau}_{\delta}}) \biggr]\qquad \forall x\in\mathbb{R}^{d}.
\]

We recall the Lyapunov function $\mathcal{V}$ from
Assumption~\ref{Ass-1}.
We have the following theorem.

\begin{theorem}\label{T-trunc}
Let Assumptions~\ref{Ass-1} and \ref{Ass-2} hold. Then
for each $l\in\mathbb{N}$ there exists a solution $V^{l}$ in
$\mathscr{W}_{\mathrm{loc}}^{2,p}(\mathbb{R}^{d})$,
for any $p>d$, with $V^{l}(0)=0$,
of the HJB equation
%
\begin{equation}
\label{E-HJB-n} \min_{u\in\mathbb{U}} \bigl[L_{l}^{u}V^{l}(x)+r_{l}(x,u)
\bigr] = \varrho_{l},
\end{equation}
where $L_{l}^{u}$ is the elliptic differential operator corresponding
to the diffusion in \eqref{E-sdeR}.
Moreover, the following hold:
\begin{longlist}[(iii)]
\item[(i)]
$\varrho_{l}$ is nonincreasing in $l$;

\item[(ii)]
there exists a constant $C_{0}$, independent
of $l$, such
that $V^{l}(x) \le C_{0}+2\varphi_{0}(x)$
for all $l\in\mathbb{N}$;

\item[(iii)]
$(V^{l})^{-}\in\mathfrak{o}(\mathcal{V}+\varphi_{0})$ uniformly
over $l\in\mathbb{N}$;

\item[(iv)]
the restriction of $V^{l}$ on $B_{l}$ is in $\mathcal{C}^{2}$.
\end{longlist}
\end{theorem}

\begin{pf}
As earlier, we can show that
\[
V^{l}_{\alpha}(x) :=\inf_{U\in\mathfrak{U}}
\mathbb{E}^{U}_{x} \biggl[\int_{0}^\infty
\mathrm{e}^{-\alpha s} r_{l}(X_{s},U_{s})
\,\mathrm{d} {s} \biggr]
\]
is the minimal nonnegative solution to
%
\begin{equation}
\label{E-disHJB-n} \min_{u\in\mathbb{U}} \bigl[L_{l}^{u}V^{l}_{\alpha
}(x)+r_{l}(x,u)
\bigr] = \alpha V^{l}_{\alpha}(x),
\end{equation}
and $V^{l}_{\alpha}\in\mathscr{W}_{\mathrm{loc}}^{2,p}(\mathbb
{R}^{d})$, $p>d$.
Moreover, any measurable selector from the minimizer in
\eqref{E-disHJB-n} is an optimal control.
A similar estimate as in Lemma~\ref{L3.4}
holds and, therefore, there exists a subsequence $\{\alpha_{n}\}$,
along which
$V^{l}_{\alpha_{n}}(x)-V^{l}_{\alpha_{n}}(0)$
converges to $V^{l}$ in $\mathscr{W}_{\mathrm{loc}}^{2,p}(\mathbb
{R}^{d})$, $p>d$, and
$\alpha_{n} V^{l}_{\alpha_{n}}(0)\to\varrho_{l}$
as $\alpha_{n}\searrow0$, and $(V^{l},\varrho_{l})$ satisfies~\eqref{E-HJB-n}
(see also \cite{ari-bor-ghosh}, Lemma~3.7.8).

To show that $\pi_{v^{l}}(r)=\varrho_{l}$, $v^{l}$ is a minimizing
selector in
\eqref{E-HJB-n}, we use the following argument.
Since $\pi_{v_{0}}(r)<\infty$, we claim
that there exists a nonnegative, inf-compact
function $g\in\mathcal{C}(\mathbb{R}^{d})$ such that
$\pi_{v_{0}} (g\cdot(1+r) )<\infty$.
Indeed, this is true
since integrability and uniform integrability of a function under any given
measure are equivalent (see also the proof of \cite{ari-bor-ghosh}, Lemma~3.7.2).
Since every control in $\mathfrak{U}_{\mathrm{SM}}(l,v_{0})$ agrees
with $v_{0}$ on $B_{l}^{c}$,
then for any $x_{0}\in\bar{B}_{l}^{c}$ the map
\[
v \mapsto\mathbb{E}^{v}_{x_{0}} \biggl[\int
_{0}^{{\breve\tau}_{l}} g(X_{s}) \bigl(1+r
\bigl(X_{s},v(X_{s}) \bigr)\bigr) \,\mathrm{d} {s} \biggr]
\]
is constant on $\mathfrak{U}_{\mathrm{SM}}(l,v_{0})$.
By the equivalence of (i) and (iii) in Lemma~3.3.4 of \cite{ari-bor-ghosh},
this implies that
\[
\sup_{v\in\mathfrak{U}_{\mathrm{SM}}(l,v_{0})} \pi_{v} \bigl(g\cdot(1+r) \bigr) <
\infty\qquad\forall l\in\mathbb{N},
\]
and thus $r$ is uniformly integrable with respect to the family
$\{\pi_{v} \dvtx\break  v\in\mathfrak{U}_{\mathrm{SM}}(l,v_{0})\}$ for
any $l\in\mathbb{N}$.
It then follows by \cite{ari-bor-ghosh}, Theorem~3.7.11,  that
%
\begin{equation}
\label{T-trunc01}
\varrho_{l} = \inf_{v\in\mathfrak{U}_{\mathrm{SM}}(l,v_{0})} \pi
_{v}(r),\qquad l\in\mathbb{N}.
\end{equation}
This yields part (i).
Moreover, in view of Lemmas~\ref{L3.4} and \ref{L3.5}, we deduce that
for any $\delta>0$ it holds that
$\sup_{B_{\delta}} \llvert V^{l}\rrvert\le\kappa_{\delta}$, where
$\kappa_{\delta}$ is a constant independent of $l\in\mathbb{N}$.
It is also evident by \eqref{T-trunc01} that
$\varrho_{l}$ is decreasing in $l$ and $\varrho_{l}\le\eta_{0}$ for
all $l\in\mathbb{N}$.
Fix $\delta$ such that
$\min_{u\in\mathbb{U}} \tilde{h}(x,u) \ge2\eta_{0}$ on $B_{\delta}^{c}$.
Since $\varphi_{0}$ is nonnegative, we obtain
%
\begin{equation}
\label{T-trunc02a}
\mathbb{E}^{v_{0}}_{x} \biggl[\int
_{0}^{{\breve\tau}_{\delta}} \bigl(\tilde{h} \bigl(X_{s},v_{0}(X_{s})
\bigr)-\eta_{0} \bigr) \,\mathrm{d} {s} \biggr] \le
\varphi_{0}(x)\qquad \forall x\in\mathbb{R}^{d}.
\end{equation}
Using an analogous argument as the one used in the proof of
\cite{ari-bor-ghosh}, Lemma~3.7.8, we have
%
\begin{equation}
\label{T-trunc02}
V^{l}(x) \le\mathbb{E}^{v}_{x}
\biggl[\int_{0}^{{\breve\tau
}_{\delta}} \bigl(r_{l}
\bigl(X_{s},v(X_{s}) \bigr)-\varrho_{l} \bigr)
\,\mathrm{d} {s} \biggr] +\kappa_{\delta} \qquad\forall v\in\mathfrak{U}_{\mathrm{SM}}(l,v_{0}).
\end{equation}
Thus, by \eqref{T-trunc02a} and \eqref{T-trunc02}, and since by the choice
of $\delta>0$, it holds that
$r\le\tilde{h} \le2(\tilde{h} -\eta_{0})$ on $B_{\delta}^{c}$, we obtain
%
\begin{eqnarray}
V^{l}(x) & \le &  \mathbb{E}^{v_{0}}_{x} \biggl[\int
_{0}^{{\breve\tau}_{\delta}} 2 \bigl(\tilde{h} \bigl(X_{s},v_{0}(X_{s})
\bigr)- \eta_{0} \bigr) \,\mathrm{d} {s} \biggr] +\kappa_{\delta}
\nonumber
\\[-8pt]
\label{T-trunc03}
\\[-8pt]
\nonumber
& \le& \kappa_{\delta}+2\varphi_{0}(x) \qquad\forall x\in
\mathbb{R}^{d}.
\end{eqnarray}
This proves part~(ii).

Now\vspace*{1pt} fix $l\in\mathbb{N}$.
Let $v^{l}_{\alpha}$ be a minimizing selector of
\eqref{E-disHJB-n}.
Note then that $v^{l}_{\alpha}\in\mathfrak{U}_{\mathrm{SM}}(l,v_{0})$.
Therefore, $v^{l}_{\alpha}$ is a stable Markov control.
Let $v^{l}_{\alpha_{n}}\to v^{l}$ in the topology of Markov controls along
the same subsequence as above.
Then it is evident that $v^{l}\in\mathfrak{U}_{\mathrm{SM}}(l,v_{0})$.
Also from Lemma~\ref{L3.7}, we have
\[
\mathbb{E}^{v^{l}_{\alpha_{n}}}_{x}[{\breve\tau}_{\delta}]\, \mathop{
\longrightarrow}_{\alpha_{n}\searrow0}\, \mathbb{E}^{v^{l}}_{x}[{ \breve
\tau}_{\delta}]\qquad \forall x\in B_{\delta}^{c}, \forall
\delta>0.
\]
Using \cite{ari-bor-ghosh}, Lemma~3.7.8, we obtain the lower bound
%
\begin{equation}
\label{T-trunc04}
V^{l}(x) \ge-\varrho_{l}
\mathbb{E}^{v^{l}}_{x}[{\breve\tau }_{\delta}]-
\kappa_{\delta}.
\end{equation}
By \cite{ari-bor-ghosh}, Theorem~3.7.12(i)
(see also (3.7.50) in \cite{ari-bor-ghosh}), it holds that
%
\begin{eqnarray}
V^{l}(x)& =&  \mathbb{E}^{v^{l}}_{x}
\biggl[\int_{0}^{{\breve\tau
}_{\delta}} \bigl(r_{l}
\bigl(X_{s},v^{l}(X_{s}) \bigr)-
\varrho_{l} \bigr) \,\mathrm{d} {s} +V^{l}(X_{{\breve\tau}_{\delta}})
\biggr]
\nonumber
\\[-8pt]
\label{T-trunc05}
\\[-8pt]
\nonumber
& \ge& \mathbb{E}^{v^{l}}_{x} \biggl[\int_{0}^{{\breve\tau
}_{\delta}}
r_{l} \bigl(X_{s},v^{l}(X_{s})
\bigr) \,\mathrm{d} {s} \biggr] - \varrho_{l}\mathbb{E}^{v^{l}}_{x}[{
\breve\tau}_{\delta}] - \kappa_{\delta} \qquad\forall x\in
B_{l}^{c}.
\end{eqnarray}
By \eqref{E-KA3}, we have
\[
2k_{0}^{-1}\tilde{h}(x,u) \mathbb{I}_{\mathcal{H}}(x,u)
\le1+r(x,u) \mathbb{I}_{\mathcal{H}}(x,u).
\]
Therefore, using the preceding inequality and \eqref{T-trunc05}, we obtain
%
\begin{eqnarray}
&&V^{l}(x)+(1+\varrho_{l})
\mathbb{E}^{v^{l}}_{x}[{\breve\tau}_{\delta}]+
\kappa_{\delta}
\nonumber
\\[-8pt]
\label{T-trunc06}
\\[-8pt]
\nonumber
&& \qquad\ge \frac{2}{k_{0}} \mathbb{E}^{v^{l}}_{x} \biggl[\int
_{0}^{{\breve
\tau}_{\delta}} \tilde{h} \bigl(X_{s},v^{l}(X_{s})
\bigr) \mathbb{I}_{\mathcal
{H}} \bigl(X_{s},v^{l}(X_{s})
\bigr) \,\mathrm{d} {s} \biggr].
\end{eqnarray}
By \eqref{E-key}, \eqref{T-trunc05} and the fact that $\mathcal{V}$
is nonnegative,
we have
%
\begin{eqnarray}
&&\frac{2}{k_{0}} \mathbb{E}^{v^{l}}_{x}
\biggl[\int_{0}^{{\breve
\tau}_{\delta}} \tilde{h}
\bigl(X_{s},v^{l}(X_{s}) \bigr)
\mathbb{I}_{\mathcal
{H}^{c}} \bigl(X_{s},v^{l}(X_{s})
\bigr) \,\mathrm{d} {s} \biggr]-\mathcal{V}(x)-\mathbb{E}^{v^{l}}_{x}[{
\breve \tau}_{\delta}]
\nonumber
\\
\label{T-trunc07}
&& \qquad  \le \mathbb{E}^{v^{l}}_{x} \biggl[\int_{0}^{{\breve\tau}_{\delta}}
r \bigl(X_{s},v^{l}(X_{s}) \bigr)
\mathbb{I}_{\mathcal
{H}} \bigl(X_{s},v^{l}(X_{s})
\bigr) \,\mathrm{d} {s} \biggr]
\\
&&\qquad \le V^{l}(x)+\varrho_{l}\mathbb{E}^{v^{l}}_{x}[{
\breve\tau}_{\delta
}] + \kappa_{\delta}.
\nonumber
\end{eqnarray}
Combining \eqref{T-trunc03}, \eqref{T-trunc06} and \eqref{T-trunc07}, we obtain
\begin{eqnarray*}
\mathbb{E}^{v^{l}}_{x} \biggl[\int_{0}^{{\breve\tau}_{\delta}}
\tilde{h} \bigl(X_{s},v^{l}(X_{s}) \bigr)
\,\mathrm{d} {s} \biggr] &\le &  k_{0} (1+\varrho_{l} )
\mathbb{E}^{v^{l}}_{x}[{\breve\tau }_{\delta}]
\\
&&{}+\frac{k_{0}}{2}\mathcal{V}(x)+2k_{0} \bigl(\varphi_{0}(x)+
\kappa_{\delta} \bigr)
\end{eqnarray*}
for all $l\in\mathbb{N}$.
As earlier, using the inf-compact property of $\tilde{h}$
and the fact that $\varrho_{l}\le\eta_{0}$ is bounded, we can
choose $\delta$ large enough such that
%
\begin{equation}
\label{T-trunc08}
\hspace*{6pt}\quad\eta_{0}\mathbb{E}^{v^{l}}_{x}[{
\breve\tau}_{\delta}] \le \mathbb{E}^{v^{l}}_{x} \biggl[
\int_{0}^{{\breve\tau}_{\delta}} \tilde{h} \bigl(X_{s},v^{l}(X_{s})
\bigr) \,\mathrm{d} {s} \biggr] \le k_{0}\mathcal{V}(x)
+4k_{0} \bigl(\varphi_{0}(x)+\kappa_{\delta} \bigr)
\end{equation}
for all $l\in\mathbb{N}$.
Since $\tilde{h}$ is inf-compact, part~(iii) follows by
\eqref{T-trunc04} and \eqref{T-trunc08}.

Part (iv) is clear from regularity theory of elliptic PDE
\cite{gilbarg-trudinger}, Theorem~9.19, page~243.
\end{pf}

Similar to Theorem~\ref{T-HJB1}, we can show that oscillations of $\{
V^{l}\}$
are uniformly bounded on compacts.
Therefore, if we let $l\to\infty$
we obtain a HJB equation
%
\begin{equation}
\label{E-HJB-hat}
\min_{u\in\mathbb{U}} \bigl[L^{u}
\hat{V}(x)+r(x,u) \bigr] = \hat\varrho,
\end{equation}
with $\hat{V}\in\mathcal{C}^{2}(\mathbb{R}^{d})$ and $\lim_{l\to
\infty}\varrho_{l}=\hat\varrho$.
By Theorem~\ref{T-trunc}, we have the bound
%
\begin{equation}
\label{estim5}
\hat{V}(x) \le C_{0}+2\varphi_{0}(x),
\end{equation}
for some positive constant $C_{0}$. This of course, implies that
$\hat{V}^{+}(x) \le C_{0}+2\varphi_{0}(x)$.
Moreover, it is straightforward to show
that for any $v\in\mathfrak{U}_{\mathrm{SSM}}$ with $\varrho
_{v}<\infty$, we have
\[
\limsup_{t\to\infty} \frac{1}{t} \mathbb{E}^{v}_{x}
\bigl[\mathcal{V}(X_{t}) \bigr] < \infty.
\]
Therefore, if in addition, we have
\[
\limsup_{t\to\infty} \frac{1}{t} \mathbb{E}^{v}_{x}
\bigl[\varphi_{0}(X_{t}) \bigr] < \infty,
\]
then it follows by Theorem~\ref{T-trunc}(iii) that
%
\begin{equation}
\label{T-trunc000}
\limsup_{t\to\infty} \frac{1}{t}
\hat{V}^{-}(X_{t}) \mathop{\longrightarrow}_{t\to\infty} 0.
\end{equation}

\begin{theorem}\label{T-trunc2}
Suppose that
$\varphi_{0}\in\mathscr{O} (\min_{u\in\mathbb{U}} \tilde
{h}(\cdot,u) )$.
Then, under the assumptions of Theorem~\ref{T-trunc}, we have
$\lim_{l\to\infty}\varrho_{l}=\hat\varrho=\varrho_{*}$, and
$\hat{V}=V_{*}$.
Moreover, $V_{*}\in\mathscr{O}(\varphi_{0})$.
\end{theorem}

\begin{pf}
Let $\{\hat{v}_{l}\}$ be any sequence of measurable selectors from the
minimizer
of \eqref{E-HJB-n} and $\{\pi_{l}\}$ the corresponding sequence of
ergodic occupation measures.
Since by Theorem~\ref{T-bst} $\{\pi_{l}\}$ is tight, then by
Remark~\ref{R-tight}
if $\hat{v}$ is a limit point of a subsequence
$\{\hat{v}_{l}\}$, which we also denote by $\{\hat{v}_{l}\}$,
then $\hat{\pi}=\pi_{\hat{v}}$
is the corresponding limit point of $\{\pi_{l}\}$.
Therefore, by the lower semi-continuity of $\pi\to\pi(r)$ we have
\[
\hat\varrho= \lim_{l\to\infty} \pi_{l}(r) \ge \hat \pi(r)
= \varrho_{\hat{v}}.
\]
It also holds that
%
\begin{equation}
\label{ET4.2a} L^{\hat{v}}\hat{V}(x)+r \bigl(x,\hat{v}(x) \bigr) = \hat{
\varrho},\qquad \mbox{a.s.}
\end{equation}
By \eqref{T-trunc000}, we have
\[
\liminf_{T\to\infty} \frac{1}{T} \mathbb{E}^{\hat{v}}_{x}
\bigl[\hat{V}(X_{T}) \bigr] = 0,
\]
and hence applying It\^o's rule on \eqref{ET4.2a} we obtain
$\varrho_{\hat{v}}\le\hat\varrho$.
On the other hand, if $v_{*}$ is an optimal stationary Markov control,
then by
the hypothesis $\varphi_{0}\in\mathscr{O} (\tilde{h} )$,
the fact
that $\pi_{v_{*}} (\tilde{h})<\infty$, \eqref{estim5} and
\cite{ichihara-sheu}, Proposition~2.6,
we deduce that $\mathbb{E}^{v_{*}}_{x} [\hat{V}^{+}(X_{t}) ]$
converges as $t\to\infty$,
which of course together with \eqref{T-trunc000} implies that
$\frac{1}{t} \mathbb{E}^{\hat{v}}_{x} [\hat{V}(X_{t}) ]$
tends to $0$ as $t\to\infty$.
Therefore, evaluating \eqref{E-HJB-hat} at $v_{*}$ and
applying It\^o's rule we obtain
$\varrho_{v_{*}}\ge\hat\varrho$.
Combining the two estimates, we have
$\varrho_{\hat{v}}\le\hat\varrho\le\varrho_{*}$,
and thus equality must hold.
Here, we have used the fact that there exists an optimal Markov control for
$r$ by Theorem~\ref{T-HJB2}.

Next, we use the stochastic representation in \eqref{T-trunc05}
%
\begin{equation}
\label{ET4.2d}
\quad V^{l}(x) = \mathbb{E}^{\hat{v}_{l}}_{x}
\biggl[\int_{0}^{{\breve\tau}_{\delta}} \bigl(r \bigl(X_{s},
\hat{v}_{l}(X_{s}) \bigr)-\varrho_{l} \bigr)
\,\mathrm{d} {s} + V^{l}(X_{{\breve\tau}_{\delta}}) \biggr],\qquad x\in
B_{\delta}^{c}.
\end{equation}
Fix any $x\in B_{\delta}^{c}$.
Since $\mathfrak{U}_{\mathrm{SM}}^{\varrho_{v_{0}}}$ is compact, it
follows that for each $\delta$
and $R$ with $0<\delta<R$, the map
$F_{\delta,R}(v)\dvtx\mathfrak{U}_{\mathrm{SM}}^{\varrho
_{v_{0}}}\to\mathbb{R}_{+}$ defined by
\[
F_{\delta,R}(v) := \mathbb{E}^{v}_{x} \biggl[\int
_{0}^{{\breve\tau}_{\delta}\wedge\tau_{R}} r \bigl(X_{s},v(X_{s})
\bigr) \,\mathrm{d} {s} \biggr]
\]
is continuous.
Therefore, the map $\bar{F}_{\delta}:=\lim_{R\nearrow\infty}
F_{\delta,R}$
is lower semi-continuous.
It follows that
%
\begin{equation}
\label{E4.17nu}
\mathbb{E}^{\hat{v}}_{x} \biggl[\int
_{0}^{{\breve\tau}_{\delta}}r \bigl(X_{s},\hat
{v}(X_{s}) \bigr) \,\mathrm{d} {s} \biggr] \le \lim
_{l\to\infty} \mathbb{E}^{\hat{v}_{l}}_{x} \biggl[\int
_{0}^{{\breve\tau}_{\delta}} r \bigl(X_{s},
\hat{v}_{l}(X_{s}) \bigr) \,\mathrm{d} {s} \biggr].
\end{equation}
On the other hand, since $\tilde{h}$ is inf-compact, it follows by
\eqref{T-trunc08} that ${\breve\tau}_{\delta}$ is uniformly
integrable with respect
to the measures $ \{\mathbb{P}^{\hat{v}_{l}}_{x} \}$.
Therefore, as also shown in Lemma~\ref{L3.7}, we have
%
\begin{equation}
\label{ET4.2f} \lim_{l\to\infty} \mathbb{E}^{\hat{v}_{l}}_{x}
[{\breve \tau}_{\delta} ] = \mathbb{E}^{\hat{v}}_{x} [{
\breve\tau}_{\delta} ].
\end{equation}
Since $V^{l}\to\hat{V}$, uniformly
on compact sets, and $\varrho_{l}\to\varrho_{*}$, as $l\to\infty$,
it follows by \eqref{ET4.2d}--\eqref{ET4.2f} that
\[
\hat{V}(x) \ge\mathbb{E}^{\hat{v}}_{x} \biggl[\int
_{0}^{{\breve\tau}_{\delta}} \bigl(r \bigl(X_{s},
\hat{v}(X_{s}) \bigr)-\varrho_{*} \bigr) \,\mathrm{d} {s} +
\hat{V}(X_{{\breve\tau}_{\delta}}) \biggr],\qquad x\in B_{\delta
}^{c}.
\]
Therefore, by Theorem~\ref{T-HJB2}(b),
for any $\delta>0$ and $x\in B^{c}_{\delta}$
we obtain
\begin{eqnarray*}
V_{*}(x) & \le &\mathbb{E}^{\hat{v}}_{x} \biggl[\int
_{0}^{{\breve\tau}_{\delta}} \bigl(r \bigl(X_{s},
\hat{v}(X_{s}) \bigr)-\varrho_{*} \bigr) \,\mathrm{d} {s} +
V_{*}(X_{{\breve\tau}_{\delta}}) \biggr]
\\
& \le & \hat{V}(x) + \mathbb{E}^{\hat{v}}_{x} \bigl[V^{*}(X_{{\breve
\tau}_{\delta}})
\bigr] - \mathbb{E}^{\hat{v}}_{x} \bigl[\hat{V}(X_{{\breve\tau}_{\delta
}})
\bigr],
\end{eqnarray*}
and taking limits as $\delta\searrow0$, using the fact that
$\hat{V}(0)=V_{*}(0)=0$,
we obtain $V_{*}\le\hat{V}$ on $\mathbb{R}^{d}$.
Since $L^{\hat{v}}(V_{*}-\hat{V})\ge0$, we must have $V_{*}= \hat{V}$.
By Theorem~\ref{T-trunc}(ii), we have $V_{*}\in\mathscr{O}(\varphi_{0})$.
\end{pf}

\begin{remark}
It can be seen from the proof of Theorem~\ref{T-trunc2}
that the assumption $\varphi_{0}\in\mathscr{O} (\tilde{h} )$
can be replaced by the weaker hypothesis that\break 
$\frac{1}{T} \mathbb{E}^{v_{*}}_{x} [\varphi_{0}(X_{T})
]\to0$ as $T\to\infty$.
\end{remark}

\begin{remark}\label{rem-trunc}
It is easy to see that if one replaces $r_{l}$ by
\[
r_{l}(x,u) = %
\cases{\ds r(x,u)+\frac{1}{l}f(u), & \quad $\mbox{for }x\in\bar{B}_{l}$,
\vspace*{2pt}\cr
\ds r \bigl(x,v_{0}(x)
\bigr)+\frac{1}{l}f \bigl(v_{0}(x) \bigr), &\quad $\mbox{otherwise}$,}
\]
for some positive valued continuous function $f$, the same conclusion of
Theorem~\ref{T-trunc2} holds.
\end{remark}

If we consider the controlled dynamics given by~\eqref{eg-sde2}, with running cost as in~\eqref{eg-cost},
then there exists a function $\mathcal{V}\sim\llvert x\rrvert^{m}$
satisfying \eqref{E-KA1}.
This fact is proved in Proposition~\ref{eg-prop}.
There also exists a Lyapunov function $\mathcal{V}_{0}\in\mathscr
{O} (\llvert x\rrvert^{m} )$,
satisfying the assumption in Theorem~\ref{T-trunc2},
relative to any control $v_{0}$ with
$\pi_{v_{0}} (\tilde{h} )<\infty$,
where $\tilde{h}$ is selected as in Remark~\ref{R3.4}.
Indeed, in order to construct $\mathcal{V}_{0}$ we
recall the function $\psi$ in \eqref{dieker-gao}.
Let $\mathcal{V}_{0}\in\mathcal{C}^{2}(\mathbb{R}^{d})$ be any function
such that $\mathcal{V}_{0}=\psi^{{m}/{2}}$ on the complement
of the unit ball centered at the origin.
Observe that for some positive constants $\kappa_{1}$ and $\kappa_{2}$
it holds that
\[
\kappa_{1}\llvert x\rrvert^{2} \le\psi(x) \le
\kappa_{2}\llvert x\rrvert^{2}.
\]
Then a straightforward calculation from \eqref{dieker-gao-ly}
shows that \eqref{E-lyap} holds with the above choice of $\mathcal{V}_{0}$.
By the stochastic representation of $\varphi_{0}$, it follows that
$\varphi_{0}\in\mathscr{O}(\mathcal{V}_{0})$.
We have proved the following corollary.

\begin{corollary}\label{C-exist}
For the queueing diffusion model with controlled dynamics given by
\eqref{eg-sde2}, and running cost given by \eqref{eg-cost},
there exists a solution (up to
an additive constant) to
the associated HJB in the class of functions in $\mathcal
{C}^{2}(\mathbb{R}^{d})$
whose positive part grows no faster
than $\llvert x\rrvert^{m}$ and whose negative part is in
$\mathfrak{o} (\llvert x\rrvert^{m} )$.
\end{corollary}

We conclude this section with the following remark.
%
\begin{remark}
Comparing the approximation technique introduced in this section with that
in Section~\ref{S-ergodic}, we see that the spatial
truncation technique relies on more restrictive assumption on the
Lyapunov function $\mathcal{V}_{0}$
and the running cost function (Theorem~\ref{T-trunc2}).
In fact, the growth of $\tilde{h}$ also restricts the growth of $r$ by~\eqref{E-KA3}.
Therefore, the class of ergodic diffusion control problems considered
in this
section is more restrictive. For example, if the running cost $r$
satisfies \eqref{eg-cost} and $\tilde{h}\sim\llvert x\rrvert^{m}$,
then it is not obvious
that one can obtain a Lyapunov function $\mathcal{V}_{0}$ with growth at
most of order $\llvert x\rrvert^{m}$.
For instance, if the drift has strictly sub-linear
growth, then it is expected that the Lyapunov function
should have growth larger than $\llvert x\rrvert^{m}$.
Therefore, the class of problems
considered in Section~\ref{S-ergodic} is larger than those
considered in this section.
\end{remark}

\section{Asymptotic convergence}\label{S-optimality}

In this section, we prove that the value of the ergodic control problem
corresponding to the multi-class $M/M/N+M$ queueing network
asymptotically converges to $\varrho_{*}$, the value of the ergodic
control for
the controlled diffusion.

Recall the diffusion-scaled processes $\hat{X}^{n}$, $\hat{Q}^{n}$ and
$\hat{Z}^{n}$ defined in \eqref{dc1}, and from \eqref{dc2} and
\eqref{dc3} that
%
\begin{eqnarray}
\hat{X}^{n}_{i}(t) &=&  \hat{X}^{n}_{i}(0)
+ \ell^{n}_{i} t - \mu_{i}^{n} \int
_{0}^{t} \hat{Z}_{i}^{n}(s)
\,\mathrm{d} {s} - \gamma_{i}^{n} \int_{0}^{t}
\hat{Q}^{n}_{i}(s) \,\mathrm{d} {s}
\nonumber
\\[-8pt]
\label{asym-1}
\\[-8pt]
\nonumber
&&{}+ \hat{M}_{A,i}^{n}(t) - \hat{M}_{S,i}^{n}(t)
- \hat{M}_{R,i}^{n}(t),
\end{eqnarray}
where\vspace*{1pt} $\hat{M}_{A,i}^{n}(t)$, $\hat{M}_{S,i}^{n}(t)$ and
$\hat{M}_{R,i}^{n}(t)$, $i=1,\ldots,d$, as defined in \eqref{dc3},
are square integrable martingales w.r.t. the filtration
$\{\mathcal{F}_{t}^{n}\}$ with quadratic variations
\begin{eqnarray*}
\bigl\langle\hat{M}_{A,i}^{n} \bigr\rangle(t) & =&
\frac{\lambda^{n}_{i}}{n} t,
\\
\bigl\langle\hat{M}_{S,i}^{n} \bigr\rangle(t) & =&
\frac{\mu^{n}_{i}}{n}\int_{0}^{t}Z^{n}_{i}(s)
\,\mathrm{d} {s},
\\
\bigl\langle\hat{M}_{R,i}^{n} \bigr\rangle(t) & =&
\frac{\gamma^{n}_{i}}{n} \int_{0}^{t}Q^{n}_{i}(s)
\,\mathrm{d} {s}.
\end{eqnarray*}

\subsection{The lower bound}

In this section, we prove Theorem~\ref{T-lowerbound}.

\begin{pf*}{Proof of Theorem~\ref{T-lowerbound}}
Recall the definition of $\hat{V}^{n}$ in \eqref{E-Vn}, and
consider a sequence such that $\sup_{n} \hat{V}^{n}(\hat
{X}^{n}(0))<\infty$.
Let $\varphi\in\mathcal{C}^{2}(\mathbb{R}^{d})$ be
any function satisfying $\varphi(x):=\llvert x\rrvert^{m}$ for $\llvert
x\rrvert\ge1$.
As defined in Section~\ref{S-notation},
$\Delta X(t)$ denotes the jump of the process $X$ at time $t$.
Applying It\^o's formula on $\varphi$ (see, e.g.,
\cite{kallenberg}, Theorem~26.7), we obtain from \eqref{asym-1} that
\begin{eqnarray*}
\mathbb{E} \bigl[\varphi \bigl(\hat{X}^{n}_{1}(t) \bigr)
\bigr] &=&  \mathbb{E} \bigl[\varphi \bigl(\hat{X}^{n}_{1}(0)
\bigr) \bigr] + \mathbb{E} \biggl[\int_{0}^{t}
\Theta^{n}_{1} \bigl(\hat{X}^{n}_{1}(s),
\hat{Z}^{n}_{1}(s) \bigr) \varphi' \bigl(
\hat{X}^{n}_{1}(s) \bigr) \,\mathrm{d} {s} \biggr]
\\
&&{}+\mathbb{E} \biggl[\int_{0}^{t}
\Theta^{n}_{2} \bigl(\hat{X}^{n}_{1}(s),
\hat{Z}^{n}_{1}(s) \bigr) \varphi''
\bigl(\hat{X}^{n}_{1}(s) \bigr) \,\mathrm{d} {s} \biggr]
\\
&&{}+\mathbb{E}\sum_{s\le t} \biggl(\Delta\varphi \bigl(
\hat {X}^{n}_{1}(s) \bigr) -\varphi' \bigl(
\hat{X}^{n}_{1}(s-) \bigr)\cdot\Delta\hat
{X}^{n}_{1}(s)
\\
&&\hspace*{56pt}{}-\frac{1}{2}\varphi'' \bigl(
\hat{X}^{n}(s-) \bigr)\Delta\hat{X}^{n}_{1}(s)
\Delta\hat{X}^{n}_{1}(s) \biggr),
\end{eqnarray*}
where
\begin{eqnarray*}
\Theta^{n}_{1}(x,z)& := & \ell^{n}_{1} -
\mu^{n}_{1} z-\gamma^{n}_{1}(x-z),
\\
\Theta^{n}_{2}(x,z)& :=& \frac{1}{2} \biggl(
\mu^{n}_{1}\rho_{1} + \frac{\lambda^{n}_{1}}{n} +
\frac{\mu^{n}_{1} z + \gamma^{n}_{1}(x-z)}{\sqrt{n}} \biggr).
\end{eqnarray*}
Since $\{\ell^{n}_{1}\}$ is a bounded sequence,
it is easy to show that for all $n$ there exist
positive constants $\kappa_{i}$, $i=1,2$,
independent of $n$, such that
\begin{eqnarray*}
\Theta_{1}^{n}(x,z) \varphi'(x)& \le&
\kappa_{1} \bigl(1+ \bigl\llvert(e\cdot x)^{+} \bigr
\rrvert^{m} \bigr)- \kappa _{2}\llvert x
\rrvert^{m},
\\
\Theta^{n}_{2}(x,z) \varphi''(x)
& \le&  \kappa_{1} \bigl(1+ \bigl\llvert(e\cdot x)^{+} \bigr
\rrvert^{m} \bigr)+\frac
{\kappa_{2}}{4}\llvert x\rrvert^{m},
\end{eqnarray*}
provided that $x-z \le(e\cdot x)^{+}$ and $\frac{z}{\sqrt{n}}\le1 $.
We next compute the terms corresponding to the jumps.
For that, first we see that the jump size is of order $\frac{1}{\sqrt{n}}$.
We can also find a positive constant $\kappa_{3}$ such that
\[
\sup_{|y-x|\le1} \bigl\llvert\varphi''(y)
\bigr\rrvert \le\kappa_{3} \bigl(1+\llvert x\rrvert^{m-2}
\bigr)\qquad \forall x\in\mathbb{R}^{d}.
\]
Using Taylor's approximation, we obtain the inequality
\[
\Delta\varphi \bigl(\hat{X}^{n}_{1}(s) \bigr)-
\varphi' \bigl(\hat {X}^{n}_{1}(s-) \bigr) \cdot
\Delta\hat{X}^{n}_{1}(s) \le \frac{1}{2}\sup
_{|y-\hat{X}^{n}_{1}(s-)|\le1} \bigl\llvert\varphi ''(y)
\bigr\rrvert \bigl[\Delta \bigl(\hat{X}^{n}_{1}(s) \bigr)
\bigr]^{2}.
\]
Hence, combining the above facts we obtain
%
\begin{eqnarray}\label{E5.2}
&&\mathbb{E}\sum_{s\le t} \biggl(\Delta
\varphi \bigl(\hat{X}^{n}_{1}(s) \bigr) -
\varphi' \bigl(\hat{X}^{n}_{1}(s-) \bigr) \cdot
\Delta\hat{X}^{n}_{1}(s)\nonumber
\\
&&\quad{}-\frac{1}{2}\varphi'' \bigl(
\hat{X}^{n}_{1}(s-) \bigr) \Delta\hat{X}^{n}_{1}(s)
\Delta\hat{X}^{n}_{1}(s) \biggr)
\nonumber
\\
&&\qquad \le\mathbb{E}\sum_{s\le t}\kappa_{3}
\bigl(1+ \bigl\llvert\hat {X}^{n}_{1}(s-) \bigr
\rrvert^{m-2} \bigr) \bigl(\Delta \bigl(\hat{X}^{n}_{1}(s)
\bigr) \bigr)^{2}
\\
&&\qquad = \kappa_{3}\mathbb{E} \biggl[\int_{0}^{t}
\bigl(1+ \bigl\llvert\hat {X}^{n}_{1}(s) \bigr
\rrvert^{m-2} \bigr) \biggl(\frac{\lambda^{n}_{1}}{n}+\frac{\mu^{n}_{1}Z_{1}^{n}(s)}{n} +
\frac{\gamma^{n}_{1} Q^n_1(s) }{n} \biggr) \,\mathrm{d} {s} \biggr]
\nonumber\\
&&\qquad \le\mathbb{E} \biggl[\int_{0}^{t} \biggl(
\kappa_{4}+\frac{\kappa_{2}}{4} \bigl\llvert\hat{X}^{n}_{1}(s)
\bigr\rrvert^{m} +\kappa_{5} \bigl( \bigl(e\cdot
\hat{X}^{n}(s) \bigr)^{+} \bigr)^{m} \biggr)
\,\mathrm{d} {s} \biggr],
\nonumber
\end{eqnarray}
for some suitable positive
constants $\kappa_{4}$ and $\kappa_{5}$, independent of $n$,
where in the second inequality we use the fact that the optional martingale
$[\hat{X}^{n}_{1}]$
is the sum of the squares of the jumps, and that
$[\hat{X}^{n}_{1}]-\langle\hat{X}^{n}_{1}\rangle$ is a martingale.
Therefore, for some positive constants $C_{1}$ and $C_{2}$ it holds that
%
\begin{eqnarray}
0 &\le & \mathbb{E} \bigl[\varphi \bigl(\hat{X}^{n}_{1}(t)
\bigr) \bigr]
\nonumber\\
\label{rev01}& \le&  \mathbb{E} \bigl[\varphi \bigl(\hat{X}^{n}_{1}(0)
\bigr) \bigr]+C_{1}t-\frac{\kappa_{2}}{2} \mathbb{E} \biggl[\int
_{0}^{t} \bigl\llvert\hat{X}^{n}_{1}(s)
\bigr\rrvert^{m} \,\mathrm{d} {s} \biggr]
\\
&&{}+ C_{2}\mathbb{E} \biggl[\int_{0}^{t}
\bigl( \bigl(e\cdot\hat{X}^{n}(s) \bigr)^{+}
\bigr)^{m} \,\mathrm{d} {s} \biggr].
\nonumber
\end{eqnarray}
By (\ref{cost1}), we have
\[
r \bigl(\hat{Q}^{n}(s) \bigr)\ge\frac{c_{1}}{d^{m}} \bigl( \bigl(e\cdot
\hat{X}^{n}(s) \bigr)^{+} \bigr)^{m},
\]
which,
combined with the assumption that $\sup_{n} \hat{V}^{n}(\hat
{X}^{n}(0))<\infty$,
implies that
\[
\sup_{n} \limsup_{T\to\infty} \frac{1}{T}
\mathbb{E} \biggl[\int_{0}^{T} \bigl( \bigl(e\cdot
\hat{X}^{n}(s) \bigr)^{+} \bigr)^{m} \,\mathrm{d}
{s} \biggr] < \infty.
\]
In turn, from \eqref{rev01} we obtain
\[
\sup_{n} \limsup_{T\to\infty} \frac{1}{T}
\mathbb{E} \biggl[\int_{0}^{T} \bigl\llvert
\hat{X}^{n}_{1}(s) \bigr\rrvert^{m} \,\mathrm{d} {s}
\biggr] < \infty.
\]
Repeating the same argument for coordinates $i=2,\ldots,d$, we obtain
%
\begin{equation}
\label{E5.4}
\sup_{n} \limsup_{T\to\infty}
\frac{1}{T} \mathbb{E} \biggl[\int_{0}^{T}
\bigl\llvert\hat{X}^{n}(s) \bigr\rrvert^{m} \,\mathrm{d} {s}
\biggr] < \infty.
\end{equation}
We introduce the process
\[
U^{n}_{i}(t) := %
\cases{ \ds\frac{\hat{X}^{n}_{i}(t)-\hat{Z}^{n}_{i}(t)}{(e\cdot\hat
{X}^{n}(t))^{+}},\qquad
i=1,\ldots,d, & \quad$\mbox{if } \bigl(e\cdot\hat{X}^{n}(t)
\bigr)^{+}>0$, \vspace*{3pt}
\cr
e_{d}, & $\quad
\mbox{otherwise}$.}
\]
Since $Z^{n}$ is work-conserving,
it follows that $U^{n}$ takes values in $\mathcal{S}$,
and $U^n_i(t)$ represents the fraction of class $i$ customers in queue.
Define the mean empirical measures
\[
\Phi^{n}_{T}(A\times B) := \frac{1}{T} \mathbb{E}
\biggl[\int_{0}^{T}\mathbb{I}_{A\times B}
\bigl(\hat{X}^{n}(s),U^{n}(s) \bigr) \,\mathrm{d} {s} \biggr]
\]
for Borel sets $A\subset\mathbb{R}^{d}$ and $B\subset\mathcal{S}$.

From \eqref{E5.4}, we see that the family $\{\Phi^{n}_{T} \dvtx T>0,
n\ge1\}$
is tight.
Hence, for any sequence $T_{k}\to\infty$, there exists a subsequence,
also denoted by $T_{k}$, such that
$\Phi^{n}_{T_{k}}\to\pi^{n}$, as $k\to\infty$.
It is evident that $\{\pi^{n} \dvtx n\ge1\}$ is tight.
Let $\pi^{n}\to\pi$ along some subsequence, with
$\pi\in\mathcal{P}(\mathbb{R}^{d}\times\mathcal{S})$.
Therefore, it is not hard to show that
\[
\lim_{n\to\infty} \hat{V}^{n} \bigl(\hat{X}^{n}(0)
\bigr) \ge\int_{\mathbb{R}^{d}\times\mathbb{U}} \tilde{r}(x,u) \pi (\mathrm{d} {x},
\mathrm{d} {u}),
\]
where, as defined earlier, $\tilde{r}(x,u)=r((e\cdot x)^{+}u)$.
To complete the proof of the
theorem, we only need to show that $\pi$ is an ergodic
occupation measure for the diffusion.
For that, consider $f\in\mathcal{C}^{\infty}_{c}(\mathbb{R}^{d})$.
Recall that $[\hat{X}^{n}_{i},\hat{X}^{n}_j]=0$
for $i\neq j$ \cite{pang-talreja-whitt}, Lemmas~9.2 and~9.3.
Therefore, using It\^o's formula and the definition of $\Phi^{n}_{T}$, we
obtain
%
\begin{eqnarray}
&& \frac{1}{T} \mathbb{E} \bigl[f \bigl(\hat{X}^{n}(T)
\bigr) \bigr]\nonumber\\
&&\qquad= \frac{1}{T} \mathbb{E} \bigl[f \bigl(\hat{X}^{n}(0)
\bigr) \bigr]
\nonumber\\
\label{asym-2}
&&\qquad\quad{}+ \int_{\mathbb{R}^{d}\times\mathbb{U}} \Biggl(\sum_{i=1}^{d}
\mathcal{A}^{n}_{i}(x,u)\cdot f_{x_{i}}(x) +
\mathcal{B}^{n}_{i}(x,u)f_{x_{i}x_{i}}(x) \Biggr) \Phi
^{n}_{T}(\mathrm{d} {x},\mathrm{d} {u})
\\
\nonumber
&&\qquad\quad{}+\frac{1}{T} \mathbb{E}\sum_{s\le T} \Biggl[
\Delta f \bigl(\hat{X}^{n}(s) \bigr)-\sum_{i=1}^{d}
f_{x_{i}} \bigl(\hat{X}^{n}(s-) \bigr)\cdot\Delta\hat
{X}^{n}_{i}(s)
\\
&&\qquad\quad\hspace*{54pt}{}- \frac{1}{2}\sum_{i,j=1}^{d}
f_{x_{i}x_{j}} \bigl(\hat {X}^{n}(s-) \bigr) \Delta
\hat{X}^{n}_{i}(s)\Delta\hat{X}^{n}_{j}(s)
\Biggr],\nonumber
\end{eqnarray}
where
\begin{eqnarray*}
\mathcal{A}^{n}_{i}(x,u)& :=& \ell^{n}_{i}
-\mu^{n}_{i} \bigl(x_{i}-(e\cdot
x)^{+}u_{i} \bigr)-\gamma^{n}_{i} (e
\cdot x)^{+}u_{i},
\\
\mathcal{B}^{n}_{i}(x,u)& := & \frac{1}{2} \biggl(
\mu^{n}_{i}\rho_{i}+ \frac{\lambda^{n}_{i}}{n} +
\frac{\mu^{n}_{i}x_{i}+(\gamma^{n}_{i}-\mu^{n}_{i})(e\cdot
x)^{+}u_{i}}{\sqrt{n}} \biggr).
\end{eqnarray*}
We first bound the last term in \eqref{asym-2}.
Using Taylor's formula, we see that
\begin{eqnarray*}
&&\Delta f \bigl(\hat{X}^{n}(s) \bigr)-\sum
_{i=1}^{d} \nabla f \bigl(\hat{X}^{n}(s-)
\bigr)\cdot\Delta\hat{X}^{n}(s)
\\
&&\quad{}
-\frac{1}{2}\sum_{i,j=1}^{d}
f_{x_{i}x_{j}} \bigl(\hat {X}^{n}(s-) \bigr) \Delta
\hat{X}^{n}_{i}(s)\Delta\hat{X}^{n}_{j}(s)
\\
&&\qquad= \frac{k\Vert f\Vert_{\mathcal{C}^{3}}}{\sqrt{n}}\sum_{i,j=1}^{d} \bigl
\llvert\Delta\hat{X}^{n}_{i}(s)\Delta\hat{X}^{n}_{j}(s)
\bigr\rrvert
\end{eqnarray*}
for some positive constant $k$, where we use the fact that the jump
size is
$\frac{1}{\sqrt{n}}$. Hence, using the fact that independent Poisson processes
do not have simultaneous jumps w.p.1, using the identity
$\hat{Q}^{n}_{i}=\hat{X}^{n}_{i}-\hat{Z}^{n}_{i}$, we obtain
%
\begin{eqnarray}
&&\frac{1}{T} \mathbb{E}\sum_{s\le T}
\Biggl[\Delta f \bigl(\hat{X}^{n}(s) \bigr) -\sum
_{i=1}^{d}\nabla f \bigl(\hat{X}^{n}(s-)
\bigr)\cdot\Delta\hat{X}^{n}(s)
\nonumber\\
\label{E5.6}
&&\hspace*{33pt}\quad{}-\frac{1}{2}\sum_{i,j=1}^{d}
f_{x_{i}x_{j}} \bigl(\hat{X}^{n}(s-) \bigr) \Delta
\hat{X}^{n}_{i}(s)\Delta\hat{X}^{n}_{j}(s)
\Biggr]
\\
&&\qquad\le\frac{k \Vert f\Vert_{\mathcal{C}^{3}}}{T\sqrt{n}} \mathbb{E} \Biggl[\int_{0}^{T}
\sum_{i=1}^{d} \biggl(\frac{\lambda
^{n}_{i}}{n} +
\frac{\mu_{i}^{n} Z^{n}_{i}(s)}{n} +\frac{\gamma_{i}^{n}Q^{n}_{i}(s)}{n} \biggr) \,\mathrm{d} {s} \Biggr].\nonumber
\end{eqnarray}
Therefore, first letting $T\to\infty$ and using \eqref{E5.2} and
\eqref{E5.4}
we see that the expectation on
the right-hand side of \eqref{E5.6} is bounded above.
Therefore, as $n\to\infty$, the left-hand side of \eqref{E5.6} tends
to $0$.
Thus, by \eqref{asym-2} and the
fact that $f$ is compactly supported, we obtain
\[
\int_{\mathbb{R}^{d}\times\mathbb{U}} L^{u} f(x) \pi(\mathrm {d} {x},
\mathrm{d} {u}) = 0,
\]
where
\[
L^{u} f(x) = \lambda_{i} \partial_{ii} f(x) +
\bigl(\ell_{i}-\mu_{i} \bigl(x_{i}-(e\cdot
x)^{+}u_{i} \bigr)-\gamma_{i} (e\cdot
x)^{+}u_{i} \bigr) \partial_{i} f(x).
\]
Therefore, $\pi\in\mathscr{G}$.
\end{pf*}

\subsection{The upper bound}

The proof of the upper bound in Theorem~\ref{T-upperbound}
is a little more involved than that of the lower bound.
Generally, it is very helpful if one has uniform stability across $n\in
\mathbb{N}$
(see, e.g., \cite{budhi-ghosh-lee}). In \cite{budhi-ghosh-lee}, uniform
stability is obtained from the reflected dynamics with the Skorohod mapping.
However, here we establish the asymptotic upper bound by using the
technique of spatial truncation that we have introduced in
Section~\ref{S-truncation}.
Let $v_{\delta}$ be any precise continuous control in $\mathfrak
{U}_{\mathrm{SSM}}$ satisfying
$v_{\delta}(x)=u_{0}=(0,\ldots,0,1)$ for $\llvert x\rrvert>K>1$.

First, we construct a work-conserving admissible policy for each $n\in
\mathbb{N}$
(see \cite{atar-mandel-rei}).
Define a measurable map
$\varpi\dvtx\{z\in\mathbb{R}^{d}_{+} \dvtx e\cdot z\in\mathbb
{Z}\}
\to\mathbb{Z}^{d}_{+}$
as follows: for $z=(z_{1},\ldots,z_{d})\in\mathbb{R}^{d}$, let
\[
\varpi(z) := \Biggl(\lfloor z_{1}\rfloor,\ldots,\lfloor
z_{d-1}\rfloor, \lfloor z_{d}\rfloor+\sum
_{i=1}^{d} \bigl(z_{i}-\lfloor
z_{i}\rfloor \bigr) \Biggr).
\]
Note that $|\varpi(z)-z|\le2d$.
Define
\begin{eqnarray*}
u_{h}(x)& :=& \varpi \bigl((e\cdot x-n)^{+}v_{\delta}
\bigl(\hat {x}^{n} \bigr) \bigr),\qquad x\in\mathbb{R}^{d},
\\
\hat{x}^{n}& := &\biggl(\frac{x_{1}-\rho_{1}n}{\sqrt{n}},\ldots, \frac{x_{d}-\rho_{d}n}{\sqrt{n}}
\biggr),
\\
A_{n}& :=&  \Bigl\{x\in\mathbb{R}^{d}_{+} \dvtx
\sup_{i} |x_{i}-\rho_{i} n|\le K\sqrt{n}
\Bigr\}.
\end{eqnarray*}
We define a state-dependent, work-conserving policy as follows:
%
\begin{equation}
\label{E-zpolicy}
Z_{i}^{n} \bigl[X^{n} \bigr] :=
\cases{\ds X_{i}^{n}-u_{h}
\bigl(X^{n} \bigr),& \quad$\mbox{if }X^{n}\in
A_{n}$,\vspace*{3pt}
\cr
\ds X_{i}^{n}\wedge
\Biggl(n-\sum_{j=1}^{i-1}X_j^{n}
\Biggr)^{+},&\quad$\mbox{otherwise}$.}
\end{equation}
Therefore, whenever the state of the system is in $A_{n}^{c}$,
the system works under the fixed priority policy with the least priority
given to class-$d$ jobs.
First, we show that this is a well-defined policy for all large $n$.
It is enough to show that
$X_{i}^{n}-u_{h}(X^{n})\ge0$ for all $i$ when $X^{n}\in A_{n}$.
If not, then for some $i$, $1\le i\le d$, we must have
$X_{i}^{n}-u_{h}(X^{n})<0$
and so $X_{i}^{n}< (e\cdot X^{n}-n)^{+} +d$.
Since $X^{n}\in A_{n}$, we obtain
\begin{eqnarray*}
-K\sqrt{n}+\rho_{i} n& \le&  X_{i}^{n}
\\
& < & \bigl(e\cdot X^{n}-n \bigr)^{+} +d
\\
& =&  \Biggl(\sum_{i=1}^{d}
\bigl(X_{i}^{n}-\rho_{i}n \bigr)
\Biggr)^{+} +d
\\
& \le &  dK\sqrt{n}+d.
\end{eqnarray*}
But this cannot hold for large $n$.
Hence, this policy is well defined for all large $n$.
Under the policy defined in \eqref{E-zpolicy}, $X^{n}$ is
a Markov process and its generator given by
\begin{eqnarray*}
\mathcal{L}_{n}f(x) &=&  \sum_{i=1}^{d}
\lambda_{i}^{n} \bigl(f(x+e_{i})-f(x) \bigr) +
\sum_{i=1}^{d} \mu_{i}^{n}Z_{i}^{n}[x]
\bigl(f(x-e_{i})-f(x) \bigr)
\\
&&{}+\sum_{i=1}^{d} \gamma^{n}_{i}
Q^{n}_{i}[x] \bigl(f(x-e_{i})-f(x) \bigr), \qquad x\in
\mathbb{Z}^{d}_{+},
\end{eqnarray*}
where $Z^{n}[x]$ is as above and $Q^{n}[x] :=x-Z^{n}[x]$.
It is easy to see that, for $x\notin A_{n}$,
\[
Q_{i}^{n}[x] = \Biggl[x_{i}- \Biggl(n-\sum
_{j=1}^{i-1} x_j
\Biggr)^{+} \Biggr]^{+}.
\]

\begin{lemma}\label{lem-uni}
Let $X^{n}$ be the Markov process corresponding to the above control.
Let $q$ be an even positive integer.
Then there exists $n_{0}\in\mathbb{N}$ such that
\[
\sup_{n\ge n_{0}} \limsup_{T\to\infty} \frac{1}{T}
\mathbb{E} \biggl[\int_{0}^{T} \bigl\llvert\hat
{X}^{n}(s) \bigr\rrvert^{q} \,\mathrm{d} {s} \biggr] < \infty,
\]
where $\hat{X}^{n}= (\hat{X}^{n}_{1},\ldots,\hat
{X}^{n}_{d} )^{\mathsf{T}}$
is the diffusion-scaled process corresponding to the process $X^{n}$,
as defined in \eqref{dc1}.
\end{lemma}

\begin{pf}
The proof technique is inspired by \cite{atar-giat-shim-2}, Lemma~3.1.
Define
\[
f_{n}(x) :=\sum_{i=1}^{d}
\beta_{i} (x_{i}-\rho_{i} n )^{q},
\]
where $\beta_{i}$, $i=1,\ldots,d$, are positive constants to be
determined later.
We first
show that for a suitable choice of $\beta_{i}$, $i=1,\ldots,d$,
there exist constants $C_{i}$,
$i=1,2$, independent of $n\ge n_{0}$, such that
%
\begin{equation}
\label{uni-06}
\mathcal{L}_{n} f_{n}(x) \le
C_{1} n^{{q}/{2}}-C_{2} f_{n}(x), \qquad x\in
\mathbb{Z}^{d}_{+}.
\end{equation}
Choose $n$ large enough so that the policy is well defined.
We define $Y^{n}_{i}:=x_{i}-\rho_{i} n$. Note that
\[
(a\pm1)^{q}-a^{q} = \pm q a\cdot a^{q-2} +
\mathscr{O} \bigl(a^{q-2} \bigr),\qquad a\in \mathbb{R}.
\]
Also, $\mu_{i}^{n} Z^{n}_{i}[x] = \mu^{n}_{i} x_{i}-\mu^{n}_{i}
Q_{i}^{n}[x]$.
Then
%
\begin{eqnarray}
\mathcal{L}_{n} f_{n}(x)& =&  \sum
_{i=1}^{d}\beta_{i}\lambda_{i}^{n}
\bigl[q Y^{n}_{i} \bigl\llvert Y^{n}_{i}
\bigr\rrvert^{q-2} +\mathscr{O} \bigl( \bigl\llvert Y^{n}_{i}
\bigr\rrvert^{q-2} \bigr) \bigr]
\nonumber\\
&&{}- \sum_{i=1}^{d} \beta_{i}
\mu_{i}^{n} x_{i} \bigl[q Y^{n}_{i}
\bigl\llvert Y^{n}_{i} \bigr\rrvert^{q-2} +
\mathscr{O} \bigl( \bigl\llvert Y^{n}_{i} \bigr
\rrvert^{q-2} \bigr) \bigr]
\nonumber
\\
&&{}-\sum_{i=1}^{d} \beta_{i}
\bigl(\gamma^{n}_{i}-\mu ^{n}_{i} \bigr)
Q^{n}_{i}[x] \bigl[q Y^{n}_{i} \bigl
\llvert Y^{n}_{i} \bigr\rrvert^{q-2} +\mathscr{O}
\bigl( \bigl\llvert Y^{n}_{i} \bigr\rrvert^{q-2}
\bigr) \bigr]
\nonumber
\\
\label{uni-10}
& \le & \sum_{i=1}^{d}\beta_{i}
\bigl(\lambda_{i}^{n}+\mu_{i}^{n}
x_{i} + \bigl\llvert\gamma^{n}_{i}-
\mu^{n}_{i} \bigr\rrvert Q_{n}^{i}[x]
\bigr) \mathscr{O} \bigl( \bigl\llvert Y^{n}_{i} \bigr
\rrvert^{q-2} \bigr)
\\
&&{}+\sum_{i=1}^{d} \beta_{i}
qY^{n}_{i} \bigl\llvert Y^{n}_{i} \bigr
\rrvert^{q-2} \bigl( \lambda_{i}^{n}-
\mu_{i}^{n} x_{i} - \bigl(\gamma^{n}_{i}-
\mu_{i}^{n} \bigr)Q^{n}_{i}[x] \bigr)
\nonumber
\\
& \le & \sum_{i=1}^{d}\beta_{i}
\bigl(\lambda_{i}^{n} + \bigl(\mu_{i}^{n}+
\bigl\llvert\gamma^{n}_{i}-\mu^{n}_{i}
\bigr\rrvert \bigr) \bigl(Y^{n}_{i}+\rho_{i} n
\bigr) \bigr) \mathscr{O} \bigl( \bigl\llvert Y^{n}_{i} \bigr
\rrvert^{q-2} \bigr)
\nonumber
\\
&&{}+\sum_{i=1}^{d} \beta_{i}
qY^{n}_{i} \bigl\llvert Y^{n}_{i} \bigr
\rrvert^{q-2} \bigl( \lambda_{i}^{n}-
\mu_{i}^{n} x_{i} - \bigl(\gamma^{n}_{i}-
\mu_{i}^{n} \bigr)Q^{n}_{i}[x] \bigr),
\nonumber
\end{eqnarray}
where in the last inequality we use the fact that $Q_{i}^{n}[x]\le
x_{i}$ for
$x\in\mathbb{Z}^{d}_{+}$.
Let
\[
\delta^{n}_{i} :=\lambda_{i}^{n}-
\mu_{i}^{n}\rho_{i}n = \mathscr{O}(\sqrt{n}).
\]
The last estimate is due to the assumptions in \eqref{HWpara}
concerning the parameters in the Halfin--Whitt regime.
Then
%
\begin{eqnarray}
&&\qquad \sum_{i=1}^{d}
\beta_{i} qY^{n}_{i} \bigl\llvert
Y^{n}_{i} \bigr\rrvert^{q-2} \bigl(
\lambda_{i}^{n}-\mu_{i}^{n}
x_{i}- \bigl(\gamma^{n}_{i}-\mu_{i}^{n}
\bigr)Q^{n}_{i}[x] \bigr)
\nonumber
\\[-8pt]
\label{uni-07}
\\[-8pt]
\nonumber
&& \qquad\qquad= - q \sum_{i=1}^{d} \beta_{i}
\mu_{i}^{n} \bigl\llvert Y^{n}_{i} \bigr
\rrvert^{q} + \sum_{i=1}^{d}
\beta_{i} qY^{n}_{i} \bigl\llvert
Y^{n}_{i} \bigr\rrvert^{q-2} \bigl(
\delta_{i}^{n}- \bigl(\gamma^{n}_{i}-
\mu_{i}^{n} \bigr)Q^{n}_{i}[x] \bigr).
\end{eqnarray}
If $x\in A_{n}$ and $n$ is large, then
\begin{eqnarray*}
Q_{i}^{n}[x] &=& u_{h}(x)  = \varpi \bigl((e
\cdot x-n)^{+}v_{\delta}(\hat{x}_{n}) \bigr)
\\
& \le & (e\cdot x-n)^{+} +d \le2dK\sqrt{n}.
\end{eqnarray*}
Let $x\in A^{c}_{n}$.
We use the fact that for any $a,b\in\mathbb{R}$ it holds that
$a^{+}-b^{+}=\xi[a-b]$ for some $\xi\in[0,1]$.
Also,
\[
\Biggl[n \rho_{i}- \Biggl(n-\sum_{j=1}^{i-1}
n \rho_j \Biggr)^{+} \Biggr]^{+} = 0,\qquad i=1,
\ldots,d.
\]
Thus, we obtain maps $\xi, \tilde{\xi}\dvtx\mathbb{R}^{d}\to
[0,1]^{d}$ such that
\begin{eqnarray*}
-Q_{i}^{n}[x]& = &\Biggl[n \rho_{i}- \Biggl(n-
\sum_{j=1}^{i-1} n \rho_j
\Biggr)^{+} \Biggr]^{+} -Q_{i}^{n}[x]
\\
& =& \xi_{i}(x) (n\rho_{i}-x_{i})-\tilde{
\xi}_{i}(x) \sum_{j=1}^{i-1}(x_j-n
\rho_j),\qquad x\in A^{c}_{n}.
\end{eqnarray*}
Hence, from \eqref{uni-07} we obtain
\begin{eqnarray*}
&& \sum_{i=1}^{d} \beta_{i}
qY^{n}_{i} \bigl\llvert Y^{n}_{i} \bigr
\rrvert^{q-2} \bigl(\lambda_{i}^{n}-
\mu_{i}^{n} x_{i}- \bigl(\gamma^{n}_{i}-
\mu_{i}^{n} \bigr)Q^{n}_{i}[x] \bigr)
\\
&&\qquad\le \mathscr{O}(\sqrt{n}) q\sum_{i=1}^{d}
\beta_{i} \bigl\llvert Y^{n}_{i} \bigr
\rrvert^{q-1}
- q \sum_{i=1}^{d} \beta_{i}
\bigl( \bigl(1-\xi_{i}(x) \bigr)\mu_{i}^{n} +
\xi_{i}(x)\gamma^{n}_{i} \bigr) \bigl\llvert
Y^{n}_{i} \bigr\rrvert^{q}
\\
&&\qquad\quad{}+ q\sum_{i=1}^{d} \beta_{i}
Y^{n}_{i} \bigl\llvert Y^{n}_{i} \bigr
\rrvert^{q-2} \Biggl(\delta_{i}^{n}- \bigl(
\gamma^{n}_{i}-\mu_{i}^{n} \bigr)
\tilde{ \xi}_{i}(x)\sum_{j=1}^{i-1}Y^{n}_{j}
\Biggr),
\end{eqnarray*}
where we used the fact that on $A_{n}$ we have
\[
\Biggl[ x_{i}- \Biggl(n-\sum_{j=1}^{i-1}
x_j \Biggr)^{+} \Biggr]^{+} = \mathscr{O}(
\sqrt{n})\qquad \forall i.
\]
Observe that there exists $\vartheta>0$,
independent of $n$ due to \eqref{HWpara}, such that
\[
\bigl(1-\xi_{i}(x) \bigr)\mu_{i}^{n}+
\xi_{i}(x)\gamma^{n}_{i} \ge\min \bigl(
\mu_{i}^{n}, \gamma^{n}_{i} \bigr) \ge
\vartheta
\]
for all $n\in\mathbb{N}$, all $x\in\mathbb{R}^{d}$, and all
$i=1,\ldots,d$.
As a result, we obtain
%
\begin{eqnarray}
&&\sum_{i=1}^{d}
\beta_{i} qY^{n}_{i} \bigl\llvert
Y^{n}_{i} \bigr\rrvert^{q-2} \bigl(
\lambda_{i}^{n}-\mu_{i}^{n}
x_{i}- \bigl(\gamma^{n}_{i}-\mu_{i}^{n}
\bigr)Q^{n}_{i}[x] \bigr)
\nonumber\\
\label{uni-08}
&& \qquad\le\mathscr{O}(\sqrt{n})q\sum_{i=1}^{d}
\beta_{i} \bigl\llvert Y^{n}_{i} \bigr
\rrvert^{q-1} - q \vartheta\sum_{i=1}^{d}
\beta_{i} \bigl\llvert Y^{n}_{i} \bigr
\rrvert^{q}
\\
&&\qquad\quad{}+ q\sum_{i=1}^{d} \beta_{i}
Y^{n}_{i} \bigl\llvert Y^{n}_{i} \bigr
\rrvert^{q-2} \Biggl(\delta_{i}^{n}- \bigl(
\gamma^{n}_{i}-\mu_{i}^{n} \bigr)
\tilde{ \xi}_{i}(x)\sum_{j=1}^{i-1}Y^j_{n}
\Biggr).\nonumber
\end{eqnarray}
We next estimate the last term on the right-hand side of \eqref{uni-08}.
Let
$\kappa:=\sup_{n,i} |\gamma^{n}_{i}-\mu_{i}^{n}|$, and
$\varepsilon_{1}:=\frac{\vartheta}{8\kappa}$.
Using Young's inequality, we obtain the estimate
\[
\bigl\llvert Y^{n}_{i} \bigr\rrvert^{q-1} \Biggl|\sum
_{j=1}^{i-1}Y^j_{n} \Biggr|
\le \varepsilon_{1} \bigl\llvert Y^{n}_{i} \bigr
\rrvert^{q} + \frac{1}{\varepsilon_{1}^{q-1}} \Biggl|\sum_{j=1}^{i-1}Y^{n}_{j}
\Biggr|^{q}.
\]
Therefore,
\begin{eqnarray*}
&& q\sum_{i=1}^{d} \beta_{i}
Y^{n}_{i} \bigl\llvert Y^{n}_{i} \bigr
\rrvert^{q-2} \Biggl(- \bigl(\gamma^{n}_{i}-
\mu_{i}^{n} \bigr) \tilde{\xi}_{i}(x)\sum
_{j=1}^{i-1}Y^{n}_{j} \Biggr)
\\
&&\qquad \le q\kappa\sum_{i=1}^{d} \Biggl(
\varepsilon_{1}\beta_{i} \bigl\llvert Y^{n}_{i}
\bigr\rrvert^{q} +\frac{\beta_{i}}{\varepsilon_{1}^{q-1}} \Biggl|\sum
_{j=1}^{i-1}Y^{n}_{j}
\Biggr|^{q} \Biggr)
\\
&&\qquad \le q\kappa\sum_{i=1}^{d} \Biggl(
\varepsilon_{1}\beta_{i} \bigl\llvert Y^{n}_{i}
\bigr\rrvert^{q}+ \frac{\beta_{i}}{\varepsilon_{1}^{q-1}}d^{q-1}\sum
_{j=1}^{i-1} \bigl\llvert Y^{n}_{j}
\bigr\rrvert^{q} \Biggr)
\\
&&\qquad = \frac{q \vartheta}{8}\sum_{i=1}^{d}
\Biggl(\beta_{i} \bigl\llvert Y^{n}_{i} \bigr
\rrvert^{q}+ \frac{\beta_{i}}{\varepsilon_{1}^{q}}d^{q-1}\sum
_{j=1}^{i-1} \bigl\llvert Y^{n}_{j}
\bigr\rrvert^{q} \Biggr).
\end{eqnarray*}
We choose $\beta_{1}=1$ and for $i\ge2$, we define $\beta_{i}$ by
\[
\beta_{i} := \frac{\varepsilon_{1}^{q}}{d^{q}}\min_{j\le i-1}
\beta_{j}.
\]
With this choice of $\beta_{i}$ it follows from above that
\[
q\sum_{i=1}^{d} \beta_{i}
Y^{n}_{i} \bigl\llvert Y^{n}_{i} \bigr
\rrvert^{q-2} \Biggl(- \bigl(\gamma^{n}_{i}-
\mu_{i}^{n} \bigr) \tilde{\xi}_{i}(x)\sum
_{j=1}^{i-1}Y^j_{n} \Biggr)
\le \frac{q \vartheta}{4}\sum_{i=1}^{d}
\beta_{i} \bigl\llvert Y^{n}_{i} \bigr
\rrvert^{q}.
\]
Using the preceding inequality in \eqref{uni-08}, we obtain
%
\begin{eqnarray}
&& \sum_{i=1}^{d}
\beta_{i} qY^{n}_{i} \bigl\llvert
Y^{n}_{i} \bigr\rrvert^{q-2} \bigl(
\lambda_{i}^{n}-\mu_{i}^{n}
x_{i}- \bigl(\gamma^{n}_{i}-\mu_{i}^{n}
\bigr)Q^{n}_{i}[x] \bigr)
\nonumber
\\[-8pt]
\label{uni-09}
\\[-8pt]
\nonumber
&& \qquad\le\mathscr{O}(\sqrt{n})q\sum_{i=1}^{d}
\beta_{i} \bigl\llvert Y^{n}_{i} \bigr
\rrvert^{q-1} -\frac{3}{4}q \vartheta\sum
_{i=1}^{d} \beta_{i} \bigl\llvert
Y^{n}_{i} \bigr\rrvert^{q}.
\end{eqnarray}
Combining \eqref{uni-10} and \eqref{uni-09}, we obtain
%
\begin{eqnarray}
 \mathcal{L}_{n} f_{n}(x) &\le& \sum
_{i=1}^{d} \mathscr{O}(\sqrt{n})\mathscr{O} \bigl(
\bigl\llvert Y^{n}_{i} \bigr\rrvert ^{q-1} \bigr)+
\sum_{i=1}^{d} \mathscr{O}(n)\mathscr{O}
\bigl( \bigl\llvert Y^{n}_{i} \bigr\rrvert^{q-2}
\bigr)
\nonumber
\\[-8pt]
\label{uni-11}
\\[-8pt]
\nonumber
&& {}-\frac{3}{4}q \vartheta\sum_{i=1}^{d}
\beta_{i} \bigl\llvert Y^{n}_{i} \bigr
\rrvert^{q}.
\end{eqnarray}
By Young's inequality, for any $\varepsilon>0$, we have the bounds
\begin{eqnarray*}
\mathscr{O}(\sqrt{n})\mathscr{O} \bigl( \bigl\llvert Y^{n}_{i}
\bigr\rrvert ^{q-1} \bigr)& \le & \varepsilon \bigl[\mathscr{O} \bigl(
\bigl\llvert Y^{n}_{i} \bigr\rrvert ^{q-1} \bigr)
\bigr]^{{q}/({q-1})} + \varepsilon^{(1-q)} \bigl[\mathscr{O}(\sqrt{n})
\bigr]^{q},
\\
\mathscr{O}(n)\mathscr{O} \bigl( \bigl\llvert Y^{n}_{i}
\bigr\rrvert^{q-2} \bigr) & \le &\varepsilon \bigl[\mathscr{O} \bigl(
\bigl\llvert Y^{n}_{i} \bigr\rrvert ^{q-2} \bigr)
\bigr]^{{q}/({q-2})}+ \varepsilon^{(1-{q}/{2})} \bigl[\mathscr{O}(n)
\bigr]^{{q}/{2}}.
\end{eqnarray*}
Thus, choosing $\varepsilon$ properly in \eqref{uni-11} we obtain
\eqref{uni-06}.

We proceed to complete the proof of the lemma by applying \eqref{uni-06}.
First, we observe that
$\mathbb{E} [\sup_{s\in[0,T]} |X^{n}(s)|^{p} ]$
is finite for any $p\ge1$ as this quantity
is dominated
by the Poisson arrival process.
Therefore, from \eqref{uni-06} we see that
\begin{eqnarray*}
\mathbb{E} \bigl[f_{n} \bigl(X^{n}(T) \bigr)
\bigr]-f_{n} \bigl(X^{n}(0) \bigr) & =& \mathbb{E} \biggl[\int
_{0}^{T}\mathcal{L}_{n}f_{n}
\bigl(X^{n}(s) \bigr)\, \mathrm{d} {s} \biggr]
\\
& \le& C_{1} n^{{q}/{2}}T -C_{2} \mathbb{E} \biggl[
\int_{0}^{T}f_{n} \bigl(X^{n}(s)
\bigr) \,\mathrm{d} {s} \biggr],
\end{eqnarray*}
which implies that
\[
C_{2}\mathbb{E} \Biggl[\int_{0}^{T}
\sum_{i=1}^{d} \beta_{i} \bigl(
\hat{X}_{i}^{n}(s) \bigr)^{q} \,\mathrm{d} {s}
\Biggr] \le C_{1} T+\sum_{i=1}^{d}
\beta_{i} \bigl(\hat{X}_{i}^{n}(0)
\bigr)^{q}.
\]
Hence, the proof follows by dividing both sides by $T$ and letting
$T\to\infty$.
\end{pf}

\begin{pf*}{Proof of Theorem~\ref{T-upperbound}}
Let $r$ be the given running cost with polynomial growth with exponent
$m$ in \eqref{cost1}.
Let $q=2(m+1)$.
Recall that $\tilde{r}(x,u)=r((e\cdot x)^{+}u)$ for $(x,u)\in\mathbb
{R}^{d}\times\mathcal{S}$.
Then $\tilde{r}$ is convex in $u$ and satisfies \eqref{eg-cost} with the
same exponent $m$.
For any $\delta>0$, we choose
$v_{\delta}\in\mathfrak{U}_{\mathrm{SSM}}$ such that $v_{\delta}$
is a continuous
precise control with invariant probability measure $\mu_{\delta}$ and
%
\begin{equation}
\label{300} \int_{\mathbb{R}^{d}}\tilde{r} \bigl(x,v_{\delta}(x)
\bigr) \mu_{\delta
}(\mathrm{d} {x}) \le \varrho_{*} +\delta.
\end{equation}
We also want the control $v_{\delta}$ to have the property that
$v_{\delta}(x)=(0,\ldots,0,1)$
outside a large ball. To obtain such $v_{\delta}$, we see that by
Theorems~\ref{T-trunc},
\ref{T-trunc2} and Remark~\ref{rem-trunc} we can find $v'_{\delta}$
and a
ball $B_{l}$ for $l$ large, such that $v'_{\delta}\in\mathfrak
{U}_{\mathrm{SSM}}$,
$v'_{\delta}(x)=e_{d}$ for $\llvert x\rrvert>l$, $v'_{\delta}$ is
continuous in $B_{l}$,
and
\[
\biggl|\int_{\mathbb{R}^{d}}\tilde{r} \bigl(x,v'_{\delta}(x)
\bigr) \mu '_{\delta}(\mathrm{d} {x}) -\varrho_{*}
\biggr| < \frac{\delta}{2},
\]
where $\mu'_{\delta}$ is the invariant probability measure corresponding
to $v'_{\delta}$.
We note that $v'_{\delta}$ might not be continuous on $\partial B_{l}$.
Let $\{\chi^{n} \dvtx n\in\mathbb{N}\}$ be a sequence
of cut-off functions such that $\chi^{n}\in[0,1]$, it vanishes on
$B^{c}_{l-({1}/{n})}$, and it takes the value $1$ on $B_{l-({2}/{n})}$.
Define the sequence
$v^{n}_{\delta}(x):=\chi^{n}(x)v'_{\delta}(x)+(1-\chi^{n}(x))e_{d}$.
Then $v^{n}_{\delta}\to v'_{\delta}$,
as $n\to\infty$, and the convergence is uniform on the
complement of any neighborhood of $\partial B_{l}$.
Also by Proposition~\ref{eg-prop} the corresponding invariant
probability measures
$\mu^{n}_{\delta}$ are exponentially tight.
Thus,
\[
\biggl\llvert \int_{\mathbb{R}^{d}}\tilde{r} \bigl(x,v'_{\delta}(x)
\bigr) \mu '_{\delta}(\mathrm{d} {x}) -\int
_{\mathbb{R}^{d}}\tilde{r} \bigl(x,v^{n}_{\delta}(x)
\bigr) \mu ^{n}_{\delta}(\mathrm{d} {x}) \biggr\rrvert \mathop{
\longrightarrow}_{n\to\infty} 0.
\]
Combining the above two expressions, we can easily find $v_{\delta}$
which satisfies \eqref{300}.
We construct a scheduling policy as in Lemma~\ref{lem-uni}.
By Lemma~\ref{lem-uni}, we see that
for some constant $K_{1}$ it holds that
%
\begin{equation}
\label{uni-04} \sup_{n\ge n_{0}} \limsup_{T\to\infty}
\frac{1}{T} \mathbb{E} \biggl[\int_{0}^{T}
\bigl\llvert\hat{X}^{n}(s) \bigr\rrvert^{q} \,\mathrm{d} {s}
\biggr] < K_{1},\qquad  q=2(m+1).
\end{equation}
Define
\begin{eqnarray*}
v_{h}(x) & :=&\varpi \bigl((e\cdot x-n)^{+}v_{\delta}(
\hat {x}_{n}) \bigr),
\\
\hat{v}_{h} \bigl(\hat{x}^{n} \bigr) & :=& \varpi \bigl(
\sqrt{n} \bigl(e\cdot\hat{x}^{n} \bigr)^{+} v_{\delta}
\bigl(\hat{x}^{n} \bigr) \bigr).
\end{eqnarray*}
Since $v_{\delta}(\hat{x}^{n})=(0,\ldots,0,1)$ when
$|\hat{x}^{n}|\ge K$, it follows that
\[
Q^{n} \bigl[X^{n} \bigr] = X^{n}-Z^{n}
\bigl[X^{n} \bigr] = v_{h} \bigl(X^{n} \bigr)
\]
for large $n$, provided that $\sum_{i=1}^{d-1}X_{i}^{n}\le n$.
Define
\[
D_{n}:= \Biggl\{x \dvtx \sum_{i=1}^{d-1}
\hat{x}_{i}^{n} > \rho_{d} \sqrt{n} \Biggr\}.
\]
Then
\begin{eqnarray*}
r \bigl(\hat{Q}^{n}(t) \bigr) &=& r \biggl(\frac{1}{\sqrt{n}}
\hat{v}_{h} \bigl(\hat{X}^{n}(t) \bigr) \biggr)+ r \bigl(
\hat{X}^{n}(t)-\hat{Z}^{n}(t) \bigr) \mathbb{I}_{\{\hat
{X}^{n}(t)\in D_{n}\}}
\\
&&{}-r \biggl(\frac{1}{\sqrt{n}}\hat{v}_{h} \bigl(\hat{X}^{n}(t)
\bigr)\biggr) \mathbb{I}_{\{\hat{X}^{n}(t)\in D_{n}\}}.
\end{eqnarray*}
Define, for each $n$, the mean empirical measure $\Psi^{n}_{T}$ by
\[
\Psi^{n}_{T}(A) :=\frac{1}{T} \mathbb{E} \biggl[\int
_{0}^{T} \mathbb{I}_A \bigl(
\hat{X}^{n}(t) \bigr) \,\mathrm{d} {t} \biggr].
\]
By \eqref{uni-04}, the family $\{\Psi^{n}_{T} \dvtx T>0, n\ge1\}$
is tight.
We next show that
%
\begin{equation}
\label{uni-05}\quad \lim_{n\to\infty} \limsup_{T\to\infty}
\frac{1}{T} \mathbb{E} \biggl[\int_{0}^{T} r
\bigl(\hat{Q}^{n}(t) \bigr) \,\mathrm{d} {t} \biggr] = \int
_{\mathbb{R}^{d}}r \bigl((e\cdot x)^{+}v_{\delta}(x)
\bigr) \mu_{\delta}(\mathrm{d} {x}).
\end{equation}
For each $n$, select a sequence
$\{T^{n}_{k} \dvtx k\in\mathbb{N}\}$ along which the ``$\limsup$''
in \eqref{uni-05} is attained.
By tightness, there exists a limit point $\Psi^{n}$
of $\Psi^{n}_{T^{n}_{k}}$.
Since $\Psi^{n}$ has support on a discrete lattice, we have
\[
\int_{\mathbb{R}^{d}} r \biggl(\frac{1}{\sqrt{n}}\hat
{v}_{h}(x) \biggr) \Psi^{n}_{T^{n}_{k}}(\mathrm{d} {x})
\mathop{\longrightarrow}_{k\to\infty} \int_{\mathbb{R}^{d}} r \biggl(
\frac
{1}{\sqrt{n}} \hat{v}_{h}(x) \biggr)\Psi^{n}(\mathrm{d}
{x}).
\]
Therefore,
\[
\limsup_{T\to\infty} \frac{1}{T} \mathbb{E} \biggl[\int
_{0}^{T} r \bigl(\hat{Q}^{n}(t) \bigr)
\,\mathrm{d} {t} \biggr] \lessgtr\int_{\mathbb{R}^{d}} r \biggl(
\frac{1}{\sqrt{n}}\hat{v}_{h}(x) \biggr) \Psi^{n}(\mathrm{d} {x})\pm\mathcal{E}^{n},
\]
where
\[
\mathcal{E}^{n} = \limsup_{T\to\infty} \frac{1}{T}
\mathbb {E} \biggl[\int_{0}^{T} \biggl(r \bigl(
\hat{Q}^{n}(t) \bigr)+ r \biggl(\frac{1}{\sqrt{n}}\hat{v}_{h}
\bigl(\hat{X}^{n}(t) \bigr) \biggr)\biggr) \mathbb{I}_{\{\hat{X}^{n}(t)\in D_{n}\}}
\,\mathrm{d} {t} \biggr].
\]
By \eqref{uni-04}, the family $\{\Psi^{n} \dvtx n\ge1\}$ is tight.
Hence, it has a limit $\Psi$.
By definition, we have
\[
\biggl\llvert \frac{1}{\sqrt{n}}\hat{v}_{h}(x) -(e\cdot
x)^{+}v_{\delta
}(x) \biggr\rrvert \le\frac{2d}{\sqrt{n}}.
\]
Thus, using the continuity property of $r$ and \eqref{cost1}
it follows that
\[
\int_{\mathbb{R}^{d}} r \biggl(\frac{1}{\sqrt{n}}\hat{v}_{h}(x)
\biggr) \Psi^{n}(\mathrm{d} {x}) \mathop{\longrightarrow}_{n\to\infty}
\int_{\mathbb{R}^{d}} r \bigl((e\cdot x)^{+}v_{\delta}(x)
\bigr) \Psi(\mathrm{d} {x}),
\]
along some subsequence.
Therefore, in order to complete the proof of \eqref{uni-05} we need to
show that
\[
\limsup_{n\to\infty} \mathcal{E}^{n} = 0.
\]
Since the policies are work-conserving, we observe that
$0\le\hat{X}^{n}-\hat{Z}^{n}\le(e\cdot\hat{X}^{n})^{+}$,
and therefore for some positive constants $\kappa_{1}$ and $\kappa
_{2}$, we have
\[
r \biggl(\frac{1}{\sqrt{n}}\hat{v}_{h} \bigl(\hat{X}^{n}(t)\bigr)
\biggr)\vee r \bigl(\hat{X}^{n}(t)-\hat{Z}^{n}(t) \bigr) \le
\kappa_{1} + \kappa_{2} \bigl[ \bigl(e\cdot
\hat{X}^{n} \bigr)^{+} \bigr]^{m}.
\]
Given $\varepsilon>0$ we can choose $n_{1}$ so that for all $n\ge n_{1}$,
\[
\limsup_{T\to\infty} \frac{1}{T} \mathbb{E} \biggl[\int
_{0}^{T} \bigl[ \bigl(e\cdot\hat
{X}^{n}(s) \bigr)^{+} \bigr]^{m}
\mathbb{I}_{ \{\llvert\hat{X}^{n}(s)\rrvert>({\rho_{d}}/{\sqrt{d}})\sqrt{n} \}} \,\mathrm{d} {s} \biggr] \le\varepsilon,
\]
where we use \eqref{uni-04}.
We observe that
$D_{n}\subset \{|\hat{x}^{n}|>\rho_{d}\sqrt{{d}/{n}} \} $.
Thus, \eqref{uni-05} holds.
In order to complete the proof, we only need to
show that $\Psi$ is the invariant probability measure corresponding
to $v_{\delta}$.
This can be shown using the convergence of generators as in the proof of
Theorem~\ref{T-lowerbound}.
\end{pf*}

\section{Conclusion}

We have answered some of the most interesting questions for the ergodic
control problem of the Markovian multi-class many-server queueing model.
This current study has raised some more questions for future research.
One of the interesting questions is to consider nonpreemptive policies and
try to establish asymptotic optimality in the class of nonpreemptive
admissible polices~\cite{atar-mandel-rei}.
It will also be interesting to study a similar control problem when
the system has multiple heterogeneous agent pools with skill-based routing.

It has been observed that customers' service requirements and patience times
are nonexponential \cite{brown-et-al}
in some situations.
It is therefore important and interesting to address similar control problems
under general
assumptions on the service and patience time distributions.

\section*{Acknowledgements}
We thank the anonymous referee for many helpful comments that
have led to significant improvements in our paper.
Ari Arapostathis acknowledges the hospitality the Department of Industrial
and Manufacturing Engineering in Penn State while he was visiting at the
early stages of this work.
Guodong Pang acknowledges the hospitality of the Department of Electrical
and Computer Engineering at University of Texas at Austin while he was visiting
for this work.
Part of this work was done while Anup Biswas was
visiting the Department of Industrial
and Manufacturing Engineering in Penn State.
Hospitality of the department is acknowledged.

%





\printaddresses
\end{document}